\documentclass[10pt,a4paper,english]{amsart}

\usepackage{amssymb}
\usepackage{amsthm}
\usepackage[T1]{fontenc}
\usepackage{lmodern}
\usepackage{ae}
\usepackage{enumerate}
\usepackage{mathtools}

\usepackage[utf8]{inputenc}
\usepackage{hyperref}
\hypersetup{
   colorlinks=true,
   linkcolor=blue,
   citecolor=blue,
   filecolor=red,
   urlcolor=green
}
\usepackage{geometry}
\usepackage{fancyhdr}
\usepackage{stmaryrd}
\usepackage{amsmath}
\theoremstyle{plain} \newtheorem{theorem}{Theorem}[section]

\newtheorem{definition}[theorem]{Definition}
\newtheorem{proposition}[theorem]{Proposition}
\newtheorem{lemma}[theorem]{Lemma}
\newtheorem{corollary}[theorem]{Corollary}
\newtheorem{remark}[theorem]{Remark}
\numberwithin{equation}{section}
\newcommand{\diag}{\operatorname{diag}}

\renewcommand{\O}{  \mathcal{O}   }

\newcommand{\s}{  \sigma   }

\renewcommand{\phi}{  \varphi  }

\newcommand{\Card}{\operatorname{Card}}
\newcommand{\be}{\begin{equation}}
\newcommand{\ee}{\end{equation}}
\newcommand{\ben}{\begin{equation*}}
\newcommand{\een}{\end{equation*}}
\newcommand{\ban}{\begin{align*}}
\newcommand{\ean}{\end{align*}}
\newcommand{\lc}{\llbracket}
\newcommand{\rc}{\rrbracket}
\author{Moudhaffar Bouthelja}
\address{Laboratoire Paul-Painlev\'e, Universit\'e de Lille 1, UMR CNRS 8524, Cit\'e Scientifique, 59655 Villeneuve-d'Ascq}
\email{moudhaffar.bouthelja@math.univ-lille1.fr}
\thanks{This work was supported in part by the CPER Photonics4Society and the Labex CEMPI (ANR-11-LABX-0007-01)}
\title[\small{KAM for the NLW equation on the circle: Small amplitude solution}]{KAM for the nonlinear wave equation on the circle:\\ Small amplitude solution}
\begin{document}
	\begin{abstract}
	In this paper we consider the nonlinear wave equation on the circle:
\begin{equation} \nonumber
u_{tt} - u_{xx} + m u = g(x,u), \quad t \in \mathbb{R},\: x \in \mathbb{S}^1,
\end{equation}
where $m \in [1,2]$ is a mass and $g(x,u)=4u^3+ O(u^4)$. This equation will be treated as a perturbation of the integrable Hamiltonian:
\begin{equation} \tag{$\ast$} \label{first equation}
u_t= v, \quad v_t = - u_{xx} + m u.
\end{equation}
Near the origin and for generic $m$, we prove the existence of small amplitude quasi-periodic  solutions close to the solution of the linear equation
\eqref{first equation}.  For the proof we use an abstract KAM theorem in infinite dimension and a Birkhoff normal form result.
	\end{abstract}
	\maketitle
	\tableofcontents

\section{Introduction and results}
\subsection{Introduction}
We consider the cubic wave equation on the circle:
\begin{equation} \label{eq onde}
u_{tt} - u_{xx} + m u = g(x,u), \quad t \in \mathbb{R},\: x \in \mathbb{S}^1,
\end{equation}
where $m\in [1,2]$ is a mass and $g$ is a real holomorphic function on $\mathbb{S}^1\times J$, for $J$ some neighborhood of the origin of $\mathbb{R}$. We suppose that the nonlinearity $g$ satisfies
\begin{equation} \label{non-linearite}
g(x,u)=4u^3+ O(u^4).
\end{equation}
We prove the existence of small amplitude quasi-periodic  solutions close to the solution of the linear equation.

Since the space variable belongs to the circle, we can diagonalize the linear part of the equation in Fourier basis. So we can study the PDE as a perturbation of an integrable Hamiltonian of the following form
\begin{equation}  \tag{$\square$} \label{hamiltonien-élémentaire}
H= \sum_{s \in \mathbb{Z}} \lambda_s \xi_s \eta_s + \mbox{ Perturbation },
\end{equation}
where $\lambda_s=\sqrt{s^2+m}$. In order to prove the existence of quasi-periodic solutions, we will use an abstract KAM theorem in infinite dimension adapted to our situation and proven in \cite{B1}. The KAM theory (Kolmogorov-Arnold-Moser) tells us that, under the effect of a small perturbation and under several conditions of non resonance, an integrable Hamiltonian system continues to exhibit finite-dimensional invariant tori in an infinite dimensional space. The existence of these invariant tori gives us the existence of quasi-periodic solutions. The main issue here is that frequencies $\lambda_s$ do not satisfy the standard non resonance hypotheses \footnote{Hypothese likes: $\vert k_1 \lambda_1 + \ldots + k_n \lambda_n \vert  \geq \frac{\gamma}{\vert k \vert^\tau}, \:\: k \in \mathbb{Z}^n \setminus\lbrace 0 \rbrace.$}. In the Dirichlet case, the sum in \eqref{hamiltonien-élémentaire} is restricted to positive indices (see \cite{Wayne90}). In this case, the standard non resonance conditions can be verified with the mass $m$. In the periodic case, both positive and negative indices are allowed. Note that $\lambda_s=\lambda_{-s}$. So we obtain a resonant Hamiltonian system. For this purpose, the KAM theorem that we will use must deal with the case of multiple eigenvalues.

The existence of quasi-periodic solutions for nonlinear Hamiltonian PDEs have interested many authors. The first result related to preserving such solutions, after the perturbation of an integrable Hamiltonian of infinite dimension, was given by Kuksin in 1987 in \cite{Kuk87,kuk89} for the Schr\"odinger equation in dimension 1 with Dirichlet conditions.

Concerning the wave equation, the first result is due to Wayne in \cite{Wayne90}. He considered the cubic-wave equation in dimension 1 with external potential in $ L ^ 2 ([0,1]) $, and with Dirichlet conditions (which leads to simplicity of the spectrum).

We can also cite the work of P{\"o}schel in \cite{Poschel_Wave}. In this paper, the author considers the wave equation in dimension 1 with mass, homogeneous Dirichlet condition, and analytical cubic nonlinearity that does not depends on the space variable.

In 1998; Chierchia and You consider in \cite{Chierchia_You} the wave equation in dimension 1 with analytic periodic potential and an analytic quadratic perturbation that does not depends on the space variable. In this case, the potential acts as an external parameter. This makes verifying the non resonance conditions possible. In particular, the authors do not authorize the case of a vanishing potential.

The most recent work is due to Berti, Biasco and Procesi in 2013 in \cite{BBM2013}. In this paper, they consider the derivative wave equation given by:
\begin{equation} \nonumber
u_{tt}-u_{xx}+mu+f(Du)=0, \quad m>0, \quad D:= \sqrt{-\partial_{xx}^2+m}, \quad (t,x) \in \mathbb{R} \times \mathbb{T},
\end{equation}
where $f(s)$ is a real analytic nonlinearity of the form
\begin{equation} \nonumber
f(s)= as^3+\sum_{k\geq5}f_ks^k, \quad a \neq 0.
\end{equation}
We remark that the nonlinearity is independent of the space variable $x$. This implies that the moment $-i \int_\mathbb{T} \bar{u} \partial_x u dx$ is preserved. This symmetry simplifies the proof of the KAM theorem.

In our case, there are no external parameters. The space variable belongs to the circle, so we are in the periodic case. The non-linearity $g$ depends on the space variable.

The plan of the paper is the following:
\begin{itemize}
\item In the first section, we give the main result of the paper (see Theorem \ref{théorème onde}).
\item In the second section, we show that, for an admissible set (see definition \ref{admissible}), the small divisors of the wave equation \eqref{eq onde} admit a positive lower bound. This is proven for $m \in [1,2] \setminus \mathcal{U}$ where $ \mathcal{U} $ is zero Lebesgue measure set.\item In the third section, using a Birkhoff normal form, we transform the resonant Hamiltonian associated to the equation \eqref{eq onde} into a Hamiltonian that satisfies the hypotheses of the KAM theorem (see Theorem~\ref{theo-F.N}).
\item In the fourth section, we state the KAM theorem (see Theorem~\ref{theoreme kam}), and we verify the non resonance hypotheses (see Lemma~\ref{separation freq}-\ref{2end Melnikov for wave}).
\item In the last part we prove the existence of quasi-periodic solutions of small amplitudes for the equation \eqref{eq onde}.
\end{itemize}
\subsection{Results}
We consider the nonlinear wave equation on the circle \eqref{eq onde} with $g$ in the form \eqref{non-linearite}.
Introducing the change of variable $v=\dot u$, the equation \eqref{eq onde} becomes:
\begin{equation*}
\begin{dcases}
\dot{u}= v,\\
\dot{v} = - \Lambda^2 u+g(x,u),
\end{dcases}
\end{equation*}
where $\Lambda:=(\sqrt{-\partial_{xx}+m})$.
Defining $\psi:=\frac{1}{\sqrt{2}} ( \Lambda^{\frac{1}{2}}u-i\Lambda^{-\frac{1}{2}}v)$, we get the following equation for $\dot \psi$:
\begin{equation} \nonumber
\dot{\psi} = \frac{1}{\sqrt{2}} \left( \Lambda^{1/2} \dot{u} - i \Lambda^{-1/2} \dot{v} \right) .
\end{equation}
We note that $u=\Lambda^{-1/2} \left( \frac{\psi+\bar{\psi}}{\sqrt{2}} \right) $; replacing $\dot{u}$ and $\dot{v}$  by their expressions yields:
\begin{equation} \nonumber
\frac{1}{i} \dot \psi = \Lambda \psi - \frac{1}{\sqrt{2}} \Lambda^{-1/2}g\left( x,\Lambda^{-1/2}\left( \frac{\psi+\bar{\psi}}{\sqrt{2}} \right)\right) .
\end{equation}
Let us endow $L^2(\mathbb{S}^1,\mathbb{C})$ with the classical real symplectic form $-i d \psi \wedge d \bar{\psi}=-du \wedge dv$ and consider the following Hamiltonian:
\begin{equation} \nonumber
H(\psi,\bar{\psi}) = \int_{\mathbb{S}^1} (\Lambda \psi)\bar{\psi} dx + \int_{\mathbb{S}^1} G \left( x,\Lambda^{-1/2} \left( \frac{\psi+\bar{\psi}}{\sqrt{2}} \right)  \right) dx,
\end{equation}
where $ G $ is a primitive of $ g $ with respect to $ u $:
\begin{equation} \nonumber
g=\partial_uG,\quad G(x,u)=u^4+ \mathcal{O}(u^5).
\end{equation}
Then, \eqref{eq onde} becomes a Hamiltonian system:
\begin{equation} \nonumber
\dot{\psi} = i \frac{\partial H}{\partial \bar{\psi}}.
\end{equation}

Consider now the complex Fourier orthonormal basis given by $\lbrace \varphi_s(x)=\frac{e^{isx}}{\sqrt{2\pi}} ,\:s\in \mathbb{Z}\rbrace$. In this base, the operator $\Lambda$ is diagonal, and we have:
\begin{equation} \nonumber
\Lambda \varphi_s=\lambda_s \varphi_s,
\end{equation}
with $\lambda_s=\sqrt{s^2+m}$. Decomposing $\psi$ and $\bar{\psi}$ in this basis yields:
\begin{equation} \nonumber
\psi=\underset{s \in \mathbb{Z}}{ \sum}\xi_s\varphi_s \quad \mbox{and}  \quad \bar{\psi}= \underset{s \in \mathbb{Z}}{ \sum} \eta_s \varphi_{-s}.
\end{equation}
By injecting this decomposition into the expression of $H$, we obtain:
\begin{equation} \label{Hamiltonien_0}
H= \underset{s \in \mathbb{Z}}{ \sum} \lambda_s \xi_s \eta_s + \int_{\mathbb{S}^1} G \left( x , \underset{s \in \mathbb{Z}}{ \sum} \frac{\xi_s\varphi_s+\eta_s \varphi_{-s}}{\sqrt{2\lambda_s}} \right) dx.
\end{equation}
Let $\mathcal{P}_{\mathbb{C}}:= \ell^2(\mathbb{Z},\mathbb{C}) \times \ell^2(\mathbb{Z},\mathbb{C})$ that we endow with the complex symplectic form
$-i{ \sum_{s \in \mathbb{Z}}} d\xi_s \wedge d\eta_s$. We define
\begin{equation} \nonumber
\mathcal{P}_{\mathbb{R}}:= \lbrace (\xi,\eta) \in \mathcal{P}_{\mathbb{C}} | \eta_s=\bar{\xi}_s \rbrace.
\end{equation}
Then, equation \eqref{eq onde} is equivalent to the following Hamiltonian system on $\mathcal{P}_{\mathbb{R}}$:
\begin{equation}  \label{onde2}
\begin{dcases}
	\dot{\xi}_s = i \frac{\partial H}{\partial \eta_s},\\
\dot{\eta}_s = -i \frac{\partial H}{\partial \xi_s},
\end{dcases}
\end{equation}
for $s \in \mathbb{Z}$.

From now, we write $H=H_2+P$, where
\begin{equation} \label{hamiltonien}
H_2=\underset{s \in \mathbb{Z}}{ \sum} \lambda_s \xi_s \eta_s, \quad \mbox{and} \quad  P=\int_{\mathbb{S}^1} G \left( x , \underset{s \in \mathbb{Z}}{ \sum} \frac{\xi_s\varphi_s+\eta_s \varphi_{-s}}{\sqrt{2\lambda_s}} \right)dx.
\end{equation}
\begin{remark} \label{forme_perturb}
Recall that $g(x,u)=4u^3+O(u^4)$ and $g=\partial_u G$, so we can decompose~$P$ into $P=P_4+R_5$ where
\begin{align*}
P_4(\xi,\eta) & = \int_{{\mathbb{S}^1}} u^4 dx= \int_{\mathbb{S}^1} \left( \underset{s \in \mathbb{Z}}{ \sum} \frac{\xi_s\varphi_s+\eta_s \varphi_{-s}}{\sqrt{2\lambda_s}} \right)^4dx,\\
R_5 (\xi,\eta,x) & =P(\xi,\eta,x)-P_4(\xi,\eta)=O(\Vert (\xi,\eta) \Vert^5) .
\end{align*}
In addition, $P_4$  reads
\begin{equation} \nonumber
P_4 = \underset{i,j,k,l \in \mathbb{Z}}{ \sum} C(i,j,k,l) (\xi_i+\eta_{-i})(\xi_j+\eta_{-j})(\xi_k+\eta_{-k})(\xi_l+\eta_{-l}),
\end{equation}
where
\begin{equation*}
C(i,j,k,l):=\int_{\mathbb{S}^1} \varphi_i(x)\varphi_j (x)  \varphi_k(x) \varphi_l(x)dx=
\begin{dcases}
\frac{1}{2\pi} \: \quad \text{if}\quad i+j+k+l=0,\\
0 \quad \quad \text{otherwise}.
\end{dcases}
\end{equation*}
\end{remark}

Let $\mathcal{A}$ be a finite set of $\mathbb{Z}$ of cardinality $n$, and a vector $I=(I_a)_{a\in\mathcal{A}}$ with positive components (i.e. $I_a > 0$ for all $ a \in \mathcal{A}$).
Let $T_I^n$ be the real torus of dimension $n$ defined by
\begin{equation} \nonumber
T_I^n=\left\{
\begin{array}{l c l}
\xi_a=\bar{\eta}_a, \: |\xi_a|^2=I_a \quad  &\text{if}&\: a \in \mathcal{A},\\
\xi_s=\eta_s=0 \quad & \text{if} & \: s \in \mathcal{L}=\mathbb{Z} \setminus  \mathcal{A}.
\end{array}
\right.
\end{equation}

This torus is stable by the Hamiltonian flow when the perturbation $ P $ is zero. We can even give the analytic expression of the solution of the linear equation.

Our purpose in all the following is to prove the persistence of the torus $T_I^n$ when the perturbation $ P $ is no longer zero, while making the crucial assumption that this torus is admissible. A torus is said to be admissible if it is constructed from an admissible  set $\mathcal{A}$.

\begin{definition} \label{admissible}
Let $\mathcal{A}$ be a finite set of $\mathbb{Z}$. $\mathcal{A}$ is admissible if, for all $j \in \mathcal{A} \setminus \lbrace 0 \rbrace$, we have $-j \notin \mathcal{A}\setminus \lbrace 0 \rbrace$ .
\end{definition}

Let us introduce the sets $\mathcal{L}=\mathbb{Z}\setminus \mathcal{A}$ and $\mathcal{A}^-=\left\lbrace j \in \mathcal{L} \: | \: -j \in \mathcal{A} \right\rbrace $.
In a neighborhood of the invariant torus $T_I^n$ in $\mathbb{C}^{2n}$, we define the action-angle variables $ (r_a, \theta_a) _ {\mathcal {A}} $ by:
\begin{equation*}
\left\{
\begin{array}{l l l}
\xi_a=\sqrt{(I_a+r_a)} e^{i\theta_a},\\
\eta_a=\bar{\xi}_a.
\end{array}
\right.
\end{equation*}
For $s \in \mathcal{A}$, we denote by $\omega_s$ (instead of $\lambda_s$) the internal frequencies. In these new variables and notations, the quadratic part $H_2$ of $H$ becomes, up to a constant,
\begin{equation} \nonumber
H_2= \underset{a \in \mathcal{A}}{ \sum} \omega_ar_a +   \underset{s \in \mathcal{L}}{ \sum} \lambda_s \xi_s\eta_s.
\end{equation}
In addition, the perturbation becomes:
\begin{equation} \nonumber
P(r,\theta,\xi,\eta)= \int_{\mathbb{S}^1} G(x,\hat{u}_{I,m}(r,\theta,\xi,\eta))dx,
\end{equation}
with
\begin{align*} \nonumber
\hat{u}_{I,m}(r,\theta,\xi,\eta)& =\underset{a \in \mathcal{A}}{ \sum} \sqrt{(I_a+r_a)} \frac{e^{-i\theta_a}\varphi_a(x)+e^{i\theta_a}\varphi_{-a}(x)}{\sqrt{2}(a^2+m)^{1/4}} \\
& + \underset{s \in \mathcal{L}}{ \sum}  \frac{\xi_s \phi_s(x)+\eta_{-s}\phi_s(x)}{\sqrt{2}(s^2+m)^{1/4}}.
\end{align*}
We set $u_{I,m}(\theta,x)=\hat{u}_{I,m}(0,\theta,0,0)$. Then, for any $I \in \mathbb{R}_+^\mathcal{A}$, $m \in \left[ 1,2 \right] $ and $\theta_0 \in \mathbb{S}^1$, the function $(t,x) \mapsto u_{I,m}(\theta_0+t\omega,x)$ is solution of the linear wave equation. In this case, the torus $T_I^n$ is stable by the Hamiltonian flow.
Our goal is to state a similar result when the perturbation is not zero (in the nonlinear case).
\begin{theorem} \label{théorème onde}
Let $\alpha>1/2$. Assume that $\mathcal{A}$ is an admissible set of cardinality~$n$. Assume also that the perturbation  $g$ is real holomorphic on a neighborhood of $\mathbb{S}^1\times J$ with $J$ some neighborhood of the origin of $\mathbb{R}$ and reads $g(x,u)=4u^3+O(u^4)$. There exists a Borel subset $\mathcal{U} \subset \left[ 1 , 2 \right]$ with zero Lebesgue measure, such that for $m \in \left( \left[ 1,2 \right] \setminus \mathcal{U} \right)$, there exists $\nu_0$ that depends on $\mathcal{A}$, $m$, and the nonlinearity $g$, such that:
\begin{itemize}
\item[]  For $0 < \nu \leq \nu_0$ there exists a Borel set $ \mathcal{D}' \subset \left[ \nu , 2 \nu \right] ^n $ asymptotically of full Lebesgue measure,i.e.
\begin{equation} \nonumber
\operatorname{mes}\left( \left[ \nu , 2 \nu \right] ^n \setminus  \mathcal{D}'  \right) \leq \nu^{n+\gamma},
\end{equation}
with $\gamma >0$ and depending on $n$. For $m \in \left( \left[ 1,2 \right] \setminus \mathcal{U} \right)$ and $I \in \mathcal{D}'$, there exists:\begin{enumerate}
\item a function $u(\theta,x)$ analytic in $\theta$ and of class $H^\alpha$ in $x \in \mathbb{S}^1$ such that:
\begin{equation} \nonumber
\underset{\theta \in \mathbb{R}}{ \sup} \parallel u(\theta,.)-u_{I,m}(\theta,.) \parallel_{H^\alpha} \leq C \nu^{4/5},
\end{equation}
with $C$ an absolute constant.
\item a mapping $\omega': \left( \left[ 1,2 \right] \setminus \mathcal{U} \right) \times \mathcal{D}' \rightarrow \mathbb{R}^n $ verifying:
\begin{equation} \nonumber
\omega' = \omega +  M I  + O(\nu^{3/2}),
\end{equation}
such that for any $m \in \left( \left[ 1,2 \right] \setminus \mathcal{U} \right)$ and $I \in \mathcal{D}'$ the function
\begin{equation}  \nonumber
t\mapsto u(\theta+t\omega',x)
\end{equation}
is solution of the wave equation \eqref{eq onde}. This solution is linear stable.
\end{enumerate}
\end{itemize}
\end{theorem}
The rest of the paper will be devoted to the proof of this theorem.
\section{Small divisors}
\subsection{Non resonance of frequencies}
In this section, we assume that $\mathcal{A}$ is an admissible set as in Definition~\ref{admissible}.\\
We consider the frequency vector
\begin{equation} \nonumber
\omega \equiv \omega(m)=(\omega_a(m))_{a\in \mathcal{A}},
\end{equation}
with $\omega_a(m)=\sqrt{a^2+m}$. The main and only result of this section is the following:
\begin{proposition}
\label{NRom}
Consider an admissible set $\mathcal{A}$ of cardinality $n$ that verifies $ \mathcal{A}  \subset \lbrace a \in \mathbb{Z} \mid |a| \leq N , N \geq 1 \rbrace$. Then, for any $k\in\mathbb{Z}^\mathcal{A}\setminus \lbrace 0 \rbrace$, any $\chi>0$ and any $c\in   \mathbb{R} $, we have
\begin{equation*}
 \operatorname{mes}\ \left\lbrace m\in[1,2]\ \mid \
\left|\sum_{a\in\mathcal{A}} k_a\omega_a(m)+c\right|\leq {\chi}\right\rbrace \leq C \frac{N^{2n^2} \chi^{1/n}}{|k|},
\end{equation*}
with $|k|:=\sum_{a\in\mathcal{A}}|k_a|$ and $C>0$ is a constant that depends only on $n$.
\end{proposition}

The proof uses the same arguments as in Theorem 6.5 of \cite{bambusi2003birkhoff} (see also \cite{bambusi1999nekhoroshev} and~\cite{bambusi2006birkhoff}). For clarity, we recall the main steps of the proof.
\begin{lemma} \label{det}
Assume that $ \mathcal{A} \subset \lbrace a \in \mathbb{Z} \mid |a|\leq N \rbrace$.
For any $p\leq n:=\operatorname{Card}( \mathcal{A})$, consider $a_1,\cdots,a_p \in \mathcal{A}$.
The the following determinant
\begin{equation*}
D:=\left|
\begin{matrix}
\frac{d \omega_{a_1}}{dm} & \frac{d \omega_{a_2}}{dm} & .& .&.&\frac{d \omega_{a_p}}{dm}
\\
\frac{d^2 \omega_{a_1}}{dm^2} & \frac{d^2 \omega_{a_2}}{dm^2} & .& .&.&\frac{d^2 \omega_{a_p}}{dm^2}
\\
.& .& .& .& .&.
\\
.& .& .& .& .&.
\\
\frac{d^p\omega_{a_1}}{dm^p}& \frac{d^p\omega_{a_2}}{dm^p}& .& .&.&\frac{d^p\omega_{a_p}}{dm^p}
\end{matrix}
\right|
\end{equation*}
verifies
\begin{equation} \label{estim-D}
|D|\geq C N^{-2p^2},
\end{equation}
where $C=C(p)>0$ is a constant that depends only on $p$.
\end{lemma}

\proof An explicit computation gives
\begin{equation}
\frac{d^j\omega_i}{dm^j}=\frac{(2j-2)!}{2^{2j-1}(j-1)!}
\frac{(-1)^{j+1}}
{(a_i^2+m)^{j-\frac12}}\ .\label{df}
\end{equation}
Inserting this formula in $D$, by factoring from each l-\textit{th} column $(a_\ell^2+m)^{-1/2}$, and from j-\textit{th} row $\frac{(2j-2)!}{2^{2j-1}(j-1)!}$, the determinant is equal up to a sign to
\begin{eqnarray*}
\left[\prod_{l=1}^{p}\omega_{a_\ell}^{-1}\right]
\left[\prod_{j=1}^{p}\frac{(2j-2)!}{2^{2j-1}(j-1)!} \right]
\times
\left|
\begin{matrix}
1& 1& 1 &. & . & . & 1
\cr
x_{a_1}& x_{a_2}& x_{a_3}&.&.&.&x_{a_p}
\cr
x_{a_1}^2& x_{ a_2}^2& x_{a_3}^2&.&.&.&x_{a_p}^2
\cr
.& .& .& .& .&.&.
\cr
.& .& .& .& .&.&.
\cr
.& .& .& .& .&.&.
\cr
x_{a_1}^{p}& x_{a_2}^{p}& x_{a_3}^{p}&.&.&.&x_{a_p}^{p}
\end{matrix}
\right|,
\end{eqnarray*}
where $x_{a}:=(a^2+m)^{-1}\equiv\omega_{a}^{-2}$.
The above determinant is a Vandermonde determinant. It is equal to
\begin{equation*}
\prod_{1\leq l<k\leq p}(x_{a_\ell}-x_{a_k})=\prod_{1\leq l<k\leq p}
\frac{a_k^2-a_\ell^2}{\omega_{a_\ell}^2\omega_{a_k}^2} .
\end{equation*}
The set $\mathcal{A}$ is admissible, so $\left|\prod_{1\leq l<k\leq p}(a_k^2-a_\ell^2)\right|\geq 1$. Since $|a|\leq N$, then for any $a\in\mathcal{A}$ we have $|\omega_a|\leq 2N$. Therefore:
\begin{equation} \nonumber
\prod_{l=1}^{p}\omega_{a_\ell}^{-1}\prod_{1\leq l<k\leq p}
\frac{1}{\omega_{a_\ell}^2\omega_{a_k}^2}\geq \frac{1}{2^{p(2p-1)}} \frac{1}{N^{p(2p-1)}} \geq \frac{1}{2^{p(2p-1)}} \frac{1}{N^{2p^2}},
\end{equation}
which leads to  \eqref{estim-D}.
\endproof

We need the following proposition, presented in appendix B of \cite{ben85}.
\begin{lemma}\label{m1.1}
Let $\left( u^{(1)}, \ldots ,u^{(p)}\right) $ be $p$ independent vector in $\mathbb{R}^n$ such that $ \Vert u^{(j)} \Vert_{\ell^1} \leq K$ for $j \in [1,\ldots,p]$. Let $w$ be a linear combination of $u^{(1)}, \ldots ,u^{(p)}$. There exists  $j\in[1,\ldots ,p]$ such that:
\begin{equation} \label{estim vol}
|u^{(j)}\cdot w|\geq\frac{\Vert w \Vert_{\ell^2} \mathcal{V}_p \left( u^{(1)},\ldots,u^{(p)}\right) }{p K^{p-1}},
\end{equation}
where $\mathcal{V}_p \left( u^{(1)},\ldots,u^{(p)} \right) $ denotes the Euclidean volume of the parallelepiped generated by the $p$ vectors $u^{(1)},\ldots,u^{(p)}$.
\end{lemma}
Recall that, for $ m \in [1,2]$, the internal frequency vector is given by
\begin{equation} \nonumber
\omega(m) \equiv ( \omega_a(m) )_{a \in \mathcal{A}} = ( \sqrt{a^2+m} )_{a \in \mathcal{A}}.
\end{equation}
\begin{corollary} \label{m1.2}
Let $n=\operatorname{Card}( \mathcal{A})$ and $w$ a nonzero vector in $\mathbb{R}^n$. Then, for any  $m\in [1,2]$, there exists  $j\in[1,...,n]$ such that
\begin{equation*}
\left|w\cdot\frac{d^j\omega}{dm^j}(m)\right|\geq
C N^{-2n^2} \Vert w \Vert_{\ell^1},
\end{equation*}
where $C>0$ is a constant that depends only on $n$.
\end{corollary}
\proof
From Lemma \ref{det}, $ \left( \frac{d\omega}{dm}(m), \ldots, \frac{d^n\omega}{dm^n}(m) \right)$ is a basis of $\mathbb{R}^n$. Therefore $w$ is as a linear combination of these vectors. According to Lemma~\ref{m1.1}, there exists $j \in [1 \ldots,n]$ such that
\begin{align*}
\left| w \cdot \frac{d^j\omega}{dm^j}(m) \right| & \geq\frac{\Vert w \Vert_{\ell^2} \mathcal{V}_n \left( \frac{d\omega}{dm}(m), \ldots, \frac{d^n\omega}{dm^n}(m)\right) }{n K^{n-1}}\\
& \geq \frac{\Vert w \Vert_{\ell^1} \mathcal{V}_n \left( \frac{d\omega}{dm}(m), \ldots, \frac{d^n\omega}{dm^n}(m)\right) }{n^{3/2} K^{n-1}}.
\end{align*}
Note that $ \left( \frac{d\omega}{dm}(m), \ldots, \frac{d^n\omega}{dm^n}(m) \right)$ is a $n$-family vector in $\mathbb{R}^n$. So
\begin{equation} \nonumber
\mathcal{V}_n \left( \frac{d\omega}{dm}(m), \ldots, \frac{d^n\omega}{dm^n}(m)\right) = D,
\end{equation}
where $D$ is the determinant defined in the Lemma~\ref{det}. Let us now give the expression of $K$. For $j \in [1,\ldots n]$ we have
\begin{align*}
\left \| \frac{d^j\omega}{dm^j}(m) \right \|_{\ell^1} & = \sum_{1 \leq k \leq n} \left| \frac{(2j-2)!}{2^{2j-1}(j-1)!}\frac{(-1)^{j+1}} {(k^2+m)^{j-\frac12}} \right| \\
& \leq \sum_{1 \leq k \leq n}  \frac{(2n-2)!}{(n-1)!} = \frac{ n (2n-2)!}{(n-1)!} =: K.
\end{align*}
Then
\begin{equation} \nonumber
\left| w \cdot \frac{d^j\omega}{dm^j}(m) \right| \geq \frac{(n-1)!}{n^{5/2}(2n-2)!} D \Vert w \Vert_{\ell^1} \geq C N^{-2n^2} \Vert w \Vert_{\ell^1},
\end{equation}
which ends the proof of the corollary.
\endproof
We need the following lemma $2.1$ from \cite{xu1997invariant}:
\begin{lemma}
\label{v.112}
Assume that $g(\tau)$ is a $p$-\textit{th} differentiable on $J\subset\mathbb{R}$. For $h>0$, we define the open set $J_h$
\begin{equation} \nonumber
J_h:=\left\lbrace  \tau \in J \mid \vert g(\tau)\vert < h \right\rbrace.
\end{equation}
If $\vert g^{(p)}(\tau) \vert \geq d>0$, for $\tau \in J$, then
\begin{equation} \label{estim mes}
\operatorname{mes} (J_h)\leq M h^{1/p},
\end{equation}
where $ M:=2(2+3+...+p+d^{-1}).$
\end{lemma}

Now we have all the tools to give a proof of the Proposition~\ref{NRom}.

\proof[Proof of Proposition~\ref{NRom}] Let $k \in \mathbb{R}^n$, where $n=\operatorname{Card}( \mathcal{A})$, and consider the function $g \in \mathcal{C}^\infty( [1,2],\mathbb{R})$ defined by:
\begin{equation*}
g(m)=  k \cdot \omega (m) + c.
\end{equation*}
From Corollary~\ref{m1.2}, there exists $j \in [1  \ldots n ]$ such that
\begin{equation} \nonumber
\left|k \cdot\frac{d^j\omega}{dm^j}(m)\right |\geq C N^{-2n^2} \vert k \vert.
\end{equation}
Using Lemma~\ref{v.112} with $h =  C N^{-2n^2} \vert k \vert$, we obtain
\begin{equation} \nonumber
 \operatorname{mes}\ \left\lbrace m\in[1,2]\ \mid \left|k \cdot \omega(m)+c\right|\leq {\chi}\right\rbrace \leq M \chi^{1/n},
\end{equation}
where $M \leq \tilde C N^{2n^2} \vert k \vert^{-1}$. The constant $\tilde C$  is strictly positive and depends only on~$n$. This ends the proof of the proposition.
\endproof

\subsection{ Small Divisors Estimates}

For $m \in [1,2]$, recall that the internal frequencies are denoted by $ \displaystyle \omega\equiv\omega(m)=(\sqrt{a^2+m})_{a\in\mathcal{A}}$, while the external frequencies are denoted by $\lambda_s\equiv\lambda_s(m) =\sqrt{s^2+m}$ for $s\in \mathcal{L}=\mathbb{Z} \setminus\mathcal{A}$. We Note that for $s \in \mathcal{L}  \setminus \lbrace 0 \rbrace$ we have
\begin{equation} \label{estim_dif_lam_enti}
\left| \lambda_s(m)-|s|\right|\leq \frac{m}{2 |s|}.
\end{equation}
We Recall that $\mathcal{A}^-:=\{s\in\mathcal{L}\mid -s \in \mathcal{A}\rbrace.$
We denote by $\mathcal{L}^\infty$ the complementary of $\mathcal{A}^-$  in $\mathcal{L}$ and $n$ the cardinality of  $\mathcal{A}$.

In this part, we will give a lower bound of the modulus of the following small divisors:
\begin{align*}
D_0=&\omega \cdot k, \quad k\in \mathbb{Z}^n\setminus\lbrace 0 \rbrace,\\
D_1=&\omega \cdot k+\lambda_a, \quad k\in \mathbb{Z}^n, \: a\in \mathcal{L},\\
D_2=&\omega\cdot k+\lambda_a+\lambda_b, \quad k\in \mathbb{Z}^n,\: a, b\in \mathcal{L},\\
D_3=&\omega\cdot k +\lambda_a-\lambda_b, \quad k\in \mathbb{Z}^n,\: a, b \in \mathcal{L}.
\end{align*}
\begin{definition}\label{Res-kab}
Let $k\in \mathbb{Z}^n,\: a, b\in \mathcal{L}$.
\begin{itemize}
\item[(i)] The vector $k$ is $D_0$ resonant if $k=0$.
\item[(ii)] The couple $(k;a)$ is $D_1$ resonant  if $|a|=|s|$ where $s\in\mathcal{A}$ and $\omega\cdot k=-\omega_s$
\item[(iii)] The triplet $(k;a,b)$ is $D_2$ resonant  if $|a|=|s|$, $|b|=|s'|$ where $s,s'\in\mathcal{A}$ and $\omega\cdot k=-\omega_s-\omega_{s'}$
\item[(iv)] The triplet $(k;a,b)$ is $D_3$ resonant  if $|a|=|s|$, $|b|=|s'|$ where $s,s'\in\mathcal{A}$ and $\omega\cdot k=-\omega_s+\omega_{s'}$.
\end{itemize}
\end{definition}
\begin{remark}\label{remL+}
Note that $(k;a,b)$ can be $D_2$ or $D_3$ resonant only if  $(a,b)\in\mathcal{A}^-\times\mathcal{A}^-$. Similarly, $(k;a)$ can be $D_1$  resonant only if $a\in\mathcal{A}^-$.
\end{remark}
 Our goal is to give a lower bound to the modulus of small divisors $D_0,\:D_1,\:D_2$ and $D_3$ when they are not resonant, for $ m \in [1,2] \setminus \mathcal{C} $, with $ \mathcal{C} $ an open set of zero Lebesgue measure.
In this section, $C$ will denotes a constant that depends only on the admissible set $ \mathcal{A}$. Let us start with the small divisors $D_0$, $D_1$ and $D_2$.
\begin{proposition} \label{final-petit diviseur}
Let $\kappa > 0$ and an integer $N>1$. Then there is an open set $\mathcal{C}\subset [1,2]$ that verifies
\begin{equation} \nonumber
\operatorname{mes} \left(  \mathcal{C} \right) \leq C \kappa^{\tau}N^{\iota},
\end{equation}
where $\tau,\iota>0$ and depend only on $n=\operatorname{Card} \left(  \mathcal{A}\right) $, such that for all $m\in \left( [1,2]\setminus \mathcal{C} \right) $, all $0<|k|\leq N$ and all $a,b \in \mathcal{L}$  we have:
\begin{equation} \label{petit div 1 bis}
|\omega \cdot k|\geq  \kappa,
\end{equation}
except when $k$  is $D_0$ resonant;
\begin{equation} \label{petit div 2 bis}
|\omega \cdot k+\lambda_a|\geq \kappa \langle a \rangle ,
\end{equation}
except when $(k;a)$ is $D_1$ resonant;
\begin{equation} \label{petit div 3 bis}
|\omega\cdot k+\lambda_a+\lambda_b|\geq \kappa\left( \langle a\rangle+\langle b\rangle\right) ,
\end{equation}
except when  $(k;a,b)$ is $D_2$ resonant.
The constant $C$ depends only on the admissible set $\mathcal{A}$.
\end{proposition}
\proof [Proof]
We start by proving \eqref{petit div 1 bis}. Let $\kappa > 0$ and an integer $N>1$. Consider
\begin{equation} \nonumber
\mathcal{U} = \lbrace m\in[1,2] \mid |\omega \cdot k| <  \kappa , \: k \in \mathbb{Z}^n \mbox{ for } 0 < \vert k \vert \leq N \rbrace .
\end{equation}
For $k \in \mathbb{Z}^n$ we consider the sets:
\begin{equation} \nonumber
\mathcal{U}_k= \lbrace m\in[1,2] \mid |\omega \cdot k| <  \kappa     \rbrace
\end{equation}
and
\begin{equation} \nonumber
\mathcal{B}_0 =\lbrace  k \in \mathbb{Z}^n | \vert k \vert \leq N \rbrace.
\end{equation}
Then
\begin{equation} \nonumber
\mathcal{U} = \bigcup_{k \in \mathcal{B}_0} \mathcal{U}_k.
\end{equation}
Thanks to Proposition~\ref{NRom}, we have $\operatorname{mes} \left( \mathcal{U}_k \right)  \leq C \frac{\kappa^{1/n}}{|k|}$. We Note that there are at most $N^n$ points in $\mathcal{B}_0$. So we obtain:
\begin{equation} \nonumber
\operatorname{mes} \left( \mathcal{U} \right)  \leq C \kappa^{1/n} N^n.
\end{equation}
Let us now look at the second small divisor \eqref{petit div 2 bis}. Consider
\begin{equation} \nonumber
C_{\mathcal{A}}= \left(  \max \lbrace \vert a \vert \mid a \in \mathcal{A} \rbrace^2   \right)^{1/2}.
\end{equation}
There are two cases: if $\vert a  \vert \geq 2 C_\mathcal{A} N$, then
\begin{equation} \nonumber
|\omega \cdot k+\lambda_a|\geq \lambda_a - \vert \omega \cdot k \vert \geq \vert a \vert -C_\mathcal{A} \vert k \vert  \geq \vert a \vert - \frac{1}{2} \vert a \vert \geq \kappa \langle a \rangle,
\end{equation}
for $0 < \kappa \leq 1/2$. If $\vert a  \vert < 2 C_\mathcal{A} N$, let
\begin{align*} \nonumber
\mathcal{V} =  \lbrace m & \in[1,2]  \mid \vert \omega \cdot k + \lambda_a \vert  <  \kappa \langle a  \rangle ,\:  (k;a) \in \mathcal{L} \times \mathbb{Z}^n \\ & \mbox{with } 0 \leq \vert k \vert \leq N  \mbox{ and } (k;a) \mbox{ non } D_1 \mbox{ resonant }\rbrace .
\end{align*}
We want to give an upper bound of the Lebesgue measure of $ \mathcal{V} $. Consider for $k \in \mathbb{Z}^n$ and $ a \in \mathcal{L}$ the sets
\begin{equation} \nonumber
\mathcal{V}_{k,a} = \lbrace m  \in[1,2]  \mid \vert \omega \cdot k + \lambda_a \vert  <  \kappa \langle a \rangle,\: (k;a) \mbox{ non } D_1 \mbox{ resonant }\rbrace
\end{equation}
and
\begin{equation} \nonumber
\mathcal{B}_1 = \lbrace  (k,a) \in \mathbb{Z}^n \times \mathcal{L} | \vert k \vert \leq N \mbox{ et } \vert a \vert < 2C_{\mathcal{A}} N    \rbrace.
\end{equation}
Then we have
\begin{equation}
\mathcal{V}  \subset \bigcup_{(k,a) \in \mathcal{B}_1} \mathcal{V}_{k,a}.\end{equation}
We note that there are at most $4C_{\mathcal{A}}N^{n+1}$ points in $\mathcal{B}_1$. It remains to give an upper bound of the Lebesgue measure of $\mathcal{V}_{k,a}$. There are two cases:
\begin{itemize}
\item If $\left\lbrace a , -a \right\rbrace \not\subset \mathcal{A}$, then $\mathcal{A}'=\mathcal{A}\cup \{a\}$ is still an admissible set of cardinality $n+1$. In addition, we have $\mathcal{A}' \subset \lbrace a \in\mathbb{Z} \mid |a|\leq C N\}$. Applying Proposition~\ref{NRom} to the new admissible set, we have
\begin{equation} \nonumber
\operatorname{mes} \left(  \mathcal{V}_{k,a} \right)  \leq C\frac{\kappa^{1/(n+1)}N^{2(n+1)^2+1/(n+1)}}{|k+1|}.
\end{equation}
\item If $|a|\in \mathcal{A}$ but $(k;a)$ is not $D_1$ resonant, then by applying Proposition~\ref{NRom} without changing $\mathcal{A}$ we have
\begin{equation} \nonumber
\operatorname{mes} \left(  \mathcal{V}_{k,a} \right) \leq C \frac{\kappa^{1/n}N^{2n^2+1/n}}{|k|}.
\end{equation}
\end{itemize}
So
\begin{equation} \nonumber
\operatorname{mes} \left( \mathcal{V} \right)   \leq C \kappa^{1/(n+1)} N^{(n+1)(2n+3)+1/(n+1)}.
\end{equation}
With the same argument we show \eqref{petit div 3 bis}. We end the proof of Proposition~\ref{final-petit diviseur} by taking $ \mathcal{C} = \mathcal{U} \cup \mathcal{V} \cup \mathcal{W}$ where $\mathcal{W}$ is the open set where \eqref{petit div 3 bis} is not verified.
\endproof
It remains to control  $D_3=\omega \cdot k+\lambda_a-\lambda_b$.
\begin{lemma} \label{Res D3}
Let $\tilde \kappa \in ]0,1]$ and  an integer $N \geqslant 1 $. We have
\begin{equation*}
\nonumber \operatorname{mes} \{m  \in[1,2] \mid |\omega \cdot k-e| < 2 \tilde \kappa, \; (k,e)\in\mathbb{Z}^{n+1}
\nonumber \mbox{ for }0<|k|\leq N\}\leq C \tilde \kappa^{\frac 1 {n}} N^{n+1} ,
\end{equation*}
where  $C>0$ and depends only on the admissible set $\mathcal{A}$.
\end{lemma}
\proof
Let $(k,e) \in \mathbb{Z}^n \times \mathbb{Z}$ such that $0<|k| \leq N$. Using Proposition~\ref{NRom}, we have
\begin{equation} \nonumber
\operatorname{mes} \lbrace m \in[1,2] \mid |\omega \cdot k-e| < 2 \tilde \kappa,  \rbrace \leq C \frac{\tilde{ \kappa}^{\frac 1 {n}}}{|k|}.
\end{equation}
Since $\kappa\leq 1$, we can restrict ourselves to
\begin{equation} \nonumber
|e|\leq |\omega \cdot k-e| + |\omega\cdot k| \leq   CN.
\end{equation}
Then
\begin{equation} \nonumber
\operatorname{mes} \bigcup_{\stackrel{\vert k \vert \leq  N}{(k,e) \in \mathbb{Z}^{n+1}}} \lbrace m \in[1,2] \mid |\omega \cdot k-e| < 2 \tilde \kappa  \rbrace \leq CN^{n+1} \tilde{\kappa}^{\frac 1 {n}}.
\end{equation}
The proof is thus concluded.
\endproof

\begin{proposition}\label{prop-D3}
Let $\kappa > 0$ and an integer $N>0$. Then there is an open set $\mathcal{C}\subset [1,2]$ satisfying
\begin{equation} \nonumber
\operatorname{mes} \left(  \mathcal{C} \right) \leq C \kappa^{\tau}N^{\iota},
\end{equation}
where $\tau$ and $ \iota$ are two strictly positive exponents which depend only on  $n=\operatorname{Card} \left(  \mathcal{A}\right) $, such that for all $m\in \left( [1,2]\setminus \mathcal{C} \right) $, all $0<|k|\leq N$ 	and all $a,b \in \mathcal{L}$ we have
\begin{equation} \label{petit div 4 bis}
|\omega\cdot k +\lambda_a-\lambda_b|\geq \kappa(1+ \vert \vert a \vert - \vert b \vert \vert),
\end{equation}
except when $(k;a,b)$ is $D_3$ resonant. The constant $C$ depends only on the admissible set $\mathcal{A}$.
\end{proposition}
\proof
Using \eqref{estim_dif_lam_enti} for $|b|\geq |a|>0$,  we remark that
\begin{equation}  \nonumber
|\lambda_a-\lambda_b-(|a|-|b|)|\leq \frac{m}{|a|}\leq 2|a|^{-1}.
\end{equation}
So we have
\begin{equation} \nonumber
|\omega\cdot k +\lambda_a-\lambda_b|\geq |\omega\cdot k +|a|-|b||-2|a|^{-1}.
\end{equation}
In Lemma~\ref{Res D3}, we denote $\tilde{\kappa} = \bar \kappa^\varrho$ where $\varrho$  is an exponent in $ ]0,1[$ which will be determined later. According to this Lemma, there is an open set $\mathcal C_1=\mathcal C_1(N, \bar\kappa^\varrho)$ whose Lebesgue measure is smaller than $C  \bar\kappa^{\frac \varrho {n}} N^{n+1}$, where $C$ is a constant that depends on $\mathcal{A}$. For all $m\in \left(  [1,2]\setminus \mathcal C_1\right) $, all $0<|k|\leq N$ and all $a,b \in \mathcal{L}$ where $|b|\geq|a|\geq 2 \bar\kappa^{-\varrho}$, we have:
\begin{equation} \label{intermediaire}
\vert \omega\cdot k +\lambda_a-\lambda_b \vert \geq  \bar \kappa^{\varrho} \geq \bar\kappa.
\end{equation}
Let us look at the remaining cases where the previous estimate does not hold. These cases are included in the set:
\begin{align*}
\mathcal{C}_2 = \left\lbrace  m \in [1,2] \mid  \vert \omega \cdot k  +  \lambda_a - \lambda_b \vert < \bar\kappa, \; (a,b) \in \mathcal{L}^2, \vert  a \vert \leq 2 \bar\kappa^{-\varrho} , 0< \vert k \vert \leq N \right\rbrace.
\end{align*}
We note that, if $|\omega\cdot k +\lambda_a-\lambda_b| < \bar\kappa$, $|a|\leq 2 \bar\kappa^{-\varrho}$ and $|k|\leq N$, then we have:
\begin{align*}
|b|&\leq \lambda_b \leq \vert \omega\cdot k +\lambda_a-\lambda_b \vert + \vert \omega \cdot k \vert + \lambda_a \\
& \leq 2 \kappa^{-\varrho} + (C_{\mathcal{A}} + 3 )N,
\end{align*}
where $C_{\mathcal{A}}= \left(  \max \lbrace \vert a \vert \mid a \in \mathcal{A} \rbrace^2   \right)^{1/2}$. Consider the set
\begin{equation} \nonumber
\mathcal B =\lbrace (a,b) \in \mathbb{Z}^2 \mid   |a|\leq |b|\leq 2 \bar \kappa^{-\varrho} + (C_{\mathcal{A}} + 3 )N \rbrace.
\end{equation}
There are at most  $4(2 \bar \kappa^{-\varrho} + (C_{\mathcal{A}} + 3 )N )^{2}$ points  in $\mathcal B$. So we have
\begin{equation} \nonumber
\mathcal{C}_2 \subset \left\lbrace  m \in [1,2] \mid  \vert \omega \cdot k  +  \lambda_a - \lambda_b \vert < \bar \kappa, \; (a,b) \in \mathcal{B}, 0< \vert k \vert \leq N \right\rbrace:= \mathcal{C}_3
\end{equation}
Recall that $\displaystyle \mathcal{L}^\infty = \mathcal{L} \setminus \mathcal{A}^-$ where $\mathcal{A}^-= -\mathcal{A}$. For  $a\in\mathcal{L}$, we define the set $]a[$  by : $]a[=\lbrace a \rbrace$ if $a\in\mathcal{L}^\infty$ and $]a[=\emptyset$ if $a\in\mathcal{A}^-$. We define
 \begin{equation} \nonumber
\mathcal{A}'= \mathcal{A}\cup]a[\cup]b[.
 \end{equation}
The set $\mathcal{A}'$ is admissible. In addition, we have
\begin{equation} \nonumber
\mathcal{A}' \subset \lbrace a \in \mathbb{Z} \mid \vert a \vert \leq C \left( 2 \bar \kappa^{-\varrho} + (C_{\mathcal{A}} + 3 )N \right)   \rbrace.
\end{equation}

The triplet $(k;a,b)$ is  $D_3$ non resonant. By applying the Proposition~\ref{NRom} with the admissible set $\mathcal{A}'$, we have:
\begin{align*}
\operatorname{mes}\left( \mathcal{C}_2\right)  & \leq \operatorname{mes}\left(\mathcal{C}_3\right)  \leq C \bar \kappa^{1/(n+2)}N^n \left( 2 \bar \kappa^{-\varrho} + (C_{\mathcal{A}} + 3 )N \right)^{2(n+2)^2}  \operatorname{Card} \mathcal B\\
&\leq C \bar \kappa^{1/(n+2)} \bar \kappa^{-2 \varrho \left( (n+2)^2+1 \right) } N^{(n+2)(2n+5)}.
\end{align*}
Let
\begin{equation}\nonumber
\varrho = \frac{1}{4 \left( (n+2)^2+1 \right) (n+2)}, \quad \tau =\frac{1}{2(n+2)},\quad \iota'=(n+2)(2n+5),
\end{equation}
and consider $\mathcal{C}= \mathcal{C}_1 \cup \mathcal{C}_2$, then we have:
\begin{equation} \nonumber
\operatorname{mes} \left(  \mathcal{C} \right)  \leq  C   \bar \kappa^{\tau} N^{\iota'}.
\end{equation}
Moreover, for all $m\in [1,2]\setminus \mathcal{C}_4$, all $0<|k|\leq N$ and all  $a,b \in \mathcal{L}$ non $D_3$ resonant, the estimation \eqref{intermediaire} is satisfied.

To conclude the proof of the proposition, we need to estimate the difference $\vert\vert a \vert - \vert b \vert\vert$. Without loss of generality, assume that $\vert a \vert > \vert b \vert$.
\begin{itemize}
\item If  $ \vert a \vert - \vert b \vert \geq 8 C_\mathcal{A} N \geq 8 \vert \omega \cdot k \vert$, then for $m \in [1,2]$ and $ 0 < \vert k \vert \leq N$, we have:
\begin{align*}
\vert \omega\cdot k +\lambda_a-\lambda_b  \vert &  \geq \lambda_a - \lambda_b - \vert \omega \cdot k \vert \geq \frac{1}{4} (\vert a \vert - \vert b \vert) - \vert  \omega \cdot k  \vert \\
& \geq \frac{1}{8} (\vert a \vert - \vert b \vert ) \geq \frac{1}{16} (1+ \vert a \vert - \vert b \vert ) \geq \kappa (1+ \vert a \vert - \vert b \vert ),
\end{align*}
if we assume that $\kappa \leq \frac{1}{16}$.
\item If $\vert a \vert - \vert b \vert < 8 C_\mathcal{A}N$, then for $ m \in [1,2] \setminus \mathcal{C}$, we have:
\begin{align*}
\vert \omega\cdot k +\lambda_a-\lambda_b  \vert &  \geq \frac{\bar\kappa}{1+8C_\mathcal{A} N} (1+ \vert a \vert - \vert b \vert ) \\
& \geq \kappa (1+ \vert a \vert - \vert b \vert ).
\end{align*}
\end{itemize}
Thus, for all $m\in [1,2]\setminus \mathcal{C}$, all $0<|k|\leq N$ and all  $a,b \in \mathcal{L}$ non $D_3$ resonant, the estimation \eqref{petit div 4 bis} holds. Moreover,
\begin{equation} \nonumber
\operatorname{mes} \left(  \mathcal{C} \right)  \leq  C \kappa^{\tau} N^{\iota'+\tau}= C \kappa^{\tau} N^{\iota}.
\end{equation}
\endproof
It remains to treat the case where $k=0$ in $D_3$.
\begin{lemma}\label{lem:D3-k=0}
Let $m\in [1,2]$ and $a,b\in\mathcal{L}$ such that $|a|\neq|b|$. Therefore
\begin{equation} \nonumber
\vert \lambda_a - \lambda_b \vert \geq \frac{1}{8}(1+\vert \vert a \vert - \vert b \vert \vert) .
\end{equation}
\end{lemma}
\proof Without loss of generality, assume that  $|a|>|b|$. Then for all $ m \in [1,2]$, we have:
\begin{equation} \nonumber
\lambda_a - \lambda_b = \frac{(\vert a \vert - \vert b \vert)(\vert a \vert + \vert b \vert) }{\sqrt{a^2+m}+\sqrt{b^2+m}} \geq \frac{1}{4} (\vert a \vert - \vert b \vert) \geq \frac{1}{8}(1+\vert \vert a \vert - \vert b \vert \vert),
\end{equation}
which concludes the proof.
\endproof

\section{Normal form}
In this section, we construct a symplectic change of  variable that puts the Hamiltonian \eqref{hamiltonien} in normal form to which we can apply our KAM theorem.
\subsection{Class of Hamiltonian function}

In this part, we begin by recalling some notations introduced in \cite{B1}. For $\mathcal{L}$ a set of $\mathbb{Z}$ and $\alpha\geq 0$, we define the $\ell_2$ weighted space:
\begin{equation} \nonumber
Y_\alpha:= \left\lbrace  \zeta= \left( \zeta_{s} = \begin{pmatrix}
\xi_s \\
\eta_s
\end{pmatrix} \in \mathbb{C}^2  , s \in \mathcal{L}  \right) | \; \Vert \zeta \Vert_\alpha < \infty \right\rbrace   ,
\end{equation}
where
\begin{equation} \nonumber
\Vert \zeta \Vert_\alpha^2 = \sum_{s \in \mathcal{L}} \vert \zeta_s \vert^2 \langle s\rangle^{2\alpha},\: \mbox{où} \quad \langle s\rangle = \max(\vert s \vert , 1 ).
\end{equation}
We endow $\mathbb{C}^2$ with the euclidean norm, i.e. if $\zeta_s={}^t(\xi_s,\eta_s)$ then $|\zeta_s|=\sqrt{|\xi_s|^2+|\eta_s|^2}$.\\
For $\beta\geq 0$, we define the $\ell_\infty$ weighted space
\begin{equation} \nonumber
L_\beta= \left\lbrace   \left( \zeta_{s}  = \begin{pmatrix}
\xi_s \\
\eta_s
\end{pmatrix} \in \mathbb{C}^2  , s \in \mathcal{L}  \right) \; | \; \vert \zeta \vert_\beta < \infty \right\rbrace   ,
\end{equation}
where
\begin{equation} \nonumber
\vert \zeta \vert_\beta= \underset{s \in \mathcal{L}}{\sup} \vert \zeta_s \vert \langle s\rangle^{\beta} .
\end{equation}
We remark that for,  $\beta \leq \alpha$, we have $Y_\alpha \subset L_\beta.$\\
\textbf{Infinite matrices.} Consider the orthogonal projector $\Pi$ defined on the set of square matrices by
\begin{equation}\nonumber
\Pi: \: \mathcal{M}_{2\times2}(\mathbb{C}) \rightarrow \mathbb{S},
\end{equation}
where
\begin{equation} \nonumber
\mathbb{S}=\mathbb{C}I+\mathbb{C}\sigma_2, \quad \text{with} \quad \sigma_2=\begin{pmatrix}
   0 & -1 \\
   1 & 0
\end{pmatrix}.
\end{equation}
We introduce $\mathcal{M}$ the set  of infinite symmetric matrices $A: \mathcal{L} \times \mathcal{L} \rightarrow \mathcal{M}_2 \left( \mathbb{R} \right)$, that verify, for any $ s, s' \in \mathcal{L}$,
\begin{center}
 $A_s^{s'} \in \mathcal{M}_2 \left( \mathbb{R} \right)$,
 $A_s^{s'} = {}^tA_{s'}^s$ and  $\Pi A_{s}^{s'}=A_{s}^{s'}.$
\end{center}
We also define $\mathcal{M}_\alpha$, a subset of $\mathcal{M}$, by:
\begin{equation} \nonumber
A \in \mathcal{M}_\alpha \Leftrightarrow \vert A \vert_\alpha:= \underset{s,s' \in \mathcal{L}}{\sup} \langle s\rangle^{\alpha} \langle s'\rangle^{\alpha} \Vert A_s^{s'} \Vert_\infty < \infty .
\end{equation}
Let $n \in \mathbb{N}$, $\rho > 0$ and  $B$ be a Banach space. We define:\begin{equation} \nonumber
\mathbb{T}^n_\rho= \lbrace  \theta \in \mathbb{C}^n/ 2\pi \mathbb{Z}^n \vert \: \vert Im\theta \vert < \rho\rbrace
\end{equation}
and
\begin{equation} \nonumber
\mathcal{O}_\rho \left( B \right) = \left\lbrace x \in B \vert \Vert x \Vert_B < \rho \right\rbrace .
\end{equation}
For $\sigma , \mu \in  \left] 0 , 1\right[ $, we define
\begin{equation} \nonumber
\mathcal{O}^\alpha (  \sigma , \mu )  = \mathbb{T}^n_\sigma \times \mathcal{O}_{\mu^2} ( \mathbb{C}^n ) \times \mathcal{O}_\mu ( Y_\alpha )= \lbrace ( \theta,r,\zeta) \rbrace,
\end{equation}
\begin{equation} \nonumber
\mathcal{O}^{\alpha, \mathbb{R}}(  \sigma , \mu )= \mathcal{O}^\alpha (  \sigma , \mu ) \cap \lbrace \mathbb{T}^n \times \mathbb{R}^n \times Y_\alpha^{\mathbb{R}} \rbrace,
\end{equation}
where $ Y_\alpha^{\mathbb{R}}= \left\lbrace  \zeta \in Y_\alpha \: | \: \zeta = \left( \zeta_{s} = \begin{pmatrix}
\xi_s \\
\eta_s
\end{pmatrix}, \xi_s=\bar{\eta}_s \: s \in \mathcal{L}  \right) \right\rbrace$.\\
Let us denote a point in $\mathcal{O}^\alpha (  \sigma , \mu )$ as $x= (\theta,r, \zeta)$. A function on $\O^\alpha(\s,\mu)$ is real if it has a real value for any real $x$. We define:
\begin{equation} \nonumber
\Vert (r,\theta,\zeta) \Vert_\alpha=\max(|r|,|\theta|,\Vert \zeta \Vert_\alpha).
\end{equation}

\textbf{Class of Hamiltonian functions.}
Let $\mathcal{D}$ be a compact set of $\mathbb{R}^p$, called the parameters set from now on. Let $f:\mathcal{O}^{\alpha}(  \delta , \mu ) \times \mathcal{D} \rightarrow \mathbb{C}$ be a $\mathcal{C}^1$ function, real and holomorphic in the first variable, such that for all $\rho \in \mathcal{D}$, the maps
\begin{equation} \nonumber
\mathcal{O}^{\alpha}(  \delta , \mu ) \ni x \mapsto \nabla_\zeta f(x,\rho) \in Y_\alpha \cap L_\beta
\end{equation}
and
\begin{equation} \nonumber
\mathcal{O}^{\alpha}(  \delta , \mu ) \ni x \mapsto \nabla_{\zeta}^2f(x,\rho) \in \mathcal{M}_\beta,
\end{equation}
are holomorphic.  We define:
\begin{align*}
\left\vert f(x,.) \right\vert_\mathcal{D} = \underset{\rho \in \mathcal{D}}{sup} \left\vert f(x,\rho) \right\vert, &\quad \left\Vert \frac{\partial f}{\partial \zeta} (x,.) \right \Vert_\mathcal{D} = \underset{\rho \in \mathcal{D}}{sup}\left\Vert  \nabla_\zeta f(x,\rho) \right\Vert_\alpha ,\\
\left\vert \frac{\partial f}{\partial \zeta} (x,.) \right \vert_\mathcal{D} = \underset{\rho \in \mathcal{D}}{sup}\left\vert \nabla_\zeta f(x,\rho) \right\vert_\beta , &\quad \left\vert \frac{\partial^2 f}{\partial \zeta^2} (x,.) \right\vert_\mathcal{D} = \underset{\rho \in \mathcal{D}}{sup}\left\vert \ \nabla_\zeta^2 f(x,\rho) \right\vert_\beta.
\end{align*}
We denote by $\mathcal{T}^{\alpha,\beta}(\mathcal{D},\sigma,\mu)$ the space of functions $f$ that verify, for all $x \in \mathcal{O}^\alpha(\sigma,\mu)$, the following estimates:
\begin{equation} \nonumber
\left\vert f(x,.)\right\vert_\mathcal{D} \leq C, \quad \left\Vert \frac{\partial f}{\partial \zeta}  (x,.) \right\Vert_\mathcal{D} \leq \frac{C}{\mu}, \quad \left\vert \frac{\partial f}{\partial \zeta}  (x,.) \right\vert_\mathcal{D} \leq \frac{C}{\mu},\quad \left\vert \frac{\partial^2 f}{\partial \zeta^2}  (x,.) \right\vert_\mathcal{D} \leq \frac{C}{\mu^2}.
\end{equation}

For $f\in \mathcal{T}^{\alpha,\beta}(\mathcal{D},\sigma,\mu)$, we denote by $\lc f \rc^{\alpha,\beta} _{\sigma,\mu,\mathcal{D}}$ the smallest constant $C$ that satisfies the above estimates. If $\partial_\rho^j f \in \mathcal{T}^{\alpha,\beta}(\mathcal{D},\sigma,\mu)$ for $j\in \lbrace 0,1 \rbrace$, then for $\gamma >0$ we define:
\begin{equation} \nonumber
\lc f \rc^{\alpha,\beta,\gamma} _{\sigma,\mu,\mathcal{D}}=\lc f \rc^{\alpha,\beta} _{\sigma,\mu,\mathcal{D}}+\gamma \lc \partial_\rho f \rc^{\alpha,\beta} _{\sigma,\mu,\mathcal{D}}.
\end{equation}
We also denote by
\begin{equation} \nonumber
\mathcal{T}^{\alpha,\beta}(\mu)= \left\lbrace f(\zeta)\: | \: f \in \mathcal{T}^{\alpha,\beta}(\mathcal{D},\sigma,\mu) \right\rbrace,
\end{equation}
the set of functions of $\mathcal{T}^{\alpha,\beta}(\mathcal{D},\sigma,\mu)$ that do not depend on $r,\theta$ and $\rho$.
The norm of such functions will be denoted by $\lc f \rc^{\alpha,\beta} _{\mu}$.

We finish this part by defining the space $\mathcal{T}^{\alpha,\beta+}(\mathcal{D},\sigma,\mu)$.
Consider the following spaces
\begin{equation} \nonumber
L_{\beta+}=\lbrace \zeta= \left( \zeta_{s} =\left( p_{s}, q_{s} \right)  , s \in \mathcal{L}  \right) | \quad \vert \zeta \vert_{\beta+} < \infty \rbrace,
\end{equation}
where $ \vert \zeta \vert_{\beta+}= \underset{s \in \mathcal{L}}{\sup} \vert \zeta_s \vert \langle s\rangle^{\beta+1},$
and
\begin{equation} \nonumber
\mathcal{M}_{\beta+}=\lbrace A \in \mathcal{M} | \: \vert A \vert_{\beta+} < \infty \rbrace,
\end{equation}
where $ \vert A \vert_{\beta+} = \underset{s,s' \in \mathcal{L}}{\sup} (1+| \: | s | - | s'| \:  |)\langle s\rangle^{\beta} \langle s'\rangle^{\beta} \Vert A_s^{s'} \Vert_\infty.$
\\We remark that $L_{\beta+} \subset L_\beta$ and  $\mathcal{M}_{\beta+} \subset \mathcal{M}_{\beta}$. We define  $\mathcal{T}^{\alpha,\beta+}(\mathcal{D},\sigma,\mu)$ the same way that we defined $\mathcal{T}^{\alpha,\beta}(\mathcal{D},\sigma,\mu)$, but replacing $L_{\beta}$ by $L_{\beta+}$ and $\mathcal{M}_{\beta}$ by $\mathcal{M}_{\beta+}$. So, we have $\mathcal{T}^{\alpha,\beta+}(\mathcal{D},\sigma,\mu) \subset \mathcal{T}^{\alpha,\beta}(\mathcal{D},\sigma,\mu).$
\\ For $f, g \in \mathcal{T}^{\alpha,\beta}(\mu)$, we define the Poisson bracket by:
\begin{equation}\nonumber
\left\lbrace f,g\right\rbrace = i \langle \nabla_\zeta f , J \nabla_\zeta g \rangle.
\end{equation}
\begin{lemma} \label{estim crochet poisson zeta}
Consider $ f \in \mathcal{T}^{\alpha,\beta}(\mu)$ and $ g \in \mathcal{T}^{\alpha,\beta+}(\mu)$. Then, for any $0< \mu' < \mu$, $ \left\lbrace f,g\right\rbrace \in \mathcal{T}^{\alpha,\beta}(\mu')$, we have:
\begin{equation} \nonumber
\lc \lbrace f,g \rbrace\rc^{\alpha,\beta} _{\mu'} \leq \frac{C}{\mu(\mu-\mu')} \lc f \rc^{\alpha,\beta} _{\mu} \lc g \rc^{\alpha,\beta+} _{\mu},
\end{equation}
where the constant $C$ depends on $\alpha$ on $\beta$.
\end{lemma}

For the proof, we recall the following lemma from \cite{Poschel1996} (appendix A).
\begin{lemma} \label{Cauchy_lineaire}
Let $E$ and $F$ be two complex Banach spaces, $f: E \to F$ and $v \in E$. Assume that there exists $r>0$ such that $f$ is holomorphic on the open ball of center $v$ and radius $r$ and satisfies $\Vert f \Vert_F \leq M$ on this ball. Then $d_vf \in \mathcal{L}(E,F)$, and we have:
\begin{equation} \nonumber
\Vert d_v f \Vert_{\mathcal{L}(E,F)} \leq \frac{M}{r}.
\end{equation}
\end{lemma}
\proof[Proof of Lemma \ref{estim crochet poisson zeta}]
Let $x \in \mathcal{O}_{\mu'} (Y_\alpha)$. Our goal is to prove that
\begin{itemize}
\item[(i)] $\displaystyle \left| \left\lbrace f,g \right\rbrace (x) \right| \leq \frac{C}{\mu(\mu-\mu')} \lc f \rc^{\alpha,\beta} _{\mu} \lc g \rc^{\alpha,\beta+} _{\mu} $;
\item[(ii)] $ \displaystyle \Vert \nabla_\zeta \left\lbrace f,g \right\rbrace (x) \Vert_\alpha \leq  \frac{C}{\mu \mu'(\mu-\mu')} \lc f \rc^{\alpha,\beta} _{\mu} \lc g \rc^{\alpha,\beta+} _{\mu} $;
\item[(iii)] $\displaystyle  \vert \nabla_\zeta \left\lbrace f,g \right\rbrace (x) \vert_\beta \leq \frac{C}{\mu \mu'(\mu-\mu')} \lc f \rc^{\alpha,\beta} _{\mu} \lc g \rc^{\alpha,\beta+} _{\mu} $;
\item[(iv)] $\displaystyle \vert \nabla_\zeta^2 \left\lbrace f,g \right\rbrace (x) \vert_\beta \leq \frac{C}{\mu \mu'^2(\mu-\mu')} \lc f \rc^{\alpha,\beta} _{\mu} \lc g \rc^{\alpha,\beta+} _{\mu} $.
\end{itemize}
\vspace{0.3cm}
Let us begin with the first estimation $(i)$. We have
\begin{align*}
\left| \left\lbrace f,g \right\rbrace (x) \right| & = \left| \langle \nabla_\zeta f(x) , J \nabla_\zeta g(x) \rangle \right| \leq \Vert \nabla_\zeta f(x) \Vert_\alpha \Vert J \nabla_\zeta g(x) \Vert_\alpha \\
& \leq \frac{1}{\mu^2} \lc f \rc^{\alpha,\beta} _{\mu} \lc g \rc^{\alpha,\beta+} _{\mu}\\
& \leq  \frac{1}{\mu(\mu-\mu')}\lc f \rc^{\alpha,\beta} _{\mu} \lc g \rc^{\alpha,\beta+} _{\mu}.
\end{align*}
Let us now turn to the $\zeta$-gradient of the Poisson bracket:
\begin{equation}\nonumber
\nabla_\zeta \left\lbrace f,g \right\rbrace (x) = \langle \nabla_\zeta^2 f(x) , J \nabla_\zeta g(x) \rangle + \langle \nabla_\zeta f(x) , J \nabla_\zeta^2 g(x) \rangle =: \Sigma_1 + \Sigma_2.
\end{equation}
For $\Sigma_1$, we have $\nabla_\zeta f: \mathcal{O}_{\mu'}(Y_\alpha) \to Y_\alpha$. Moreover, $x \mapsto \nabla_\zeta f(x)$ is holomorphic, so   $\nabla_\zeta^2 f(x) \in \mathcal{L}(Y_\alpha, Y_\alpha)$ for $x \in \mathcal{O}_{\mu'}(Y_\alpha)$. Using Lemma~\ref{Cauchy_lineaire}, we have:
\begin{align*}
\Vert \Sigma_1 \Vert_\alpha & \leq \Vert \nabla^2_\zeta f(x) \Vert_{\mathcal{L}(Y_\alpha, Y_\alpha)} \Vert J \nabla \zeta g(x) \Vert_\alpha\\
& \leq \frac{1}{\mu-\mu'}\sup_{y  \in  \mathcal{O}_\mu(Y_\alpha)} \left(\Vert \nabla_\zeta f(y) \Vert_\alpha \right) \Vert \nabla_\zeta g(x) \Vert_\alpha \\
& \leq \frac{1}{\mu^2(\mu-\mu')} \lc f \rc^{\alpha,\beta} _{\mu} \lc g \rc^{\alpha,\beta+} _{\mu} \\
& \leq \frac{1}{\mu \mu'(\mu-\mu')} \lc f \rc^{\alpha,\beta} _{\mu} \lc g \rc^{\alpha,\beta+} _{\mu}.
\end{align*}
We use the same arguments for $\Sigma_2$, which ends the proof of (ii). To prove (iii), we use estimations 2. and 3. from Lemma 2.1 in \cite{B1}.
So we have:
\begin{align*}
\left| \nabla_\zeta \left\lbrace f,g \right\rbrace (x) \right|_\beta & \leq \left| \langle \nabla_\zeta^2 f(x) , J \nabla_\zeta g(x) \rangle  \right|_\beta  + \left| \langle \nabla_\zeta f(x) , J \nabla_\zeta^2 g(x) \rangle  \right|_\beta \\
& \leq C \left| \nabla_\zeta^2 f (x) \right|_\beta \left| J \nabla_\zeta g(x) \right|_{\beta+} + C \left| \nabla_\zeta f(x) \right|_\beta \left| J \nabla^2_\zeta g(x) \right|_{\beta+}   \\
& \leq  \frac{C}{\mu^3} \lc f \rc^{\alpha,\beta} _{\mu} \lc g \rc^{\alpha,\beta+} _{\mu} \\
& \leq \frac{C}{\mu\mu'(\mu-\mu')} \lc f \rc^{\alpha,\beta} _{\mu} \lc g \rc^{\alpha,\beta+} _{\mu}.
\end{align*}
\vspace{0.3cm}
It remains to prove estimation (iv). We start by computing the second derivative of the Poisson bracket:
\begin{align*}
\nabla_\zeta^2 \left\lbrace f,g \right\rbrace (x) & = \langle \nabla_\zeta^3 f(x) , J \nabla_\zeta g(x) \rangle + \langle \nabla_\zeta^2 f(x) , J \nabla_\zeta^2 g(x) \rangle \\
& + \langle \nabla_\zeta^2 f(x) , J \nabla_\zeta^2 g(x) \rangle + \langle \nabla_\zeta f(x) , J \nabla_\zeta^3 g(x) \rangle \\
& =: \Gamma_1 +\Gamma_2 + \Gamma_3 + \Gamma_4.
\end{align*}
For $\Gamma_1$, we have $\nabla^2_\zeta f: \mathcal{O}_{\mu'}(Y_\alpha) \to \mathcal{M}_\beta$. Moreover, $x \mapsto \nabla^2_\zeta f(x)$ is holomorphic,  so $\nabla_\zeta^3 f(x) \in \mathcal{L}(Y_\alpha, \mathcal{M}_\beta)$ for $x \in \mathcal{O}_{\mu'}(Y_\alpha)$. By Lemma~\ref{Cauchy_lineaire}, we have
\begin{align*}
\left| \Gamma_1 \right| & \leq \Vert \nabla_\zeta^3 f(x) \Vert_{\mathcal{L}(Y_\alpha,\mathcal{M}_\beta)} \vert J \nabla_\zeta g(x) \Vert_\alpha \\& \leq  \frac{1}{\mu-\mu'} \sup_{y  \in  \mathcal{O}_\mu(\mathcal{M}_\beta)} \left| \nabla^2_\zeta f (y) \right|_\beta \Vert \nabla_\zeta g(x) \Vert_\alpha \\
& \leq \frac{1}{\mu^3 (\mu-\mu')} \lc f \rc^{\alpha,\beta}_{\mu} \lc g \rc^{\alpha,\beta+}_{\mu} \\
& \leq \frac{1}{\mu \mu'^2  (\mu-\mu')} \lc f \rc^{\alpha,\beta}_{\mu} \lc g \rc^{\alpha,\beta+}_{\mu}.
\end{align*}
We use the same arguments for $\Gamma_4$. It remains to estimate $\Gamma_2$ and $\Gamma_3$. The two cases are treated in the same way. Let us look for example at $\Gamma_3$. Using the first estimation from Lemma 2.1 in \cite{B1}, we have
\begin{align*}
\left| \Gamma_3 \right| & \leq C \left| \nabla_\zeta^2 f(x)\right|_\beta \left| \nabla^2_\zeta g(x)  \right|_{\beta+}\\
& \leq \frac{C}{\mu^4} \lc f \rc^{\alpha,\beta}_{\mu} \lc g \rc^{\alpha,\beta+}_{\mu} \\
& \leq \frac{C}{\mu \mu'^2  (\mu-\mu')} \lc f \rc^{\alpha,\beta}_{\mu} \lc g \rc^{\alpha,\beta+}_{\mu}.
\end{align*}
The proof is thus concluded.
\endproof
Consider $g$ a $\mathcal{C}^1$-function on $\mathcal{O}_{\mu}(Y_\alpha)$. We denote by $\Phi_g$ the Hamiltonian  flow of $g$ at $t=1$, i.e.
\begin{equation} \nonumber
\zeta = \zeta(1) = \Phi_g(\zeta_0)=\Phi^{t=1}_g(\zeta_0),
\end{equation}
where
\begin{equation} \nonumber
\dot \zeta = i J \nabla_\zeta g(\zeta) \quad \mbox{and} \quad \zeta(0)=\zeta_0.
\end{equation}
\begin{corollary} \label{estim flot zeta}
Let $ f \in \mathcal{T}^{\alpha,\beta}(\mu)$, $ g \in \mathcal{T}^{\alpha,\beta+}(\mu)$ and $0 < \mu' < \mu $ such that:
\begin{equation} \nonumber
\lc g \rc^{\alpha,\beta+}_{\mu} \leq \frac{1}{C} \mu (\mu-\mu').
\end{equation}
Then $\Phi_g: \mathcal{O}_{\mu'}(Y_\alpha) \to \mathcal{O}_{\mu}(Y_\alpha)$ and $f \circ \Phi_g \in \mathcal{T}^{\alpha,\beta}(\mu')$. In addition, we have:
\begin{equation} \label{flot zeta}
\lc f \circ \Phi_g \rc^{\alpha,\beta}_{\mu'} \leq C' \lc f \rc^{\alpha,\beta}_{\mu},
\end{equation}
where $C$ is a constant that depends on $\alpha$ and $\beta$ while $C'$ is an absolute constant.

\end{corollary}
\proof
Let us first prove that $\Phi_g$ is well-defined and that $\Phi_g: \mathcal{O}_{\mu'}(Y_\alpha) \to \mathcal{O}_{\mu}(Y_\alpha)$.
Consider
\begin{equation} \nonumber
\bar{t} = \sup \lbrace t> 0 \: | \: \zeta(a) \mbox{ is well-defined for } 0\leq a \leq t\mbox{ and } \Vert \zeta ( a ) \Vert_\alpha < \mu   \rbrace .
\end{equation}
Let us prove that $\bar t\geq 1$. We have:
\begin{equation} \nonumber
\zeta(a)= \zeta_0 + i \int_0^a J \nabla_\zeta g(\zeta(s))ds.
\end{equation}
Recall that $g \in \mathcal{T}^{\alpha,\beta+}(\mu)$, then for $\zeta \in \mathcal{O}_{\mu}(Y_\alpha)$ we have
\begin{equation} \nonumber
\Vert \nabla_\zeta g(\zeta)  \Vert_\alpha \leq \frac{\lc g \rc^{\alpha,\beta+}_{\mu}}{\mu} \leq \frac{1}{C}(\mu - \mu').
\end{equation}
It follows that, for $a=1$,
\begin{align*}
\Vert \zeta ( a ) \Vert_\alpha &\leq \Vert \zeta_0 \Vert_\alpha + \frac{a}{C}(\mu - \mu') \\
& \leq \mu' +  \frac{a}{C}(\mu - \mu') < \mu,
\end{align*}
We deduce that $\bar t\geq1$, $\Phi_g$ is well-defined and $\Phi_g: \mathcal{O}_{\mu'}(Y_\alpha) \to \mathcal{O}_{\mu}(Y_\alpha)$.
Let us prove now that  $f \circ \Phi_g \in \mathcal{T}^{\alpha,\beta}(\mu')$, as well as the estimation \eqref{flot zeta}. We recall that, for two holomorphic functions $f$ and  $g$, we have:
\begin{align*}
f \circ \Phi_g & = f + \left\lbrace f,g\right\rbrace + \frac{1}{2!} \left\lbrace \left\lbrace f,g\right\rbrace , g \right\rbrace  +\frac{1}{3!} \left\lbrace \left\lbrace \left\lbrace f,g\right\rbrace ,g \right\rbrace ,g \right\rbrace  + \ldots \\
& = \sum_{n \geq 0} \frac{1}{n!} P^n_g f,
\end{align*}
where $P^0_gf=f$, $P^1_gf=\left\lbrace f,g\right\rbrace $ et $P^2_gf= \left\lbrace \left\lbrace f,g\right\rbrace , g \right\rbrace \ldots$
\\ Using Lemma~\ref{estim crochet poisson zeta}, we have:
\begin{equation} \nonumber
\lc P^n_g f \rc^{\alpha,\beta} _{\mu'} \leq C^n \lc f \rc^{\alpha,\beta} _{\mu} \left( \mu (\mu - \mu')  \lc g \rc^{\alpha,\beta+} _{\mu}\right)^n.
\end{equation}
Thanks to the assumption made on $g$, we obtain that $f \circ \Phi_g$ is a convergent series. So $f \circ \Phi_g \in \mathcal{T}^{\alpha,\beta}(\mu')$ and satisfies estimation \eqref{flot zeta}.
\endproof
We define the real finite-dimensional torus:
\begin{equation} \label{tore inv}
T_\rho= \left\lbrace  \zeta=((\xi_s,\eta_s),\:s \in \mathbb{Z}) | \: \xi_s=\bar{\eta}_s ,\: \vert \zeta_a \vert^2= \nu \rho_a  \mbox{ si }  a \in \mathcal{A}, \: \xi_s=0 \mbox{ si } s  \in \mathcal{L} \right\rbrace ,
\end{equation}
 where $\nu >0$ is small and $\rho = (\rho_a , a \in \mathcal{A})$ is a parameter vector that belongs to $\mathcal{D}=\left[1,2 \right] ^\mathcal{A}$.

Let $n=\operatorname{Card}( \mathcal{A})$. The $n$-dimensional torus $T_\rho$ is invariant for the linear wave equation. We wish to put the Hamiltonian  $H_2+P$ in a suitable normal form to which we will apply the KAM theorem \ref{theoreme kam}. This normal form will be defined on tori constructed on the space $Y_\alpha$ and in the vicinity of the real torus $T_\rho$.\\
In the vicinity of the real torus $T_\rho$, we change from variables $(\zeta_a, a \in \mathcal{A})$ to action-angle variables defined by:
\begin{equation} \nonumber
\xi_a=\sqrt{I_a} e^{i\theta_a},\quad \eta_a=\sqrt{I_a} e^{-i\theta_a}, \quad a \in \mathcal{A}.
\end{equation}
So we change from variables $(\xi,\eta)$ to the new variables $(I,\theta, \xi ,\eta)$ where $I=(I_a, a \in \mathcal{A})$, $\theta=(\theta_a, \in \mathcal{A})$, $\xi=(\xi_s, \in \mathcal{L})$ and $\eta=(\eta_a, a \in \mathcal{L})$. The new variable vector is real if $I=\bar{I}$, $\theta=\bar{\theta}$, and $\xi=\bar{\eta}$.\\
We now define a complex toroidal vicinity of the real torus $T_\rho$ by
\begin{equation} \label{tore-complex}
\mathfrak{T}_\rho=\mathfrak{T}_\rho(\nu,\sigma,\mu,\alpha)=\left\lbrace (I,\theta,\xi,\eta) | \:  \vert I-\nu\rho \vert < \nu \mu^2, \vert Im \theta \vert < \sigma, \Vert \zeta^\mathcal{L} \Vert_\alpha < \nu^{1/2}\mu  \right\rbrace,
\end{equation}
where $I=(I_a, a \in \mathcal{A})$, $\theta=(\theta_a, \in \mathcal{A})$, $\zeta^\mathcal{L}=(\zeta_s, \in \mathcal{L})$.
\begin{lemma} \label{Perturbation-in-tore}
Assume that $g$ is a real holomorphic fuction on $\mathbb{S}^1 \times J$, for $J$ some neighborhood of the origin on $ \mathbb{R}$. Let $\alpha > 0$ and $\nu>0$ small. There exist $\sigma^* >0$ and $\mu^* >0$  such that the perturbation $P$ is well-defined and analytic on $\mathfrak{T}_\rho(\nu,\sigma,\mu,\alpha)$ for $0< \sigma \leq \sigma^*$ and $0<\mu\leq \mu^*$. The parameters $\sigma^*$ and $\mu^*$ depend on the nonlinearity $g$, the admissible set $\mathcal{A}$, $\nu$ and $\alpha$.
\end{lemma}
\begin{remark}
\leavevmode
\begin{itemize}
\item[-] We can fix $\sigma^*$ ( $\sigma^*=1 $ for example) and explicitly determine $\mu^*$.
\item[-] For small $\nu$, we have:
\begin{equation} \nonumber
\mathfrak{T}_\rho(\nu,\sigma,\mu,\alpha) \subset \mathcal{O}^\alpha \left(   \sigma , \mu^* \right).
\end{equation}
\end{itemize}
\end{remark}
\proof
The nonlinearity $g$ is real holomorphic on $\mathbb{S}^1 \times J$, for $J$ some neighborhood of the origin on $\mathbb{R}$. Assume that $I=\left[ -M,M \right]$ for $M>0$. We can extend holomorphically $g$ on $\mathbb{S}^1  \times I_{\mathbb{C}}$, for some $I_{\mathbb{C}}$ of the form
\begin{equation} \nonumber
I_{\mathbb{C}}=\left\lbrace u \in \mathbb{C} | \: \vert Re(u) \vert \leq M, \: \vert Im(u)  \vert < K  \right\rbrace.
\end{equation}
Recall that
\begin{align*}
u(\zeta)(x) & = \sum_{s \in \mathbb{Z}} \frac{\xi_s\varphi_s+\eta_s \varphi_{-s}}{\sqrt{2\lambda_s}} \\
& = \sum_{a \in \mathcal{A}} \frac{\xi_a\varphi_a+\eta_a \varphi_{-a}}{\sqrt{2\lambda_a}}  + \sum_{s \in \mathbb{Z}} \frac{\xi_s\varphi_s+\eta_s \varphi_{-s}}{\sqrt{2\lambda_s}} .
\end{align*}
We want to control these two sums. For the first sum we have:
\begin{align*}
\left| \sum_{a \in \mathcal{A}} \frac{\xi_a\varphi_a+\eta_a \varphi_{-a}}{\sqrt{2\lambda_a}} \right| & \leq \sum_{a \in \mathcal{A}} |\xi_a| +|\eta_a| \\ & \leq 2 \sum_{a \in \mathcal{A}} \sqrt{I_a}e^{|Im(\theta_a)|} \leq C_\mathcal{A} |I|^{1/2} e^{|Im(\theta)|},
\end{align*}
where $ C_\mathcal{A}$ is a constant that depends on $\mathcal{A}$. For the second sum, using the Cauchy-Schwarz inequality, we have:
\begin{equation} \nonumber
\left| \sum_{s \in \mathcal{L}} \frac{\xi_s\varphi_s+\eta_s \varphi_{-s}}{\sqrt{2\lambda_s}} \right| \leq \sum_{s \in \mathcal{L}} \frac{\vert \xi_s \vert +\vert \eta_s \vert }{\sqrt{2\lambda_s}} \leq C(\alpha) \Vert \zeta \Vert_\alpha.
\end{equation}
So we have:
\begin{align*}
\vert u(\zeta)\vert & \leq C(\alpha) \Vert \zeta \Vert_\alpha + C_\mathcal{A} |I|^{1/2} e^{|Im(\theta)|},\\
& \leq C(\alpha) \Vert \zeta \Vert_\alpha + C_\mathcal{A} \vert I-\nu\rho \vert^{1/2} e^{|Im(\theta)|} + C_\mathcal{A}(\nu\rho)^{1/2} e^{|Im(\theta)|}.
\end{align*}
We want to prove that, if $(I,\theta,\zeta) \in \mathfrak{T}_\rho(\nu,\sigma,\mu,\alpha)$, then $u(\zeta)(x) \in I_{\mathbb{C}}$ for all $x \in \mathbb{S}^1$. This is true if we assume
\begin{equation} \nonumber
C(\alpha) \nu^{1/2} \mu + C_{\mathcal{A}} \nu^{1/2 }\mu e^{\sigma} + C_{\mathcal{A}} \sqrt{2} \nu^{1/2}e^{\sigma} \leq \min (M,K).
\end{equation}
For example, if we assume that $\sigma^*=1$, we have
\begin{equation} \nonumber
\mu^* \leq \frac{\min(M,K)-C_\mathcal{A}\sqrt{2} \nu^{1/2}e}{\nu^{1/2}(C(\alpha)+C_\mathcal{A})}.
\end{equation}
The proof is thus achieved.
\endproof
Now, we are interested in the perturbation $P$. We will prove that $P$ belongs to the right class of Hamiltonian functions. Recall that
\begin{equation} \nonumber
P(\zeta)=\int_{\mathbb{S}^1} G(x,u(\zeta)(x))dx,
\end{equation}
where $g=\partial_uG$ and $g(x,u)=4u^3+O(u^4)$.
\begin{lemma} \label{Clas pertur}
Assume that $(x,u) \mapsto g(x,u)$ is real holomorphic on a neighborhood of $\mathbb{S}^1 \times J$, for $J$ some neighborhood of the origin of $\mathbb{R}$. Then for $ \alpha > 0 $, there exists $\mu^*>0$ such that, for $0< \mu \leq \mu^*$, the perturbation
\begin{align*}
P: \mathcal{O}_\mu(Y_\alpha)  &\to \mathbb{C} \\
 \zeta & \mapsto P( \zeta )
\end{align*}
belongs to $\mathcal{T}^{\alpha,1/2}(\mu)$.
\end{lemma}
\proof
Recall that for $x \in \mathbb{S}^1$, we have:
\begin{equation} \nonumber
u(x) = \sum_{s \in \mathbb{Z}}\frac{\xi_s\varphi_s(x)+\eta_s \varphi_{-s}(x)}{\sqrt{2\lambda_s}}=u(\zeta)(x).
\end{equation}
Using the Cauchy-Schwarz inequality and the fact that $\alpha>0$, there exists a constant $C_\alpha$ that depends on $\alpha$, such that for $\zeta  \in \mathcal{O}_\mu(Y_\alpha)$ we have
\begin{equation} \nonumber
\vert u(\zeta)(x) \vert \leq C_\alpha \Vert \zeta \Vert_\alpha \leq C_\alpha \mu .
\end{equation}
To prove that $P \in \mathcal{T}^{\alpha,1/2}(\mu)$, it is enough to show that
\begin{equation} \nonumber
\nabla P \in Y^\alpha \cap L_{1/2}\quad \text{ and } \quad \nabla^2P \in \mathcal{M}_{1/2}.
\end{equation}
For $\alpha\geq 0$, we define the following space:
\begin{equation} \nonumber
Z_\alpha= \left\lbrace v=\left( v_s \in \mathbb{C}, \: s \in \mathbb{Z} \right) \:| \: \left( \vert v_s \vert \langle s\rangle^{\alpha}\right)_s \in \ell^2\left( \mathbb{Z} \right)  \right\rbrace .
\end{equation}
For $v \in Z_\alpha$, we define the Fourier transform $\mathcal{F}(v)$ of $v$ by  $u(x)=\mathcal{F}(v):=~ \sum v_s e^{isx}$. We also define the discrete Sobolev space by
\begin{equation} \nonumber
H^\alpha (\mathbb{S}^1)=\left\lbrace u \,|\, u(x)=\sum_{s\in \mathbb{Z}} \hat{u}(s)e^{isx} | \: \left( \vert\hat{u}(s) \vert \langle s\rangle^{\alpha}\right)_s \in \ell^2\left( \mathbb{Z} \right) \right\rbrace .
\end{equation}
If $\alpha \in \mathbb{N}$, then
\begin{equation} \nonumber
H^\alpha (\mathbb{S}^1)=\left\lbrace u \, | \, u(x)=\sum_{s\in \mathbb{Z}} \hat{u}(s)e^{isx} | \: \left( \widehat{ \partial^\alpha u}(s)\right) _s \in \ell^2\left( \mathbb{Z} \right) \right\rbrace .
\end{equation}
So we have the following equivalence:
\begin{equation} \label{eqiv sob}
u \in H^\alpha(\mathbb{S}^1) \Longleftrightarrow \left( \hat{u}(s)\right)_s \in Z_\alpha.
\end{equation}
\begin{itemize}
\item To prove that $\nabla_\zeta P \in Y_\alpha$, it is sufficient to prove, for example, that $\frac{\partial P}{\partial \xi} \in Z_\alpha$.
We have
\begin{align*}
\frac{\partial P}{ \partial \xi_s}  (\zeta) & = \frac{1}{\sqrt{2 \lambda_s}} \int_{\mathbb{S}^1} \partial_u G \left( x, u(\zeta)(x) \right) \phi_s(x)dx.
\end{align*}
The map  $(x,u) \mapsto g(x,u)$ is real holomorphic on a neighborhood of $\: \mathbb{S}^1 \times J$, so $x \mapsto  \partial_u f \left( x, u(\zeta)(x) \right) \in H^\alpha (\mathbb{S}^1)$. We deduce from equivalence \eqref{eqiv sob} that $\frac{\partial P}{\partial \xi} \in Z_\alpha$.

\item Let us prove now that $\nabla^2P \in \mathcal{M}_{1/2}$. Recall that:
\begin{equation} \nonumber
\vert  \nabla^2 P \vert_{1/2} = \sup_{s,s' \in \mathbb{Z}} \langle s \rangle^{1/2}\langle s' \rangle^{1/2} \left\Vert \frac{\partial^2 P}{\partial\zeta_s \partial \zeta_{s'}} \right\Vert_\infty .
\end{equation}
We have
\begin{equation} \nonumber
\frac{\partial^2 P}{\partial \xi_s \xi_{s'}}= \frac{1}{2 \lambda_s^{1/2} \lambda_{s'}^{1/2}} \int_{\mathbb{S}^1} \partial_u^2 G (x,u(\zeta)(x)) \varphi_s(x) \phi_{s'}(x)dx.
\end{equation}
Then
\begin{equation} \nonumber
\frac{\partial^2 P}{\partial \zeta_s \zeta_{s'}}=  \frac{1}{2 \lambda_s^{1/2} \lambda_{s'}^{1/2}} \begin{pmatrix}
\widehat{\partial^2_u G} (s+s')& \widehat{\partial^2_u G} (s-s') \\
\widehat{\partial^2_u G} (-s+s')& \widehat{\partial^2_u G} (-s-s')
\end{pmatrix},
\end{equation}
which leads to
\begin{equation} \nonumber
\vert  \nabla^2 P \vert_{1/2} = \sup_{s \in \mathbb{Z}} \left\vert \widehat{\partial^2_u G} (s) \right\vert  < \infty.
\end{equation}
\item To conclude the proof, we have to show that $\nabla P \in L_{1/2}$. Recall that for  $\beta \leq \alpha$, we have $Y_\alpha \subset L_\beta$. So $\nabla P \in Y_1 \subset L_{1/2}$, and the proof is achieved. \endproof
\end{itemize}
Now, we are able to give the symplectic change of variable which puts the Hamiltonian \eqref{hamiltonien} into a normal form that satisfies the hypotheses of the KAM theorem.
\\Let $0< \mu < \mu^*$ and $0 < \sigma < \sigma^*$ where  $\mu^*$ and $\sigma^*$ are defined in Lemma~\ref{Cauchy_lineaire}.
\begin{theorem} \label{theo-F.N}
Let $\mathcal{A}$ be an admissible set. There exists a Borel set of zero Lebesgue measure $\mathcal{U} \subset \left[ 1,2 \right]  $, such that for any $m\in  \left(  \left[ 1,2\right]  \setminus \mathcal{U} \right) $ there exists $\nu_0$ that depends on $\mathcal{A}$, $m$ and $g$ such that:\begin{itemize}
\item[(i)] For $0<\nu \leq \nu_0$, $\alpha> 1/2$ and $\rho \in \mathcal{D}$ there exists a real symplectic holomorphic change of variable
\begin{equation} \nonumber
\Psi_\rho: \mathcal{O}^\alpha \left(   \frac{\sigma}{2} , \frac{\mu}{2} \right)  \rightarrow \mathfrak{T}_\rho(\nu,\sigma,\mu,\alpha),
\end{equation}
that transforms the symplectic form $-id\xi \wedge d\eta$ on $\mathfrak{T}_\rho(\nu,\sigma,\mu,\alpha)$ into
\begin{equation} \nonumber
 -\nu \sum_{a\in \mathcal{A}} dr_a \wedge d \theta_a - i \nu  \sum_{s\in \mathcal{L}} d\xi_s \wedge d\eta_s.
 \end{equation}
\item[(ii)]For $c> \frac{1}{2}$, the change of variable $\Phi_\rho$ can be extended holomorphically on the following complex domain
\begin{equation} \nonumber
\mathcal{D}_c = \left\lbrace \rho \in\mathbb{C}^{\mathcal{A}}  \mid  \left| \rho_j - \frac{3}{2} \right| \leq c , \; 1 \leq j \leq \operatorname{Card}( \mathcal{A})   \right\rbrace .
\end{equation}
It transforms the perturbed Hamiltonian $H=H_2+P$ into the following normal form
\begin{equation} \label{Final FN}
\nu^{-1}H \circ \Psi_\rho=\Omega(\rho)\cdot r + \sum_{a\in\mathcal{L}}\Lambda_a ( \rho) \xi_a\eta_a + f(r, \theta, \zeta, \rho),
\end{equation}
 for all $\rho \in \mathcal{D}_c $. The internal frequency vector $\Omega$ and the external frequencies $\Lambda_a$, $a\in \mathcal{L}$, are given explicitly in  \eqref{Om} and \eqref{Lam}. Moreover, they are linear in $\rho$ and verify:
\begin{equation} \label{new_frequencies}
\vert \Omega (\rho) - \omega (\rho) \vert \leq C \nu , \quad \vert \Lambda_a(\rho) - \lambda_a (\rho) \vert \leq C \nu \vert a \vert^{-1},
\end{equation}
 for all $\rho \in \mathcal{D}_c $. The constant $C$ depends on the admissible set $\mathcal{A}$.
\item[(iii)] The perturbation $f$ is real holomorphic on $ \mathcal{D}_c $, belongs to $\mathcal{T}^{\alpha,1/2}(\mathcal{D},\frac{\sigma}{2},\frac{\mu}{2})$ and satisfies
\begin{align*}
\lc f \rc^{\alpha,1/2,\gamma} _{\frac\sigma2,\frac\mu2,\mathcal{D}} & \leq C_1 (1+\gamma) \nu\mu^4, \\
\lc f^T \rc^{\alpha,1/2,\gamma} _{\frac\sigma2,\frac\mu2,\mathcal{D}} & \leq C_1 (1+\gamma) \nu^{3/2}\mu^5.
\end{align*}
The constant $C_1$ depends on the admissible set $\mathcal{A}$, the mass $m$ and the nonlinearity $g$.
\end{itemize}
\end{theorem}
\begin{remark}
In $(iii)$, we need to estimate the derivative of the perturbation $f$ with respect to the parameter $\rho$. This is possible thanks to the Cauchy estimate. So we need to extend $\Psi_\rho$ holomorphically onto $\mathcal{D}_c$ (see \ref{Re-scaling des variables}).
\end{remark}
All the rest of this section will be dedicated to the proof of the previous theorem.
\subsection{Resonance}
We consider $H_4$, the quartic part of the Hamiltonian $H$ defined by:
\begin{equation} \nonumber
H_4=H_2+P_4,
\end{equation}
where
\begin{equation} \nonumber
H_2 =  \underset{s \in \mathbb{Z}}{\sum} \lambda_s \xi_s \eta_s,
\end{equation}
and
\begin{equation} \nonumber
P_4 =\underset{(i,j,k,l) \in \mathbb{Z}}{\sum}C(i,j,k,l)\frac{(\xi_i+\eta_{-i})(\xi_j+\eta_{-j})(\xi_k+\eta_{-k})(\xi_l+\eta_{-l})}
{4\sqrt{\lambda_i\lambda_j\lambda_k\lambda_l}}.
\end{equation}
The constant $C(i,j,k,l) $ is defined by:
\begin{equation*}
C(i,j,k,l):=\int_{\mathbb{S}^1} \varphi_i(x)\varphi_j (x)  \varphi_k(x) \varphi_l(x)dx= \left\{
\begin{array}{r c l}
\frac{1}{2\pi} &\text{if}\:\:i+j+k+l=0,\\
0&\text{if}\:\:i+j+k+l\neq0.
\end{array}
\right.
\end{equation*}
We define the following subset of $\mathbb{Z}^4$:
\begin{equation} \nonumber
\mathcal{J}:= \lbrace (i,j,k,l) \in \mathbb{Z}^4 | i+j=k+l \rbrace.
\end{equation}
So
\begin{equation} \nonumber
P_4 =\frac{1}{2\pi}\underset{(i,j,k,l) \in \mathcal{J}}{\sum}\frac{(\xi_i+\eta_{-i})(\xi_j+\eta_{-j})(\xi_k+\eta_{-k})(\xi_l+\eta_{-l})}
{4\sqrt{\lambda_i\lambda_j\lambda_k\lambda_l}}.
\end{equation}
We can decompose $P_4$ in three parts $P_4=P_4^0+P_4^1+P_4^2$ where:
\begin{align*}
P_4^0 &= \frac{1}{8\pi} \underset{(i,j,k,l) \in \mathcal{J}}{\sum} \frac{\xi_i\xi_j\xi_k\xi_l+\eta_i\eta_j\eta_k\eta_l}{\sqrt{\lambda_i\lambda_j\lambda_k\lambda_l}},\\
P_4^1 &= \frac{1}{2\pi} \underset{(i,j,k,-l) \in \mathcal{J}}{\sum} \frac{\xi_i\xi_j\xi_k\eta_l+\eta_i\eta_j\eta_k\xi_l}{\sqrt{\lambda_i\lambda_j\lambda_k\lambda_l}},\\
P_4^2 &= \frac{3}{4\pi}\underset{(i,j,-k,-l) \in \mathcal{J}}{\sum}\frac{\xi_i\xi_j\eta_k\eta_l}{\sqrt{\lambda_i\lambda_j\lambda_k\lambda_l}}.
\end{align*}
For $(i,j,k,l) \in \mathbb{Z}^4$, we define the small divisors:
\begin{equation*}
\Omega_0(i,j,k,l) = \lambda_i+\lambda_j+\lambda_k+\lambda_l; \hspace{0.5cm}
\Omega_1(i,j,k,l) = \lambda_i+\lambda_j+\lambda_k-\lambda_l; \hspace{0.5cm}
\Omega_2(i,j,k,l) = \lambda_i+\lambda_j-\lambda_k-\lambda_l.
\end{equation*}
\begin{definition}
A monomial $\xi_i\xi_j\xi_k\eta_l$ or $\eta_i\eta_j\eta_k\xi_l$ is resonant if $\Omega_1(i,j,k,l)=0$. In this case, we denote $\mathcal{R}_1:=\lbrace (i,j,k,l) \in \mathbb{Z}^4\: | \: \Omega_1(i,j,k,l)=0 \rbrace$. A monomial $\xi_i\xi_j\eta_k\eta_l$ is resonant if $\Omega_2(i,j,k,l)=0$. In this case, we denote $\mathcal{R}_2:=\lbrace (i,j,k,l) \in \mathbb{Z}^4 \:|\: \Omega_2(i,j,k,l)=0 \rbrace$. Let $\mathcal{R}$ be the union of $\mathcal{R}_1$ and $\mathcal{R}_2$.
\end{definition}
We define also
\begin{equation} \nonumber
\mathcal{J}_2= \lbrace (i,j,k,l) \in \mathcal{J}\: |\: \#\lbrace i,j,k,l \rbrace \cap \mathcal{A} \geq 2 \rbrace,
\end{equation}
and
\begin{equation} \nonumber
\mathcal{J}_2^c=\mathcal{J} \setminus \mathcal{J}_2 = \lbrace (i,j,k,l) \in \mathcal{J}\: |\: \#\lbrace i,j,k,l \rbrace \cap \mathcal{L} \geq 3 \rbrace
\end{equation}
\begin{lemma} \label{lem_ens}
There exists a Borel set $ \mathcal{U} \subset \left[ 1,2 \right] $ of full Lebesgue measure, such that for $m \in \mathcal{U} $ we have:
\begin{itemize}
\item[(i)] $\mathcal{R}_1=\emptyset$,
\item[(ii)] $\mathcal{R}_2 \subset \lbrace (i,j,k,l) \in \mathbb{Z}^4\: | \: \lbrace \vert i\vert, \vert j \vert \rbrace =\lbrace \vert k\vert, \vert l \vert \rbrace \rbrace$,
\item[(iii)] There exists $\gamma(m)>0$ such that for any $(i,j,k,l) \in \mathcal{J}_2 \setminus \mathcal{R}$, we have:
\begin{equation} \nonumber
\vert \Omega_1(i,j,k,l) \vert \geq \gamma(m); \hspace{1cm}\vert \Omega_2(i,j,k,l) \vert \geq \gamma(m).
\end{equation}
\end{itemize}
\end{lemma}
\proof We begin by proving assertions (i) and (ii). Let us fix $(i,j,k,l)\in \mathbb{Z}^4$. Consider, for $\delta=\pm1$, the function:
\begin{equation} \nonumber
f(m)= \sqrt{i^2+m}+\sqrt{j^2+m}+\delta\sqrt{k^2+m}-\sqrt{l^2+m}.
\end{equation}
The  function $f$ is analytic on $\left[ 1,2 \right] $, and it can be extended to an algebraical multi-valued function of  $m \in \mathbb{C}$. If $j=k=l=0$, then $-i^2$ is a branching point of $f$. Branching points for $f$ are $-i^2,-j^2,-k^2$ and $-l^2$.

If $\delta=1$, then $f$ is not identically zero on $\left[ 1,2 \right] $, and admits finitely many zeros. It follows that, there is a Borel set $\mathcal{U}_1 \subset \left[ 1,2 \right]$ of full Lebesgue measure, such that for $m$ $\in \mathcal{U}_1$, we have $\mathcal{R}_1=\emptyset$.

Now, for $\delta=-1$, if $f\mid_{\left[ 1,2 \right]}\equiv0$, then $f\equiv0$, and we have $\lbrace i^2,j^2 \rbrace=\lbrace k^2,l^2 \rbrace$. If $f$ is not identically zero on $\left[ 1,2 \right]$, then $f$ admits finitely many zeros. So there is a Borel set $\mathcal{U}_2 \subset \left[ 1,2 \right]$ of full Lebesgue measure, such that for any $m$ $\in \mathcal{U}_2$, we have $\mathcal{R}_2 \subset \lbrace (i,j,k,l) \in \mathbb{Z}^4\: | \: \lbrace \vert i\vert, \vert j \vert \rbrace =\lbrace \vert k\vert, \vert l \vert \rbrace \rbrace$.

It remains to prove the last assertion. Let us begin with $\Omega_2$. Thanks to Proposition~\ref{prop-D3}, for $\kappa>0$, there is an open set $\mathcal{C}_\kappa \subset \left[ 1,2 \right]$ such that
\begin{equation} \nonumber
mes (\mathcal{C}_\kappa)<C\kappa^\tau,
\end{equation}
where $\tau>0$ and depends on $\operatorname{Card}( \mathcal{A})$. The constant $C$ depends on the admissible set $\mathcal{A}$. For $m \in \left( \left[ 1,2 \right] \setminus\mathcal{C}_\kappa \right)  $ and any $(i,j,k,l) \in \mathcal{J}_2\setminus\mathcal{R}$, we have:
\begin{equation} \nonumber
\vert\Omega_2(i,j,k,l)\vert> \kappa.
\end{equation}
If $ \kappa' \leq \kappa$, then $\mathcal{C}_{\kappa'} \subset \mathcal{C}_\kappa$. So $\mathcal{C}:=\underset{0<\kappa<1}{ \cap}\mathcal{C}_\kappa$ is a Borel set, and we have
\begin{equation} \nonumber
mes( \mathcal{C)}=0.
\end{equation}
Moreover, for $m \in \mathcal{U}_3\equiv \left( \left[ 1,2 \right]\setminus\mathcal{C} \right) =\underset{0<\kappa<1}{ \cup}(\left[ 1,2 \right] \setminus\mathcal{C}_\kappa)$, there exists a constant $\gamma(m)$ such that for each $(i,j,k,l) \in \mathcal{J}_2\setminus\mathcal{R}$, we have:
\begin{equation} \nonumber
\vert \Omega_2(i,j,k,l)\vert>\gamma(m).
\end{equation}
To control $\Omega_1(i,j,k,l)$ we follow the same procedure, but we use Proposition~\ref{final-petit diviseur} instead of Proposition~\ref{prop-D3}. Finally we denote $\mathcal{U}=\mathcal{U}_1\cap\mathcal{U}_2\cap\mathcal{U}_3\cap\mathcal{U}_4$, where $\mathcal{U}_4$ is the Borel set of full Lebesgue measure that we obtain after controlling $\Omega_1(i,j,k,l)$.
\endproof
\subsection{Birkhoff's Precedure}
For $\alpha>0$, we recall the definition of the following space:
\begin{equation} \nonumber
Z_\alpha= \left\lbrace v=\left( v_s \in \mathbb{C}, \: s \in \mathbb{Z} \right) \:| \: \left( \vert v_s \vert \langle s\rangle^{\alpha}\right)_s \in \ell^2\left( \mathbb{Z} \right)  \right\rbrace .
\end{equation}
We endow $Z_\alpha$ with the norm:
\begin{equation} \nonumber
\Vert v \Vert_\alpha^2 = \sum_{s \in \mathbb{Z}} \vert v_s \vert^2 \langle s\rangle^{2\alpha}, \quad \langle s\rangle = max(\vert s \vert , 1 ).
\end{equation}
We denote by $v * y$ the convolution in $\ell^2(\mathbb{Z})$ defined by $(v * w)_l= \sum_{i+j=l} v_iw_j$.
We recall Lemma 2 from~\cite{kuksin1996invariant}.
\begin{lemma} \label{convolution}
Consider $v ,w \in Z_\alpha$ for $\alpha > \frac{1}{2}$. Then, $v*w \in Z_\alpha$, and
\begin{equation}
\| v * w \|_\alpha < C(\alpha)  \| v\|_\alpha  \| w\|_\alpha,
\end{equation}
where $C$ is a constant that depends only on $\alpha$.
\end{lemma}
\proof Consider $v,w \in Z_\alpha$. We have:
\begin{align*}
\| v*w \|^2_\alpha & = \underset{s \in \mathbb{Z}}{ \sum} \langle s \rangle^{2\alpha} \left| \underset{i+j=s}{ \sum} v_iw_j \right|^2 \\
& = \sum_{s \in \mathbb{Z}} \langle s \rangle^{2\alpha} \left| \sum_{i+j=s} \frac{\langle s \rangle^{\alpha}}{\langle i \rangle^{\alpha}\langle j \rangle^{\alpha}} \frac{\langle i \rangle^{\alpha}\langle j \rangle^{\alpha}}{\langle s \rangle^{\alpha}} v_iw_j \right|^2 \\
& \leq \sum_{s \in \mathbb{Z}} \langle s \rangle^{2\alpha} \left( \sum_{i+j=s} \left(  \frac{\langle s \rangle}{\langle i \rangle\langle j \rangle} \right)^{2 \alpha}   \right) \left(  \sum_{i+j=s} \frac{\langle i \rangle^{2\alpha}\langle j \rangle^{2\alpha}}{\langle s \rangle^{2\alpha}} |v_i|^2 |w_j|^2 \right),
\end{align*}
and
\begin{equation} \nonumber
\sum_{i+j=s} \left(  \frac{\langle s \rangle}{\langle i \rangle\langle j \rangle} \right)^{2 \alpha}  \leq \sum_{i,j \in \mathbb{Z}} \left(     \frac{\langle i \rangle + \langle j \rangle}{\langle i \rangle\langle j \rangle} \right)^{2\alpha} \leq 4^\alpha \sum_{i,j \in \mathbb{Z}}  \frac{1}{\langle i \rangle^{2\alpha}} +  \frac{1}{\langle j \rangle^{2\alpha}} \leq C^2(\alpha).
\end{equation}
Then
\begin{align*}
\Vert v * w \Vert_\alpha^2 & \leq  C(\alpha)^2 \sum_{s\in \mathbb{Z}} \sum_{i+j = s}  \langle i \rangle^{2\alpha} |v_i|^2 \langle j \rangle^{2\alpha} |w_j|^2 \\
& \leq \sum_{i,j \in \mathbb{Z}} \langle i \rangle^{2\alpha} |v_i|^2 \langle j \rangle^{2\alpha} |w_j|^2\\
& = C(\alpha)^2 \Vert v \Vert_\alpha^2 \Vert w \Vert_\alpha^2.
\end{align*}
The proof is thus achieved.
\endproof
We endow the phase space with symplectic structure $-i\sum d\xi_k \wedge d\eta_k$. For $\alpha > 1/2$, $Y_\alpha$ is an algebra for the convolution.
\begin{lemma}  \label{XPanalytique}
Let $\alpha > 1/2$ and $P^4$ a real homogeneous polynomial on $Y_\alpha$, of degree $4$, indexed by $\mathcal{J}$. We assume that $P^4$ is of the form:
\begin{equation} \nonumber
P^4(\zeta)= \underset{(j_1,j_2,j_3,j_4) \in \mathcal{J}}{ \sum} \: \underset{1\leq r \leq 4}{ \sum} a^r_{j_1,j_2,j_3,j_4} \xi_{j_1}\ldots \xi_{j_r}\eta_{j_{r+1}} \ldots \eta_{j_4},
\end{equation}
where $|a^r_{j_1,j_2,j_3,j_4}|<M$ for any $(j_1,j_2,j_3,j_4) \in \mathcal{J}$. Then we have:
\begin{equation} \nonumber
\Vert \nabla P^4 \Vert_\alpha \leq C(\alpha,M) \Vert \zeta \Vert_\alpha^3.
\end{equation}
So, for $t\leq 1$, the flow  $\Phi^t_{P^4}$ of the Hamiltonian vector $X_{P^4}=i J \nabla P^{4}$ is well-defined, real and analytic on the ball:\begin{equation} \nonumber
\mathcal{O}_\delta(Y_\alpha)=\left\lbrace \zeta \in Y_\alpha | \Vert \zeta \Vert_\alpha < \delta=\delta(M) \right\rbrace.
\end{equation}
Moreover, for any $\zeta \in \mathcal{O}_\delta(Y_\alpha)$,
\begin{equation} \nonumber
\Vert \Phi^t_{P_4}(\zeta)-\zeta \Vert_\alpha \leq C(M) \Vert \zeta \Vert_\alpha^3.
\end{equation}
\end{lemma}
\proof
We recall that $X_{P^4}= i \begin{pmatrix}
\nabla_\xi P^4 \\
\nabla_\eta P^4
\end{pmatrix}$. Since $|a^r_{j_1,j_2,j_3,j_4}| < M$, we have:
\begin{equation} \nonumber
\left| \dfrac{\partial P}{\partial \eta_l} \right| \leq M \underset{(i,j,k,l) \in \mathcal{J}}{ \sum} | \xi_i\xi_j\xi_k | + |\xi_i\xi_j \eta_k | +|\xi_i\eta_j\eta_k |+ |\eta_i\eta_j\eta_k|.
\end{equation}
We remark that
\begin{equation} \nonumber
\underset{(i,j,k,l) \in \mathcal{J}}{ \sum} | \xi_i\xi_j\xi_k |=(\xi\star\xi\star\xi)_l+(\xi\star\xi\star\tilde{\xi})_l+(\xi\star\tilde{\xi}\star\tilde{\xi})_l+(\tilde{\xi}\star\tilde{\xi}\star\tilde{\xi})_l+(\xi\star\tilde{\xi}\star\xi)_l+(\tilde{\xi}\star\tilde{\xi}\star\xi)_l,
\end{equation}
where  $\tilde{\xi}=(\tilde{\xi}_j)_{j\in \mathbb{Z}}$ and $\tilde{\xi}_j=\xi_{-j}$. Using Lemma~\ref{convolution}, we have
\begin{align*}
\left\Vert \frac{\partial P^4}{\partial \eta} \right\Vert^2_\alpha & = \sum_{l \in \mathbb{Z}} \langle l \rangle^{2\alpha} \left\vert \frac{\partial P^4}{\partial \eta_l} \right\vert^2  \\
& \leq C(M) \left(  \Vert \xi \star \xi \star \xi \Vert_\alpha^2 + \Vert \xi \star \xi \star \eta \Vert_\alpha^2 + \Vert \xi \star \eta \star \eta \Vert_\alpha^2 + \Vert \eta \star \eta \star \eta \Vert_\alpha^2   \right)  \\
& \leq C(\alpha,M) \left( \Vert \xi \Vert^6_\alpha + \Vert \xi \Vert^4_\alpha \Vert \eta \Vert^2_\alpha + \Vert \xi \Vert^2_\alpha \Vert \eta \Vert^4_\alpha + \Vert \eta \Vert^6_\alpha  \right)  \\
& \leq C(\alpha, M) \left(  \Vert \xi \Vert^2_\alpha + \Vert \eta \Vert^2_\alpha \right)^3 = C(\alpha, M ) \Vert \zeta \Vert_\alpha^6.
\end{align*}
We prove the same way that:
\begin{equation} \nonumber
\left\Vert \frac{\partial P^4}{\partial \xi} \right\Vert^2_\alpha \leq C(\alpha, M ) \Vert \zeta \Vert_\alpha^6.
\end{equation}
So, we have:
\begin{equation} \nonumber
\left\Vert \nabla_\zeta P^4 \right\Vert_\alpha = \left(  \left\Vert \frac{\partial P^4}{\partial \xi} \right\Vert^2_\alpha + \left\Vert \frac{\partial P^4}{\partial \eta} \right\Vert^2_\alpha  \right) ^{1/2} \leq C(\alpha,M)\Vert \zeta \Vert_\alpha^3.
\end{equation}
This concludes the proof of the lemma.
\endproof
\begin{lemma} \label{Hess poly}
Consider $ D^-$ the following bounded operator on $Y_\alpha$
\begin{equation} \nonumber
D^-= \diag \left\lbrace \lambda_s^{-1/2} I_2, \, s \in \mathbb{Z}  \right\rbrace.
\end{equation}
Define $Q^4(\zeta):=P^4(D^-(\zeta))$, where $P^4$ is the polynomial defined in the Lemma~\ref{XPanalytique}. Then $\nabla^2_\zeta Q^4 \in \mathcal{M}_{1/2}$, and we have:
\begin{equation} \nonumber
\left\vert \nabla^2_\zeta Q^4 \right\vert_{1/2} \leq C( \alpha , M) \Vert \zeta \Vert_\alpha^2.
\end{equation}
\end{lemma}
\proof
We recall that
\begin{equation} \nonumber
P^4(\zeta)= \underset{(j_1,j_2,j_3,j_4) \in \mathcal{J}}{ \sum} \: \underset{1\leq r \leq 4}{ \sum} a^r_{j_1,j_2,j_3,j_4} \xi_{j_1}\ldots \xi_{j_r}\eta_{j_{r+1}} \ldots \eta_{j_4},
\end{equation}
where $|a^r_{j_1,j_2,j_3,j_4}|<M$ for $(j_1,j_2,j_3,j_4) \in \mathcal{J}$. So we have
\begin{equation} \nonumber
Q^4(\zeta)= \underset{(j_1,j_2,j_3,j_4) \in \mathcal{J}}{ \sum} \: \underset{1\leq r \leq 4}{ \sum} \frac{a^r_{j_1,j_2,j_3,j_4}}{\sqrt{\lambda_{j_1}\lambda_{j_2}\lambda_{j_3}\lambda_{j_4}}}  \xi_{j_1}\ldots \xi_{j_r}\eta_{j_{r+1}} \ldots \eta_{j_4},
\end{equation}
and
\begin{equation} \nonumber
\vert \nabla_\zeta^2 Q^4 \vert_{1/2} = \sup_{s,s' \in \mathbb{Z}} \langle s  \rangle^{1/2} \langle s'  \rangle^{1/2} \left\Vert \frac{\partial^2 Q^4}{\partial \zeta_s \partial \zeta_{s'}} \right\Vert_\infty.
\end{equation}
For any $s$ and $s'$ in $\mathbb{Z}$, we have:
\begin{equation} \nonumber
 \langle s  \rangle^{1/2} \langle s'  \rangle^{1/2} \left\Vert \frac{\partial^2 Q^4}{\partial \zeta_s \partial \zeta_{s'}} \right\Vert_\infty \leq \Vert A \Vert_\infty,
\end{equation}
where $A$ is a real square matrix of size 2, whose coefficients are homogeneous polynomials of $Y_\alpha$ of degree~2, of the form:
\begin{equation} \nonumber
P^2(\zeta)= \sum_{\bar{j}\in \mathcal{J}} \sum_{k,l=1}^4 \frac{1}{\sqrt{\lambda_{j_{\sigma(k)}}\lambda_{j_{\sigma(l)}}}} \left( a_1 \xi_{j_{\sigma(k)}}\xi_{j_{\sigma(l)}} + a_2 \xi_{j_{\sigma(k)}}\eta_{j_{\sigma(l)}} + a_3 \eta_{j_{\sigma(k)}}\eta_{j_{\sigma(l)}}  \right),
\end{equation}
where $\bar{j}=(j_{\sigma(1)},j_{\sigma(2)},j_{\sigma(3)},j_{\sigma(4)})$ and $\sigma$ is a permutation from the symmetric group $ S_4$. Using Lemma~\ref{convolution}, we have:
\begin{equation} \nonumber
\left\vert P^2(\zeta) \right\vert \leq C(\alpha,M) \Vert  \zeta \Vert_\alpha^2 .
\end{equation}
So
\begin{equation} \nonumber
\vert \nabla^2_\zeta Q \vert_{1/2} \leq C(\alpha,M)\Vert  \zeta \Vert_\alpha^2.
\end{equation}
\endproof
\begin{remark}
We recall that $P_4=P_4^0+P_4^1+P_4^2$ where
\begin{align*}
P_4^0 &= \frac{1}{8\pi} \underset{(i,j,k,l) \in \mathcal{J}}{\sum} \frac{\xi_i\xi_j\xi_k\xi_l+\eta_i\eta_j\eta_k\eta_l}{\sqrt{\lambda_i\lambda_j\lambda_k\lambda_l}},\\
P_4^1 &= \frac{1}{2\pi} \underset{(i,j,k,-l) \in \mathcal{J}}{\sum} \frac{\xi_i\xi_j\xi_k\eta_l+\eta_i\eta_j\eta_k\xi_l}{\sqrt{\lambda_i\lambda_j\lambda_k\lambda_l}},\\
P_4^2 &= \frac{3}{4\pi}\underset{(i,j,-k,-l) \in \mathcal{J}}{\sum}\frac{\xi_i\xi_j\eta_k\eta_l}{\sqrt{\lambda_i\lambda_j\lambda_k\lambda_l}}.
\end{align*}
The coefficients of each monomial are bounded by $3/4\pi$. Using \ref{XPanalytique} and \ref{Hess poly}, we have $P_4 \in \mathcal{T}^{\alpha,1/2}(\mu)$, for $ \alpha > 1/2 $ and $\zeta  \in \mathcal{O}_\mu(Y_\alpha)$.
\end{remark}

Let $ \mathcal{U} \subset [1,2] $ be the Borel full Lebesgue measure set from Lemma~\ref {lem_ens}. For $ m \in \mathcal{U} $, we want to construct a holomorphic real symplectic change of variable in the neighborhood of the origin of $ Y_\alpha $ which transforms the quartic part  of the Hamiltonian $ H $ into a Birkhoff normal form up to order 5. This transformation extracts the integrable terms from the quartic part of the perturbation $P$ and cubic terms in the direction of $ \mathcal{L} = \mathbb{Z} \setminus \mathcal{A} $.
\begin{proposition} \label{Bir.Nor.For}
For $m \in \mathcal{U}$, there is a holomorphic real symplectic change of variable $\tau$ on $\mathcal{O}_{\delta(m)}(Y_\alpha)$, for some $\delta(m)>0$ and $ \alpha > 1/2$.  The change of variable $\tau$ satisfies:
\begin{equation} \label{chang.birk.id}
\Vert \tau(\zeta)-\zeta \Vert_\alpha \leq C(m) \Vert \zeta \Vert_\alpha^3,\quad \forall \zeta \in  \mathcal{O}_{\delta(m)}(Y_\alpha).
\end{equation}
The mapping $ \tau $ tranforms the Hamiltonian $H=H_2+P=H_2+P_4+R_5$ into :
\begin{equation} \label{Birkhoff-FN}
H \circ \tau = (H_2+P) \circ \tau= H_2 + Z_4+Q_4+R_6+R_5\circ \tau,
\end{equation}
where
\begin{equation} \nonumber
Z_4 = \frac{3}{\pi} \underset{(i,j,k,l) \in \mathcal{J}_2\cap\mathcal{R}_2}{\sum} \frac{\xi_i\xi_j\eta_k\eta_l}{\lambda_i\lambda_j},
\end{equation}
and $Q_4=Q_{4,1}+Q_{4,2}$ for
\begin{align*}
Q_{4,1} & = \frac{1}{2\pi} \underset{(i,j,-k,l) \in \mathcal{J}_2^c}{\sum} \frac{\xi_i\xi_j\xi_k\eta_l+\eta_i\eta_j\eta_k\xi_l}{\sqrt{\lambda_i\lambda_j\lambda_k\lambda_l}}, \\
Q_{4,2} & = \frac{3}{4\pi} \underset{(i,j,k,l) \in \mathcal{J}_2^c}{\sum} \frac{\xi_i\xi_j\eta_k\eta_l}{\sqrt{\lambda_i\lambda_j\lambda_k\lambda_l}}.
\end{align*}
The polynomial $Z_4$ contains integrable terms while $Q_4$ is cubic or quartic in the direction of $\mathcal{L}$. Moreover, $Z_4$, $Q_4$, $R_6$  and $R_5 \circ \tau$ are real holomorphic on $\mathcal{O}_{\delta(m)}(Y_\alpha)$. The remainder terms $R_5 \circ \tau$ and $R_6 $ are respectively of order 5 and 6 at the origin of $Y_\alpha$. Moreover, for any $0 < \mu \leq \delta (m) $,  $Z_4$, $Q_4$, $R_5 \circ \tau$ and $R_6 $ belong to $\mathcal{T}^{\alpha,1/2}(\mu)$ and satisfy:
\begin{equation} \label{estim-Z-Q}
\lc Z_4 \rc^{\alpha,1/2} _{\mu} + \lc Q_4 \rc^{\alpha,1/2} _{\mu}  \leq C \mu^4,
\end{equation}
\begin{equation} \label{estim-R-6}
\lc R_6 \rc^{\alpha,1/2} _{\mu} \leq C \mu^6,
\end{equation}
\begin{equation} \label{estim-R-5}
\lc R_5\circ \tau \rc^{\alpha,1/2} _{\mu} \leq C \mu^5,
\end{equation}
where the constant $C$ depends on $m$, the nonlinearity $g$ and the admissible set $\mathcal{A}$.
\end{proposition}

We recall that the Poisson bracket associated to the symplectic form$-i\underset{s \in \mathbb{Z}}{ \sum} d\xi_s \wedge d\eta_s$  is
\begin{equation} \nonumber
\lbrace f,g \rbrace (\xi,\eta)=i\underset{j \in \mathbb{Z}}{\sum} \frac{\partial f}{\partial \eta_j}\frac{\partial g}{\partial \xi_j}-\frac{\partial f}{\partial \xi_j}\frac{\partial g}{\partial \eta_j} ,
\end{equation}
for $f,g \in \mathcal{C}^1(Z_\alpha \times Z_\alpha)$.
\begin{lemma} \label{calcul}
Let $P$ be a homogeneous polynomial of degree 4 defined by:
\begin{equation} \nonumber
P(\xi,\eta) = \sum _{\vert \alpha\vert + \vert \beta \vert =4} = a_{\alpha , \beta} \xi^\alpha \eta^\beta,
\end{equation}
where $\xi^\alpha=\xi^{\alpha_1}_1\xi^{\alpha_2}_2\xi^{\alpha_3}_3\xi^{\alpha_4}_4$. Then:
\begin{equation} \nonumber
\left\lbrace H_2, P \right\rbrace(\xi,\eta)= i \sum _{\vert \alpha\vert + \vert \beta \vert =4}  a_{\alpha , \beta} \Omega_{\min (\vert \alpha\vert , \vert \beta \vert)} (\alpha,\beta)  \xi^\alpha \eta^\beta.
\end{equation}
\end{lemma}
\proof
We prove the previous lemma by using the expression of the Hamiltonian $H_2$, the Poisson bracket and the frequencies $\Omega_p(\alpha,\beta)$ for $0\leq p\leq 2$.
\endproof
\proof[Proof of the Proposition~\ref{Bir.Nor.For}]
We want to construct a holomorphic real symplectic change of variable $ \tau$ in the neighborhood of the origin of $Y_\alpha$ for $\alpha>1/2$. The mapping $\tau$ puts the Hamiltonian $H$ into a Birkhoff normal form up to order 5. To do this, we use a classical method: $\tau$ will be the time one flow of a Hamiltonian $\chi_4$ (ie $\tau=\Phi_{\chi_4}^1$ where $\Phi_{\chi_4}^t$ is the flow of $\chi_4$ at time $t$). The Hamiltonian $ \chi_4$  will be a solution of a certain homological equation.
Using the Taylor formula, we obtain:
\begin{align*}
 \left( H_2+P_4+R_5  \right) \circ \tau & = \left( H_2+P_4 \right) \circ \tau + R_5  \circ \tau \\
 & = H_2+ P_4  + \lbrace H_2, \chi_4  \rbrace + \lbrace P_4,\chi_4 \rbrace  \\
 & + \int_0^1 (1-t) \left\lbrace \left\lbrace H_2+P_4,\chi_4 \right\rbrace ,\chi_4 \right\rbrace \circ \Phi_{\chi_4}^tdt + R_5  \circ \tau.
\end{align*}
We want that
\begin{equation} \nonumber
\left( H_2+P_4+R_5  \right) \circ \tau = H_2+Z_4+Q_4+R_6+R_5 \circ \tau.
\end{equation}
So, by taking
\begin{equation} \nonumber
R_6=\left\lbrace P_4,\chi_4 \right\rbrace + \int_0^1 (1-t) \left\lbrace \left\lbrace H_2+P_4,\chi_4 \right\rbrace ,\chi_4 \right\rbrace \circ \Phi_{\chi_4}^tdt,
\end{equation}
the Hamiltonian $\chi_4$ satisfies the following homological equation:
\begin{equation} \label{eq-homo}
\left\lbrace H_2,\chi_4 \right\rbrace = Z_4+Q_4-P_4.
\end{equation}
Using Lemma~\ref{calcul}, the Hamiltonian $\chi_4$ is given by:
\begin{align*}
\chi_4 & =   \frac{i}{8\pi}  \underset{(i,j,-k,-l) \in \mathcal{J}}{\sum} \frac{\xi_i\xi_j\xi_k\xi_l-\eta_i\eta_j\eta_k\eta_l}{\Omega_0(i,j,k,l)\sqrt{\lambda_i\lambda_j\lambda_k\lambda_l}} \\
& + \frac{i}{2\pi} \underset{(i,j,-k,l) \in \mathcal{J}_2}{\sum} \frac{\xi_i\xi_j\xi_k\eta_l-\eta_i\eta_j\eta_k\xi_l}{\Omega_1(i,j,k,l)\sqrt{\lambda_i\lambda_j\lambda_k\lambda_l}}\\
& + \frac{3i}{4\pi} \underset{(i,j,k,l) \in \mathcal{J}_2\setminus\mathcal{R}_2}{\sum} \frac{\xi_i\xi_j\eta_k\eta_l}{\Omega_2(i,j,k,l)\sqrt{\lambda_i\lambda_j\lambda_k\lambda_l}}.
\end{align*}
By Lemma~\ref{lem_ens}, there is a Borel set $\mathcal{U} \subset [1,2]$ of full Lebesgue  measure set, such that for $m\in \mathcal{U}$, there is a constant $\gamma(m)>0$ smaller than  $|\Omega_1(i,j,k,l)|$ and $|\Omega_2(i,j,k,l)|$. Remark also that $\Omega_0(i,j,k,l)>4$. Then $\chi_4$ is a homogeneous polynomial of degree 4 indexed by $J$, with bounded coefficients. So, the Hamiltonian vector field $X_{\chi_4}$ is real and holomorphic on  $Y_\alpha$. Using Lemma~\ref{XPanalytique}, for $m \in \mathcal{U}$  and $\alpha >1/2$, there exists $C(\alpha,m)>0$ such that:
\begin{equation} \nonumber
\Vert X \chi_4 \Vert_\alpha \leq C(\alpha,m) \Vert \zeta \Vert_\alpha^3.
\end{equation}
So, there is $\delta(m)>0$ such that $\tau$ is real holomorphic symplectic change of variable on $\mathcal{O}_{\delta(m)}(Y\alpha)$.

By Lemma~\ref{calcul}, we have
\begin{align*}
\left\lbrace H_2,\chi_4\right\rbrace   = & - \frac{1}{8\pi}  \underset{(i,j,-k,-l) \in \mathcal{J}}{\sum} \frac{\xi_i\xi_j\xi_k\xi_l+\eta_i\eta_j\eta_k\eta_l}{\sqrt{\lambda_i\lambda_j\lambda_k\lambda_l}} \\
& - \frac{1}{2\pi} \underset{(i,j,-k,l) \in \mathcal{J}_2}{\sum} \frac{\xi_i\xi_j\xi_k\eta_l+\eta_i\eta_j\eta_k\xi_l}{\sqrt{\lambda_i\lambda_j\lambda_k\lambda_l}}\\
& - \frac{3}{4\pi} \underset{(i,j,k,l) \in \mathcal{J}_2\setminus\mathcal{R}_2}{\sum} \frac{\xi_i\xi_j\eta_k\eta_l}{\sqrt{\lambda_i\lambda_j\lambda_k\lambda_l}}.
\end{align*}
So
\begin{align*}
(H_2+P_4) \circ \tau & =  H_2 + Z_4 + Q_4 + R_6,
\end{align*}
where $Z_4$, $Q_4$ is defined as in the proposition.  They are two homogeneous polynomials of degree 4 with bounded coefficients. From Lemma \ref{convolution} and \ref{XPanalytique}, for $\alpha > 1/2$ and $0< \mu \leq \delta(m)$,  these two polynomials belong to $\mathcal{T}^{\alpha,1/2}(\mu)$ and satisfy \ref{estim-Z-Q}.

Let us study the remainder terms $R_6$ and $R_5 \circ \tau$. Concerning  $R_6$, by construction, $R_6$ is a holomorphic Hamiltonian of order 6 in the neighborhood of the origin of $Y_\alpha$. We recall that
\begin{equation} \nonumber
R_6= H \circ \tau - H_2 - Z_4 - Q_4 - R_5 \circ \tau.
\end{equation}
The right-hand side of the equation is real, so $R_6$ is also real. Let us prove that $R_6$ belongs to $\mathcal{T}^{\alpha,1/2}(\mu)$, with $0 < \mu \leq \delta (m)$.

We begin by proving that $\chi_4 \in \mathcal{T}^{\alpha,1/2+}(\mu)$. We remark that, for $ i \in \mathbb{Z}$ such that $(i,j,k,l) \in \mathcal{J}$, we have
\begin{equation} \nonumber
\frac{\langle i  \rangle^{3/2}}{\vert\Omega_\iota(i,j,k,l)\vert \sqrt{\lambda_i\lambda_j\lambda_k\lambda_l}} \leq C(\mathcal{A},m), \quad \iota=0,1,2.
\end{equation}
Using this estimate and the same method as in the proof of the Lemma~\ref{XPanalytique}, we get that $\nabla_\zeta \chi_4 \in L_{\frac{1}{2}+}$.
\\It remains to prove that $\nabla^2_\zeta\chi_4 \in \mathcal{M}_{\frac{1}{2}+}$. The first terms of $\chi_4$ are indexed by $\mathcal{J}$. For $i,j \in \mathbb{Z}$, we have
\begin{equation} \nonumber
\frac{\langle i \rangle^{1/2}\langle j \rangle^{1/2} (1+ \vert \vert i \vert - \vert j \vert \vert)}{\Omega_0(i,j,k,l)\sqrt{\lambda_i\lambda_j\lambda_k\lambda_l}}  \leq 1 .
\end{equation}
The next terms of $\chi_4$ are indexed by $\mathcal{J}_2$. By Proposition~\ref{final-petit diviseur}, we have:
\begin{equation} \nonumber
\frac{\langle i \rangle^{1/2}\langle j \rangle^{1/2} (1+ \vert \vert i \vert - \vert j \vert \vert)}{\Omega_\iota(i,j,k,l)\sqrt{\lambda_i\lambda_j\lambda_k\lambda_l}}  \leq C'(m,\mathcal{A}), \quad \iota=1,2.
\end{equation}
Using these two estimates and the same method as in the proof of the Lemma~\ref{Hess poly}, we get that $\nabla^2_\zeta\chi_4 \in \mathcal{M}_{\frac{1}{2}+}$. So we proved that $\chi_4 \in \mathcal{T}^{\alpha,1/2+}(\mu)$ for $\alpha >1/2$ and $0<\mu \leq \delta(m)$.
\\By Lemma~\ref{estim crochet poisson zeta}, we have $ \left\lbrace \mathcal{T}^{\alpha,1/2}(\delta(m)),\mathcal{T}^{\alpha,1/2+}(\delta(m)) \right\rbrace \in \mathcal{T}^{\alpha,1/2}(\frac{1}{2}\delta(m)) $. So $\left\lbrace P_4, \chi_4  \right\rbrace  \in \mathcal{T}^{\alpha,1/2}(\frac{1}{2}\delta(m)) $ and for $0< \mu \leq \frac{1}{2} \delta(m) $, we have:
\begin{equation} \nonumber
\lc \left\lbrace P_4, \chi_4  \right\rbrace \rc^{\alpha,1/2} _{\mu} \leq C \mu^{-2} \lc P_4 \rc^{\alpha,1/2} _{\mu} \lc \chi_4 \rc^{\alpha,1/2+} _{\mu} \leq C\mu^6.
\end{equation}
Due to the homological equation \eqref{eq-homo}, we have:
\begin{equation} \nonumber
\lbrace H_2+P_4, \chi_4 \rbrace =Z_4+Q_4-P_4+ \lbrace P_4,\chi_4 \rbrace \in  \mathcal{T}^{\alpha,1/2}(\frac{1}{2}\delta(m)).
\end{equation}
Using Lemma~\ref{estim crochet poisson zeta} again, for $0 < \mu \leq \frac{1}{4} \delta(m)$, we have:
\begin{equation} \nonumber
\left\lbrace \left\lbrace H_2+P_4,\chi_4 \right\rbrace ,\chi_4 \right\rbrace \in \mathcal{T}^{\alpha,1/2}(\frac{1}{4}\delta(m)) \quad \text{ and } \quad\lc \left\lbrace \left\lbrace H_2+P_4,\chi_4 \right\rbrace ,\chi_4 \right\rbrace \rc^{\alpha,1/2} _{\mu}  \leq C  \mu^6.
\end{equation}
Since $\chi_4 \in \mathcal{T}^{\alpha,1/2+}(\delta(m))$ and $\lc \chi_4 \rc^{\alpha,1/2+} _{\mu} \leq C\mu^4$, we have by Corollary~\ref{estim flot zeta}
\begin{equation} \nonumber
\mathcal{T}^{\alpha,1/2}(\frac{1}{4}\delta(m)) \circ \Phi_{\chi_4}^t \in \mathcal{T}^{\alpha,1/2}(\frac{1}{8}\delta(m)).
\end{equation}
So, $R_6 \in \mathcal{T}^{\alpha,1/2}(\mu)$ and satisfies \eqref{estim-R-6} for $0<\mu \leq \frac{1}{8}\delta(m)$.

Now, consider the remainder term $R_5 \circ \tau$. Recall that $R_5=P-P_4$, so $R_5$ is real and holomorphic, of order 5 at the origin and belongs to $\mathcal{T}^{\alpha,1/2}(\delta(m))$. Using Corollary~\ref{estim flot zeta} again, we obtain that $R_5 \circ \tau \in \mathcal{T}^{\alpha,1/2}(\frac{1}{2}\delta(m))$, so $R_5 \circ \tau \in \mathcal{T}^{\alpha,1/2}(\mu)$ and satisfies for $0 < \mu \leq \frac{1}{2}\delta(m)$:
\begin{equation} \nonumber
\lc R_5\circ \tau \rc^{\alpha,1/2} _{\mu} \leq \lc P-P_4  \rc^{\alpha,1/2} _{\mu} \leq C \mu^5.
\end{equation}
We finish the proof by replacing $\frac{1}{8} \delta(m)$ by $\delta(m)$.
\endproof
\begin{lemma} \label{impage-tore-Birkhoff}
For $m\in \mathcal{U}$ and $\alpha > 1/2$, the change of variable $\tau$ defined on the Proposition~\ref{Bir.Nor.For} satisfies:
\begin{equation} \nonumber
\tau \left( \mathfrak{T}_\rho(\nu,\frac{\sigma}{2},\frac{\mu}{2},\alpha) \right)  \subset \mathfrak{T}_\rho(\nu,\sigma,\mu,\alpha),
\end{equation}
for $0<\sigma\leq 1$,  $0<\mu\leq 1$ and  $\nu \leq 4 \mu^{-2} e^{-\sigma} \delta^2(m)$.
\end{lemma}
\proof
Consider $ m\in \mathcal{U} \subset [1,2]$  and $\alpha > 1/2$, the change of variable $\tau$ satisfies:
 \begin{equation} \nonumber
\Vert \tau(\zeta)-\zeta \Vert_\alpha \leq C(m) \Vert \zeta \Vert_\alpha^3,\quad \forall \zeta \in  \mathcal{O}_{\delta(m)}(Y_\alpha) .
\end{equation}
Recall that
\begin{equation} \nonumber
\mathfrak{T}_\rho(\nu,\sigma,\mu,\alpha)=\left\lbrace (I,\theta,\xi,\eta) \mid \:  \vert I-\nu\rho \vert < \nu \mu^2, \, \vert Im \theta \vert < \sigma, \, \Vert \zeta^\mathcal{L} \Vert_\alpha < \nu^{1/2}\mu  \right\rbrace.
\end{equation}
Let $\tilde{\zeta}=\tau(\zeta)$. Then, for $\nu \leq 4 \mu^{-2} e^{-\sigma} \delta^2(m)$, we have:
\begin{equation} \label{chang_birk_id_tore}
\Vert \tilde{\zeta} - \zeta  \Vert_\alpha \leq C'(m) \mu^3 e^{\frac{3\delta}{2}} \nu^{\frac32},
\end{equation}
where $C'(m)$ is a multiple constant of $C(m)$. Using the previous estimate, let us prove that:
\begin{equation} \nonumber
\tau \left( \mathfrak{T}_\rho(\nu,\frac{\sigma}{2},\frac{\mu}{2},\alpha) \right)  \subset \mathfrak{T}_\rho(\nu,\sigma,\mu,\alpha).
\end{equation}
\begin{itemize}
\item On $\mathcal{L}$, we have:
\begin{equation} \nonumber
\Vert \tilde{\zeta} \Vert_\alpha \leq \Vert \zeta \Vert_\alpha + C'(m) \mu^3 e^{\frac{3\delta}{2}} \nu^{\frac32} < \frac12 \nu^{\frac{1}{2}} + C'(m) \mu^3 e^{\frac{3\delta}{2}} \nu^{\frac32} < \nu \mu.
\end{equation}
\item For $a \in \mathcal{A}$, we have:
\begin{align*}
\vert \tilde{I}_a - \nu \rho_a \vert & \leq \vert I_a - \nu \rho_a \vert + \vert \tilde{I}_a - I_a \vert \\
& < \frac{1}{4} \nu \mu^2 + \vert \tilde{\xi}_a \tilde{\eta}_a - \xi_a \eta_a \vert \\
& < \frac{1}{4} \nu \mu^2+ \vert \tilde{\xi}_a \vert \vert \tilde\eta_a - \eta_a\vert + \vert \eta_a \vert \vert \tilde{\xi}_a-\xi_a \vert.
\end{align*}
As $\zeta \in \mathfrak{T}_\rho(\nu,\frac{\sigma}{2},\frac{\mu}{2},\alpha)$, we have
\begin{equation} \nonumber
\vert \eta_a \vert < ( \frac{1}{2} \mu + \sqrt{2}) e^{\frac{\sigma}{2}} \nu^{\frac{1}{2}}.
\end{equation}
Using estimate \eqref{chang_birk_id_tore}, we have:
\begin{align*}
\vert \tilde{\xi}_a \vert   & < ( \frac{1}{2} \mu + \sqrt{2}) e^{\frac{\sigma}{2}} \nu^{\frac{1}{2}} +  C'(m) \mu^3 e^{\frac{3\delta}{2}} \nu^{\frac32},\\
\vert \tilde\eta_a - \eta_a\vert  + \vert \tilde{\xi}_a-\xi_a \vert  & < C'(m) \mu^3 e^{\frac{3\delta}{2}} \nu^{\frac32}.
\end{align*}
So
\begin{equation} \nonumber
\vert \tilde{I}_a - \nu \rho_a \vert < \nu \mu^2.
\end{equation}
\item It remains to verify that $\vert Im(\tilde{\theta}_a) \vert < \sigma$ for $a \in \mathcal{A}$. On the one hand we have:
\begin{equation} \nonumber
\vert \tilde{\xi}_a \vert = \left| \sqrt{\tilde{I}_a} e^{i \tilde{\theta}_a} \right| \leq \sqrt{ \tilde{I}_a}e^{\vert Im(\tilde{\theta}_a) \vert }  < \nu^{\frac{1}{2}}(\mu+\sqrt{2}) e^{\vert Im(\tilde{\theta}_a) \vert }.
\end{equation}
On the other hand, using estimate \eqref{chang_birk_id_tore}, we have:
\begin{equation} \nonumber
\vert \tilde{\xi}_a \vert    < ( \frac{1}{2} \mu + \sqrt{2}) e^{\frac{\sigma}{2}} \nu^{\frac{1}{2}} +  C'(m) \mu^3 e^{\frac{3\delta}{2}} \nu^{\frac32}.
\end{equation}
So
\begin{equation} \nonumber
e^{\vert Im(\tilde{\theta}_a) \vert } \leq \frac{\mu + 2 \sqrt{2}}{2(\mu+\sqrt{2})} e^{\frac{\sigma}{2}} + \nu \frac{C'(m)\mu^3}{\mu+\sqrt{2}} e^{\frac{3\sigma}{2}},
\end{equation}
and we get that $\vert Im(\tilde \theta) \vert \leq \sigma$.
\end{itemize}
\endproof
\subsection{Normal form on admissible sets}
We recall that $Z_4$ is given by:
\begin{equation} \nonumber
Z_4=3/4\pi \underset{(i,j,k,l) \in \mathcal{J}_2\cap\mathcal{R}_2}{\sum} \frac{\xi_i\xi_j\eta_k\eta_l}{\lambda_i\lambda_j},
\end{equation}
where
\begin{equation} \nonumber
\mathcal{J}_2= \lbrace (i,j,k,l) \in \mathcal{J}\: |\: \#\lbrace i,j,k,l \rbrace \cap \mathcal{A} \geq 2 \rbrace.
\end{equation}
We note that $Z_4$ contains integrable terms formed by the monomials of the form $\xi_i\xi_j\eta_i\eta_k=I_iI_j$. Those monomials depend only on actions defined by $I_n=\xi_n\eta_n$ for $n \in \mathbb{Z}$. We denote those terms by $Z_4^+$ and $Z_4^-=Z_4-Z_4^+$.  After straightforward computations, we obtain that:
\begin{equation} \label{Z++}
 Z_4^+=\frac{3}{4\pi} \underset{l \in \mathcal{A},\:k \in \mathbb{Z}}{ \sum} \frac{4-3\delta_{l,k}}{\lambda_l\lambda_k} I_lI_k.
\end{equation}
Concerning $Z_4^-$, we have $Z_4^-=\sum_{0\leq r \leq 4} Z_4^{-r}$ where  $r=\Card \left( \lbrace i,j,k,l \rbrace\cap\mathcal{A} \right) $. Using the definition of $\mathcal{J}_2$ yields $Z_4^{-0}=Z_4^{-1}= 0$.
\begin{lemma} \label{Z--4}
Assume that $\mathcal{A}$ is an admissible set. Then, for $m\in\mathcal{U}$, we have
\begin{equation} \nonumber
Z_4^{-4} = 0.
\end{equation}
\end{lemma}
\proof
Consider $\mathcal{F}=\left\lbrace  (i,j,k,l) \in \mathcal{J}_2 \cap \mathcal{R}_2 \cap \mathcal{A} \right\rbrace $. Then, for $(i,j,k,l) \in \mathcal{F}$ and $m\in\mathcal{U}$,
\begin{itemize}
\item[$\star$] $i+j=k+l$ using the definition of $\mathcal{J}$,
\item[$\star$] $\left\lbrace |i|,|j| \right\rbrace = \left\lbrace |k|,|l| \right\rbrace$ using Lemma~\ref{lem_ens},
\item[$\star$]  $\left\lbrace i,j \right\rbrace \neq \left\lbrace |k|,|l| \right\rbrace$ since  $\mathcal{A}$ is an admissible set.
\end{itemize}
So $\mathcal{F}=\emptyset$ and $Z_4^{-4}=0$.
\endproof
\begin{lemma}\label{Z--3}
Assume that $\mathcal{A}$ is an admissible set. Then for $m\in\mathcal{U}$
\begin{equation} \nonumber
Z_4^{-3} = 0.
\end{equation}
\end{lemma}
\proof
Suppose that there is $(i,j,k,l)\in \mathcal{J}_2 \cap \mathcal{R}_2$ such that $\#\lbrace i,j,k,l \rbrace \cap \mathcal{A}=3$. Without loss of generality, we can suppose that $i,j,k \in \mathcal{A}$ and $l\in \mathcal{L}$. Due to Lemma~\ref{lem_ens}, we have $|i|=|k|$ or $|j|=|k|$. Moreover,
 since $\mathcal{A}$ is an admissible set, we have $i=k$ or $j=k$. Suppose that $i=k$. Since $i+j=k+l$, we have $j=l$ and $l\in \mathcal{A}$, which contradicts the fact that $\mathcal{A}$ is an admissible set. So $Z_4^{-3} = 0$.
\endproof
\begin{lemma} \label{Z--2}
Assume that $\mathcal{A}$ is an admissible set. Then, for any $m\in\mathcal{U}$,
\begin{equation} \nonumber
Z_4^{-2} =  0.
\end{equation}
\end{lemma}
\proof
Consider $\aleph = \lbrace (i,j,k,l)\in \mathcal{J}_2 \cap \mathcal{R}_2 \;|\; \#\lbrace i,j,k,l \rbrace \cap \mathcal{A}=2 \rbrace$.
\begin{itemize}
\item[$\star$] Assume that $i,j \in \mathcal{A}$ and $k,l \in \mathcal{L}$. Then, by Lemma~\ref{lem_ens}, we have $i=-k$ and $j=-l$ or $i=-l$ and $j=-k$. Without loss of generality, we can suppose that $i=-k$ and $j=-l$. We have $(i,j,k,l) \in \mathcal{J}$,  so $i+j=k+l=-k-l$ and $i=-j$. It contradicts the fact that $\mathcal{A}$ is an admissible set. The case where $k,l \in \mathcal{A}$ and $i,j \in \mathcal{L}$ is treated in the same way.
\item[$\star$] Assume now that $i,l \in \mathcal{A}$ and $j,k \in \mathcal{L}$.  By Lemma~\ref{lem_ens}, we have $|i|=|k|$ and $|j|=|l|$ or $i=l$ and $|j|=|k|$. Let us first consider the case where $i=l$ and $|j|=|k|$. since $i+j=k+l$, we have $j=k$. So the monomial $\xi_i\xi_j\eta_j\eta_i$ will be in $Z_4^-$, which leads to a contradiction. Consider now the case where $|i|=|k|$ and $|j|=|l|$. Since $\mathcal{A}$ is an admissible set, we have $i=-k$ and $l=-j $. As $i+j=k+l$, necessarily  $i=l$ and $j=k$, which lead to the previous case. The cases where $k,l \in \mathcal{A}$ and $i,j \in \mathcal{L}$ or $i,k \in \mathcal{A}$ and $j,l \in \mathcal{L}$ are treated in the same way.
\end{itemize}
Thus, we deduce that $\aleph=\emptyset$ and $Z_4^{-2} = 0$
\endproof
\subsection{Action-angle variables}
As in \eqref{tore-complex}, we pass from variables $(\xi,\eta)$ to $(I,\theta,\zeta^{\mathcal{L}})$, where $I=(I_a, a \in \mathcal{A})$, $\theta=(\theta_a, \in \mathcal{A})$ and $\zeta^\mathcal{L}=(\zeta_s, \in \mathcal{L})$. We recall that, for $a \in \mathcal{A}$, action-angle variables $I$ and $ \theta$ are given by:
\begin{equation} \nonumber
\xi_a=\sqrt{I_a} e^{i\theta_a}, \quad \eta_a = \sqrt{I_a} e^{-i\theta_a}.
\end{equation}
In these new variables, the symplectic form $-id\xi \wedge d\eta$ becomes:
\begin{equation} \label{symplect-action-angle}
-\sum_{a\in \mathcal{A}} dI_a \wedge d\theta_a - i  \sum_{s\in \mathcal{L}} d\xi_s \wedge d\eta_s.
\end{equation}
Moreover $I$ is of order 2, $\theta$ is of order zero, $\xi$ and $\eta$ are of order 1.

Using expression \eqref{Z++} $Z_4^+$ and  the lemmas  \ref{Z--4}-\ref{Z--2}, the Hamiltonian \ref{Birkhoff-FN} becomes:
\begin{align*}
H\circ\tau &= \underset{a \in \mathcal{A}}{ \sum} \omega_aI_a + \frac{3}{4\pi} \underset{l,a \in \mathcal{A}}{ \sum} \frac{4-3\delta_{l,a}}{\lambda_l\lambda_a} I_lI_a +  \underset{s \in \mathcal{L}}{ \sum} \lambda_s\xi_s\eta_s + \frac{3}{\pi} \underset{l \in \mathcal{A},\:s \in \mathcal{L}}{ \sum} \frac{1}{\lambda_l\lambda_s} I_l\xi_s\eta_s
\\ &+ Q_4 + R_5 \circ \tau + R_6.
\end{align*}
The first line contains the integrable terms. The second one contains:
\begin{itemize}
\item $Q_4$,  of order 4 and at least of order 3 in the direction of $\mathcal{L}$;
\item $R_6$, which comes from the Birkhoff normal form and of order 6;
\item $R_5 \circ \tau$, which comes from term of order 5 of the nonlinearity  \eqref{non-linearite}.
\end{itemize}

The Hamiltonian $H \circ \tau$ depends on variables $\left( I,\theta,\zeta^{\mathcal{L}} \right) $. For the rest of the paper, we will drop the multi-index $\mathcal{L}$ ( i.e. $ \zeta^{\mathcal{L}} $ will be replaced by $\zeta$).
\subsection{Rescaling the variables} \label{Re-scaling des variables}
We want to study the Hamiltonian $H_1$. To do this, we will rescale the variables $(I,\theta,\xi,\eta)$ by considering the following change of variables:
\begin{equation} \label{res-calling}
\chi_\rho: \left( \tilde{r},\tilde{\theta},\tilde{\xi},\tilde{\eta}\right)  \mapsto \left( I,\theta,\xi,\eta \right)  ,
\end{equation}
where
\begin{align*}
 I=&\nu\rho+ \nu \tilde{r}, \quad \theta=\tilde{ \theta},\\
 \xi=&\nu^{1/2} \tilde \xi,\quad \quad  \eta=\nu^{1/2} \tilde{\eta}.
\end{align*}
We have
\begin{equation} \nonumber
\chi_\rho: \mathcal{O}^\alpha\left( \frac{\sigma}{2}, \frac{\mu}{2}  \right) \rightarrow \mathfrak{T}_\rho \left( \nu,\frac{\sigma}{2}, \frac{\mu}{2},\alpha \right).
\end{equation}
In these new variables, the symplectic form \eqref{symplect-action-angle} becomes
\begin{equation}
-\nu \sum_{a\in \mathcal{A}} d\tilde{r}_a \wedge d \tilde{\theta}_a - i \nu  \sum_{s\in \mathcal{L}} d\tilde\xi_s \wedge d\tilde\eta_s.
\end{equation}
Consider
\begin{equation} \nonumber
\breve{\Psi}=\breve{\Psi}_\rho=\tau \circ \chi_\rho.
\end{equation}
The change of variables $\chi_\rho$ is linear on $\rho$. We can extend $\Phi$ holomorphically on
\begin{equation} \nonumber
\mathcal{D}_c = \left\lbrace \rho \in\mathbb{C}^{\mathcal{A}}  \mid  \left| \rho_j - \frac{3}{2} \right| \leq c , \; 1 \leq j \leq \operatorname{Card}( \mathcal{A})   \right\rbrace .
\end{equation}
 To simplify notations, we will drop the tilde. In these rescaled variables, the Hamiltonian $H$ becomes up to a constant:
 \begin{equation*}
\begin{split}
H \circ \breve{\Psi}& = \nu \sum_{\substack{a \in \mathcal{A}}} \omega_a r_a +  \nu^2 \frac{3}{2\pi}    \sum_{\substack{a,l \in \mathcal{A}}}  \frac{4-3\delta_{a,l}}{\lambda_a\lambda_l} \rho_l r_a \\
&+  \nu \underset{s \in \mathcal{L}}{ \sum} \lambda_s \xi_s \eta_s + \nu^2 \frac{3}{\pi}  \sum_{\substack{l \in \mathcal{A} \\ s\in \mathcal{L}}} \frac{1}{\lambda_l\lambda_s} \rho_l \xi_s \eta_s\\
&+ \nu^2 \frac{3}{4\pi} \sum_{\substack{a,l \in \mathcal{A}}}  \frac{4-3\delta_{a,l}}{\lambda_a\lambda_l} r_l r_a + \nu^2 \frac{3}{\pi}  \sum_{\substack{l \in \mathcal{A} \\ s\in \mathcal{L}}} \frac{1}{\lambda_l\lambda_s} r_l \xi_s \eta_s \\
&+\left( Q_4 + R_5 \circ \tau  + R_6 \right)  \circ \chi_\rho.
\end{split}
\end{equation*}
By dividing by $\nu$, we can rewrite the previous Hamiltonian under the following form:
\begin{equation} \label{new-hamiltonian}
\nu^{-1}H \circ \breve{\Psi}=h_0 + f,
\end{equation}
where $h_0 \equiv h_0(r, \xi , \eta;\rho, \nu)$ and contains linear terms in $r$, quadratic terms in $\xi,\eta$ and independent from the angle variable $\theta$. The new perturbation $f$ contains all the rest and depends on the angle variable. More precisely
\begin{equation} \label{new-perturbation}
f \equiv \nu \frac{3}{4\pi} \sum_{\substack{a,l \in \mathcal{A}}}  \frac{4-3\delta_{a,l}}{\lambda_a\lambda_l} r_l r_a + \nu \frac{3}{\pi}  \sum_{\substack{l \in \mathcal{A} \\ s\in \mathcal{L}}} \frac{1}{\lambda_l\lambda_s} r_l \xi_s \eta_s + \nu^{-1} \left( Q_4 + R_5 \circ \tau  + R_6 \right)  \circ \chi_\rho.
\end{equation}
We can rewrite the new Hamiltonian $h$ under the following form
\begin{equation} \label{ham}
h_0=\Omega\cdot r + \sum_{a\in\mathcal{L}}\Lambda_a \xi_a\eta_a,
\end{equation}
where $\Omega=(\Omega_k)_{k\in\mathcal{A}}$, and
\begin{align}\label{Om}\Omega_k = \Omega_k (\rho,\nu) & = \omega_k + \nu \tilde{\omega}_k = \omega_k + \nu \frac{3}{2\pi} \frac{1}{\lambda_k} \sum_{l \in \mathcal{A}} \frac{4-3\delta_{l,k}}{\lambda_l}\rho_l,\\
\label{Lam}\Lambda_a = \Lambda_a(\rho,\nu) & = \lambda_a + \nu \tilde{\lambda}_a = \lambda_a + \nu \frac{3}{\pi} \frac{1}{\lambda_a}\sum_{l\in\mathcal{A}}\frac{\rho_l}{\lambda_l},
\end{align}
for $\rho \in \mathcal{D}_c$. We remark that for the internal frequencies, we have:
\begin{equation} \label{matrice frequence}
\tilde{\omega}_k=\sum_{l \in \mathcal{A}} M_k^l \rho_l, \quad M_k^l= \frac{3}{2\pi}\frac{4-3\delta_{l,k}}{\lambda_ k\lambda_l}.
\end{equation}
$M$ is an invertible matrix, since
\begin{equation} \nonumber
\det M = \left( \frac{3}{2\pi}\right)^{n} \left( \prod_{l \in \mathcal{A}} \lambda_l^{-2} \right) \left( 4n-3\right) \left( -3\right)^{n-1}, \quad n=\operatorname{Card}( \mathcal{A}).
\end{equation}
Recall that $\rho \in \left[ 1,2 \right] ^{\mathcal{A}}$ and $\mathcal{A} \subset \lbrace a \in \mathbb{Z}^d \mid |a| \leq N \rbrace$. For $\rho \in \mathcal{D}_c$, we have
\begin{equation} \nonumber
\left| \Omega(\rho) -\omega(\rho) \right| \leq \nu \left| \frac{3}{2\pi} \sum_{l \in \mathcal{A}} \frac{4-3\delta_{l,k}}{\lambda_l}\rho_l \right|, \quad \text{and} \quad \vert \Lambda_a(\rho) - \lambda_a (\rho)\vert \leq \nu \vert a \vert^{-1} \left| \frac{3}{\pi} \sum_{l\in\mathcal{A}}\frac{\rho_l}{\lambda_l} \right|.
\end{equation}
This proves estimations \eqref{new_frequencies} and concludes the proof of the first and second points of Theorem~\ref{theo-F.N}.

We recall that $\zeta=(\zeta_a)_{a\in\mathcal{L}}$. The quadratic part of the Hamiltonian $h$ is given by the following infinite matrix:
\begin{equation} \nonumber
A(\rho,\nu)=\diag \left( \left(\begin{array}{ll}
0& \Lambda_a(\rho,\nu)\\
\Lambda_a(\rho,\nu)& 0
\end{array}\right), \ a\in\mathcal{L}\right) .
\end{equation}
We can put the Hamiltonian $h$ under the following form
\begin{equation} \nonumber
h=\Omega (\rho,\nu) \cdot r + \frac{1}{2} \langle A(\rho,\nu)\zeta,\zeta \rangle.
\end{equation}
The Hamiltonian operator is given by:
\begin{equation} \nonumber
iJA(\rho,\mu)=\diag \Big(\left(\begin{array}{rr}
-i \Lambda_a(\rho,\nu) &0\\
0 &i \Lambda_a(\rho,\nu)
\end{array}\right), \ a\in \mathcal{L}\Big).
\end{equation}
The spectrum of the Hamiltonian operator is:
\begin{equation} \nonumber
\sigma(iJA)=\{\pm i \Lambda_a(\rho,\nu),\ a\in\mathcal{L}\}.
\end{equation}

Let us study now the perturbation \eqref{new-perturbation}. Due to Proposition~\ref{Bir.Nor.For}, the perturbation $f$  is real holomorphic and belongs to $\mathcal{T}^{\alpha,1/2}(\mathcal{D},\frac{\sigma}{2},\frac{\mu}{2})$. Using estimations \eqref{estim-Z-Q}-\eqref{estim-R-5} and for $x=(r,\theta,\zeta) \in \mathcal{O}^\alpha\left( \frac{\sigma}{2},\frac{\mu}{2} \right)$, we have:
\begin{align*}
\vert f \vert & \leq C\nu \mu^4, \\
\Vert \nabla_\zeta f \Vert_\alpha & \leq C \nu \mu^3,\\
\vert \nabla_\zeta f \vert_{1/2} & \leq C \nu \mu^3,\\
\vert \nabla^2_\zeta f \vert_{1/2} & \leq C \nu \mu^2,
\end{align*}
where $C$ is a constant that depends on the admissible set $\mathcal{A}$, the mass $m$ and the nonlinearity $g$. So we have:
\begin{equation} \nonumber
\lc f \rc^{\alpha,1/2} _{\frac\sigma2,\frac\mu2,\mathcal{D}} \leq C \nu\mu^4.
\end{equation}
Let us study now the jet of the perturbation $f$. Recall that the jet function $f^T$ is defined by:
\begin{equation} \nonumber
f^T  =  f(\theta,0,0,\rho)+\nabla_r f(\theta,0,0,\rho)r +\langle\nabla_\zeta f(\theta,0,0,\rho),\zeta\rangle +\frac{1}{2} \langle \nabla^2_{\zeta\zeta}f(\theta,0,0,\rho)\zeta,\zeta\rangle.
\end{equation}
We look for the terms of \eqref{new-perturbation} which can contribute to $f^T$. Clearly the first two terms do not contribute to $f^T$. The third term is indexed by $\mathcal{J}_2^c$ and  does not contribute to $f^T$. Let us now look at $R_5 \circ \tau$. According to \ref{Bir.Nor.For}, $R_5\circ \tau$ is of order 5. Moreover $R_5\circ \tau$ depends on the action variable (of order 2), on the angle $\theta$, on $\xi$ and $\eta$ (of order 1).
So $R_5\circ \tau$ can contains terms like:
\begin{itemize}
\item[-] $I^{5/2}$ will contribute to $f(\theta,0,0,\rho)$ and $\nabla_r f(\theta,0,0,\rho)$;
\item[-] $I^2 \xi $, $I^2 \eta$ will contribute to $\nabla_\xi f(\theta,0,0,\rho)$ or $\nabla_\eta f(\theta,0,0,\rho)$;
\item[-] $I^{3/2} \xi \eta$, $I^{3/2} \xi \xi$, $I^{3/2} \eta \eta$ will contribute to $\nabla^2_{\xi\eta}f(\theta,0,0,\rho)$,$\nabla^2_{\xi\xi}f(\theta,0,0,\rho)$ or $\nabla^2_{\eta\eta}f(\theta,0,0,\rho)$.
\end{itemize}
So $R_5\circ \tau \circ \chi_\rho$ will contribute to $f^T$. A similar result holds for $R_6\circ \chi_\rho$. We deduce that:
\begin{equation} \nonumber
\lc f^T \rc^{\alpha,1/2} _{\frac\sigma2,\frac\mu2,\mathcal{D}} \leq C \nu^{3/2}\mu^5.
\end{equation}
To finish the proof of the third point of Theorem~\ref{theo-F.N}, we have to look at the derivative in $\rho$ of the perturbation $f$ and its jet $f^T$. Note that, according to \eqref{new-perturbation}, the dependence of the perturbation $f$ on $\rho$ comes from the change of variable $ \chi_\rho $ via the relation $ I = \nu \rho + \nu \tilde r$. We can extend $f$ holomorphically on $\mathcal{D}_c$ with the same estimates. Using a Cauchy estimate on $\mathcal{D}_c$, we obtain
\begin{equation} \nonumber
\lc \partial_\rho f \rc^{\alpha,1/2} _{\frac\sigma2,\frac\mu2,\mathcal{D}} \leq c^{-1} \lc f \rc^{\alpha,1/2} _{\frac\sigma2,\frac\mu2,\mathcal{D}}, \quad \lc \partial_\rho f^T \rc^{\alpha,1/2} _{\frac\sigma2,\frac\mu2,\mathcal{D}}\leq c^{-1} \lc  f^T \rc^{\alpha,1/2} _{\frac\sigma2,\frac\mu2,\mathcal{D}}.
\end{equation}
Then $\partial_\rho f$ and $\partial_\rho f^T$ satisfy the same estimates as $f$ and $f^T$.
\subsection{Real variables}
In the normal form \eqref{Final FN}, the quadratic part is expressed in complex variables. However, the KAM theorem is expressed in real variables. In order to remedy to this, we consider the following symplectic change of variable:
\begin{equation} \nonumber
\Upsilon(r,\theta, \zeta) = (r, \theta, \breve{\zeta}),
\end{equation}
where $\breve{\zeta}=(p,q)=\left( \breve{\zeta}_s=(p_s,q_s), \: s \in \mathcal{L} \right)$. Variables $p$ and $q$ are given by:
\begin{equation} \nonumber
\xi_s=\frac{1}{\sqrt{2}}(p_s+iq_s), \quad \eta_s=\frac{1}{\sqrt{2}}(p_s-iq_s).
\end{equation}
Under the hypotheses of Theorem~\ref{theo-F.N}, consider $\Psi_\rho= \breve{\Psi}_\rho \circ \Upsilon$. This change of variables is real holomorphic and satisfies:
\begin{equation} \nonumber
\Psi_\rho: \mathcal{O}^\alpha \left(   \frac{\sigma}{2} , \frac{\mu}{2} \right)  \rightarrow \mathfrak{T}_\rho(\nu,\sigma,\mu,\alpha).
\end{equation}
It transforms the symplectic form $ -id\xi \wedge d\eta$ into $-dr \wedge d\theta - dp \wedge dq$. In the new variables, the normal form \eqref{Final FN} becomes:
\begin{equation} \label{real final NF}
\begin{array}{rcl}
h(r,\theta,\breve{\zeta};\rho)& = & H\circ \Psi_\rho \\
& = & \displaystyle \Omega(\rho)\cdot r + \frac{1}{2}  \sum_{a\in\mathcal{L}}\Lambda_a ( \rho) (p_a^2+q_a^2) + f(r, \theta, \breve{\zeta}; \rho)\\
& = & \displaystyle \Omega(\rho)\cdot r + \frac{1}{2} \sum_{a\in\mathcal{L}}\langle \breve{A}(\rho,\nu)\breve{\zeta},\breve{\zeta} \rangle + f(r, \theta, \breve{\zeta}; \rho),\\
\end{array}
\end{equation}
where
\begin{equation*}
\breve{A}(\rho,\nu)=\diag \left( \left(\begin{array}{ll}
\Lambda_a(\rho,\nu)& 0 \\
0& \Lambda_a(\rho,\nu)
\end{array}\right), \ a\in\mathcal{L}\right) .
\end{equation*}
In these new variables, the perturbation $f$ satisfies the estimates of the third point of Theorem~\ref{theo-F.N}.
\section{KAM for the wave equation}
\subsection{Abstract KAM theorem} In this section we state a KAM theorem and adapt the notations for the cubic nonlinear wave equation on the circle.

Consider a real Hamiltonian $h_\rho$ on normal form and depends on a parameter $\rho$, given by:
\begin{equation}\label{hom}
h(r,p,q;\rho)=\Omega(\rho) \cdot r+ \frac{1}{2} \sum_{a\in\mathcal{L}}\Lambda_a (\rho)\left(  p_a^2+q^2_a\right) ,
\end{equation}
with
\begin{itemize}
\item $\rho \in \mathcal{D}$, a compact set of $\mathbb{R}^p$;
\item $\Omega:  \mathcal{D}\to\mathbb{R}^n$ a $C^1$ internal frequency vector;
\item $\mathcal{L}$ a set of $\mathbb{Z}$;
\item for $a\in\mathcal{L}$, $\Lambda_a$ an external  frequency of class $C^1$ on $\mathcal{D}$.
\end{itemize}
The internal frequencies $\Omega$ and the external frequencies $\Lambda$ satisfy certain hypotheses which will be stated in the following paragraph. Let us fix two parameters $0 < \delta_0 \leq \delta \leq 1 $ and consider  $\mathcal{A}^-$ a finite set of $\mathcal{L}$. The set $\mathcal{L} \setminus \mathcal{A}^-$ shall be denoted by $\mathcal{L}^\infty$.

\textbf{Hypothesis A1: Separation condition.} Assume that, for all $\rho \in \mathcal{D}$, we have:
\begin{itemize}
\item[$\star$] for all $a \in \mathcal{L}$,
\begin{equation}
 \Lambda_{a} (\rho)  \geq c_0 \langle a \rangle;
 \end{equation}
\item[$\star$] for all $a$, $b$ $\in \mathcal{L}$ and $|a| \neq |b|$, we have
\begin{equation}
 \vert \Lambda_{a} (\rho) - \Lambda_{b} ( \rho) \vert \geq c_1\left| \left| a \right| - \left| b \right| \right|.
\end{equation}
\end{itemize}
\textbf{Hypothesis A2: Transversality condition.} Assume that for all $\Omega' \in \mathcal{C}^1 (\mathcal{D},\mathbb{R}^n)$ that satisfies
\begin{equation} \nonumber
\vert \Omega - \Omega' \vert_{\mathcal{C}^1 \left( \mathcal{D} \right)  } < \delta_0,
\end{equation}
\emph{for all $k \in \mathbb{Z}^n $, there exists a unit vector $z_k \in \mathbb{R}^p$, and all $a$, $b$ $\in \mathcal{L}$ with $|a|>|b|$ the following holds}:
\begin{itemize}
\item[$\star$]
\begin{equation} \nonumber
\vert k\cdot\Omega' (\rho) \vert  \geq \delta,\quad \forall \rho \in \mathcal{D},
\end{equation}
or
\begin{equation} \nonumber
\langle \partial_\rho ( k\cdot\Omega'(\rho)),z_k \rangle \geq \delta \quad \forall \rho \in \mathcal{D};
\end{equation}
where $k \neq 0$
\item[$\star$]
\begin{equation} \nonumber
\vert k\cdot\Omega' (\rho)   \pm \Lambda_{a}(\rho) \vert  \geq \delta \langle a \rangle, \quad \forall \rho \in \mathcal{D},
\end{equation}
or
\begin{equation} \nonumber
\langle \partial_\rho ( k\cdot\Omega'(\rho)\pm \Lambda_{a}(\rho)),z_k \rangle \geq \delta \quad \forall \rho \in \mathcal{D};
\end{equation}
\item[$\star$]
\begin{equation} \nonumber
\vert k\cdot\Omega' (\rho)  + \Lambda_{a}(\rho) + \Lambda_{b} ( \rho)  \vert \geq \delta( \langle a \rangle + \langle b \rangle), \quad \forall \rho \in \mathcal{D},
\end{equation}
or
\begin{equation} \nonumber
\langle \partial_\rho ( k\cdot\Omega'(\rho)+\Lambda_{a}(\rho)+\Lambda_{b}(\rho)),z_k \rangle \geq \delta \quad \forall \rho \in \mathcal{D};
\end{equation}
\item[$\star$]
\begin{equation} \nonumber
\vert k\cdot\Omega' (\rho) + \Lambda_{a}(\rho) - \Lambda_{b} ( \rho) \vert \geq \delta ( 1 + \vert  \vert a \vert - \vert b \vert \vert ), \quad \forall \rho \in \mathcal{D},
\end{equation}
or
\begin{equation} \nonumber
\langle \partial_\rho ( k\cdot\Omega'(\rho)+\Lambda_{a}(\rho) - \Lambda_{b} ( \rho)),z_k \rangle \geq \delta \quad \forall \rho \in \mathcal{D};
\end{equation}
\end{itemize}

\textbf{Hypothesis A3: Second Melnikov condition.}
Assume that for all $\Omega' \in \mathcal{C}^1 (\mathcal{D},\mathbb{R}^n)$ that satisfies
\begin{equation} \nonumber
\vert \Omega - \Omega' \vert_{\mathcal{C}^1 \left( \mathcal{D} \right)  } < \delta_0,
\end{equation}
the following holds:
\\ for each $0 < \kappa < \delta$ and $N>1$ there exists a closed set $\mathcal{D}'\subset \mathcal{D}$ that satisfies
\begin{equation}
mes ( \mathcal{D} \setminus \mathcal{D}') \leq C (\delta^{-1} \kappa)^\tau N^\iota;
\end{equation}
for some $\tau,\: \iota >0$, such that for all  $\rho \in \mathcal{D}'$, all $ 0< \vert k \vert < N$ and all $a$,$b\in \mathcal{L}$ with $|a| \neq |b|$ we have:
\begin{equation}
\vert \Omega'( \rho) \cdot k+\Lambda_{a}(\rho)-\Lambda_{b}(\rho) \vert \geq \kappa(1+  \vert \vert a \vert - \vert b \vert \vert).
\end{equation}
\smallskip

We denote $A_0= \diag\lbrace \Lambda_a I_2, a \in \mathcal{L} \rbrace$. Now we are able to state the abstract KAM theorem proved in \cite{B1}:

\begin{theorem}\label{theoreme kam}
Assume that $h$ is a Hamiltonian given by \eqref{hom} and satisfies hypotheses A1, A2 and A3 for fixed $\delta$ and $\delta_0$ and all $\rho \in \mathcal{D}$. Fix $\alpha,\beta>0$ and $0<\sigma,\mu \leq1$. Then there is $\varepsilon_0$ depending on $\alpha,\beta, \sigma ,n,\mu, |\omega_0|_{\mathcal{C}^1 \left( \mathcal{D} \right) }$ and $\vert  A_0 \vert_{\beta,{\mathcal{C}^1 \left( \mathcal{D} \right) }}$ such that, if  $\partial^j_\rho f\in \mathcal{T}^{\alpha,\beta}(\mathcal{D},\sigma, \mu)$ for $j=0,1$,  if $$\lc f^T \rc^{\alpha,\beta,\kappa} _{\sigma,\mu,\mathcal{D}}=\varepsilon < \min(\varepsilon_0, \frac{1}{8} \delta_0) \mbox{ and  } \: \lc f \rc^{\alpha,\beta,\kappa} _{\sigma,\mu,\mathcal{D}} 
=O(\varepsilon^\tau),$$
for $\: 0 < \tau <1$, then there is a Borel set $\mathcal{D}' \subset \mathcal{D}$ with $mes ( \mathcal{D} \setminus \mathcal{D}') \leq  c \varepsilon^\gamma$ such that for all $\rho \in \mathcal{D}'$:
\begin{itemize}
\item there is a symplectic analytical change of variable
\begin{equation} \nonumber
\Phi=\Phi_\rho: \mathcal{O}^\alpha(\frac{\sigma}{2}, \frac{\mu}{2}) \rightarrow \mathcal{O}^\alpha(\sigma,\mu)
\end{equation}
\item there is a new internal frequency vector $\tilde{\Omega}(\rho) \in \mathbb{R}^n$, a matrix $A \in \mathcal{M}_\beta$ and a perturbation $\tilde{f} \in \mathcal{T}^{\alpha,\beta}(\mathcal{D}',\sigma/2, \mu/2)$ such that
\begin{equation}  \nonumber
(h_\rho + f ) \circ \Phi = \tilde{\Omega}(\rho) \cdot r + \frac{1}{2}  \langle \zeta,A(\rho)\zeta \rangle + \tilde{f} (\theta,r,\zeta; \rho),
\end{equation}
where $ A:\mathcal{L}\times\mathcal{L} \to \mathcal{M}_{2\times 2}(\mathbb{R})$ is a block diagonal symmetric infinite matrix in $\mathcal{M}_\beta$ (ie $A_{[a]}^{[b]}=0$ if $[a] \neq [b]$). Moreover $\partial_r \tilde{f}=\partial_\zeta \tilde{f}=\partial_{\zeta\zeta}^2 \tilde{f}=0$ for $r=\zeta=0$. The change of variables  $\Phi = (\Phi_\theta, \Phi_r,\Phi_\zeta)$ is close to identity, and for all $x \in  \mathcal{O}^\alpha(\frac{\sigma}{2}, \frac{\mu}{2})$ and all $\rho \in \mathcal{D}'$, we have:
\begin{equation} \label{kam_id}
 \Vert \Phi - Id \Vert_\alpha \leq C \varepsilon^{4/5}.
\end{equation}
For all $\rho \in \mathcal{D}'$, the new frequencies $\tilde{\omega}$ and matrix $A$ satisfy
\begin{equation} \label{freq_kam}
\left|  A(\rho)-A_0(\rho))\right|_\alpha\leq C \varepsilon, \qquad \vert \tilde{\Omega} (\rho) - \Omega (\rho) \vert_{\mathcal{C}^1(\mathcal{D}')} \leq C \varepsilon,
\end{equation}
where $C$ is a constant that depends on $\varepsilon_0$.
\end{itemize}
\end{theorem}
\subsection{Verification of the hypotheses of the KAM theorem}
\subsubsection{Non resonance}
In this section, we verify that the real normal form \eqref{real final NF} satisfies the hypotheses of Theorem~\ref{theoreme kam}. We start by verifying the separation hypothesis $A1$, then the transversality condition $A2$ and finally the second Melnikov condition $A3$.
\begin{lemma} \label{separation freq}
For all $\rho \in \mathcal{D}$, and all $a$ , $b \in \mathcal{L}$, we have
\begin{align*}
(i) & \:\:\:  \Lambda_{a} (\rho)   \geq  \langle a \rangle ;  \\
(ii) & \:\:\:  \vert \Lambda_{a} (\rho) - \Lambda_{b} ( \rho) \vert \geq \frac{1}{8} \left| \left| a \right| - \left| b \right| \right|,  \quad \text{with} \: |a| \neq |b|.
\end{align*}
\end{lemma}
\proof
Recall that, for $a \in \mathcal{L}$, the external frequencies are given by:
\begin{equation} \nonumber
\Lambda_a(\rho)  = \lambda_a + \nu \frac{3}{\pi} \frac{1}{\lambda_a}\sum_{l\in\mathcal{A}}\frac{\rho_l}{\lambda_l}  = \lambda_a + \nu  \frac{C}{\lambda_a}.
\end{equation}
Estimation $(i)$ is obvious. For $(ii)$, we remark that, for $\nu$ small enough, we have:
\begin{equation}\nonumber
 C\nu \left| \frac{1}{\lambda_a}- \frac{1}{\lambda_b}  \right| = \frac{C\nu}{\lambda_a \lambda_b} \left| \lambda_a - \lambda_b \right|  \leq \frac{1}{2} \left| \lambda_a - \lambda_b \right|.
\end{equation}
So
\begin{equation} \nonumber
\vert \Lambda_a - \Lambda_b \vert \geq \frac{1}{2} \vert \lambda_a - \lambda_b \vert \geq \frac{1}{8} \vert \vert a \vert - \vert b \vert \vert,
\end{equation}
which concludes the proof of the lemma.
\endproof

The non resonance hypothesis $A2$ will be verified in three steps. We begin by recalling the results obtained in the Propositions~\ref{final-petit diviseur} and \ref{prop-D3}. For $ \kappa = \nu^{1/2} $, we have the following lemma:
\begin{lemma}
For $\gamma >0$ small enough, $|k| \leq \nu^{-\gamma}$ and  $(a,b) \in \mathcal{L}^2$ we have:
\begin{itemize}
\item[] \begin{equation} \nonumber
|\omega \cdot k|\geq 2 \nu^{1/2},
\end{equation}
except when $k$ is  $D_0$  resonant.
\item[] \begin{equation} \nonumber
|\omega \cdot k+\lambda_a|\geq 2 \nu^{1/2}\langle a \rangle,
\end{equation}
except when $(k,a)$ is $D_1$ resonant.
\item[] \begin{equation} \nonumber
|\omega\cdot k+\lambda_a+\lambda_b|\geq 2 \nu^{1/2}(\langle a \rangle+\langle b \rangle),
\end{equation}
except when $ (k,a,b)$ when $D_2$ resonant.
\item[] \begin{equation} \nonumber
|\omega\cdot k+\lambda_a-\lambda_b|\geq 2 \nu^{1/2}(1+ \vert \vert a \vert - \vert b \vert \vert ),
\end{equation}
where $ \vert a \vert \neq \vert b \vert$  and $(k,a,b)$ is not $ D_3$  resonant.
\end{itemize}
\end{lemma}
\begin{remark}
In the previous lemma, we have applied the Propositions~\ref{final-petit diviseur} and \ref{prop-D3} with $\kappa=\nu^{1/2}$, $N= \nu^{-\gamma}$ and $m \in [1,2] \setminus \mathcal{C}$. The Lebesgue measure of $\mathcal{C}$ satisfies:
\begin{equation} \nonumber
\operatorname{mes} \left(  \mathcal{C} \right) \leq C \kappa^{\tau}N^{\iota},
\end{equation}
where $\tau=O(\dfrac{1}{n})$ and $\iota=O(n^2)$. With this choice of parameter, the Lebesgue measure of $\mathcal{C}$ still small, if we assume that $\gamma < O(\dfrac{1}{n^3})$.
\end{remark}
Now we will verify the transversality hypothesis $A2$ for $k$ small. Recall that the internal frequencies are given by:
\begin{equation} \nonumber
 \Omega= \omega + \nu M \rho,
\end{equation}
where $M$ is the symmetric invertible matrix defined in \eqref{matrice frequence}. We denotes $C_\mathcal{A} = \Vert  M^{-1} \Vert_2$.
\begin{lemma} \label{Cond non resonance petit}
For $\gamma >0$ small enough, $k \in \mathbb{Z}^n$ with  $|k| \leq \nu^{-\gamma}$,  $(a,b) \in \mathcal{L}^2$, consider $\delta_0 = \frac{1}{2} C_{\mathcal{A}}^{-1} \nu$. Then for all $\Omega' \in \mathcal{C}^1(\mathcal{D}, \mathbb{R}^n)$ that satisfies
\begin{equation} \nonumber
\vert \Omega - \Omega' \vert_{\mathcal{C}^1 \left( \mathcal{D} \right)  } < \delta_0,
\end{equation}
and all $\rho \in \mathcal{D}$, we have
\begin{itemize}
\item[] \begin{equation} \nonumber
|\Omega' \cdot k|\geq \nu^{1/2},
\end{equation}
except when $k$ is  $D_0$  resonant.
\item[] \begin{equation} \nonumber
|\Omega' \cdot k+\Lambda_a|\geq  \nu^{2/3}\langle a \rangle,
\end{equation}
except when $(k,a)$ is $D_1$ resonant.
\item[] \begin{equation} \nonumber
|\Omega'\cdot k+\Lambda_a+\Lambda_b|\geq  \nu^{2/3}(\langle a \rangle+\langle b \rangle),
\end{equation}
except when $ (k,a,b)$ is $D_2$ resonant.
\item[] \begin{equation} \nonumber
|\Omega'\cdot k+\Lambda_a-\Lambda_b|\geq  \nu^{2/3}(1+\vert \vert a \vert - \vert b \vert \vert),
\end{equation}
whith $ \vert a \vert \neq \vert b \vert$  and $(k,a,b)$ is not $ D_3$  resonant.
\end{itemize}
\end{lemma}
\proof
Let $k \in \mathbb{Z}^n$ such that $\vert k \vert \leq \nu^{-\gamma}$ for $\gamma>0$ small enough. We begin with the first estimate. We have
\begin{equation} \nonumber
\vert \Omega'-\omega \vert \leq \vert \Omega' - \Omega \vert + \vert \Omega - \omega \vert \leq \frac{1}{2} C_\mathcal{A}^{-1} \nu + c \nu \leq C \nu.
\end{equation}
So, for $a \in \mathcal{A}$
\begin{equation} \nonumber
\vert \Omega'-\omega \vert \leq C \nu\langle a \rangle.
\end{equation}
By the Cauchy-Schwarz inequality, for all $(a,b) \in \mathcal{L}$, we have:
\begin{align*}
& \vert \Omega'\cdot k-\omega\cdot k \vert \leq C \nu^{1-\gamma} \leq \nu^{1/2}.\\
& \vert \Omega'\cdot k-\omega\cdot k \vert \leq \nu^{1/2} \langle a \rangle.\\
& \vert \Omega'\cdot k-\omega\cdot k \vert  \leq \nu^{1/2}(\langle a \rangle+\langle b \rangle).  \\
& \vert \Omega'\cdot k-\omega\cdot k \vert \leq \nu^{1/2}(1+\vert \vert a \vert - \vert b \vert \vert), \quad \vert a \vert \neq \vert b \vert.
\end{align*}
To conclude the proof of the first case, we use the fact that:
\begin{equation*}
\vert\Omega' \cdot k \vert\geq \vert \omega \cdot k\vert- \vert \Omega'\cdot k-\omega\cdot k \vert \geq 2 \nu^{1/2} - \nu^{1/2}=\nu^{1/2}.
\end{equation*}
Let us now look at the second estimate. We note that for $a \in \mathcal{L}$, we have:
\begin{equation} \nonumber
|\lambda_a - \Lambda_a| \leq \tilde{c} \nu \langle a \rangle.
\end{equation}
So
\begin{align*}
|\Omega' \cdot k+\Lambda_a| & \geq |\Omega' \cdot k+\lambda_a| - |\lambda_a - \Lambda_a| \\
& \geq |\omega \cdot k+\lambda_a| - \vert \Omega'\cdot k-\omega\cdot k \vert - |\lambda_a - \Lambda_a| \\
& \geq 2 \nu^{1/2} \langle a \rangle - \nu^{1/2} \langle a \rangle - \tilde{c} \nu \langle a \rangle \geq \nu^{2/3} \langle a \rangle.
\end{align*}
Consider now the third estimate. For $(a,b) \in \mathcal{L}^2$, we have:
\begin{align*}
\vert \Omega'\cdot k+\Lambda_a+\Lambda_b \vert & \geq \vert \Omega' \cdot k+\lambda_a + \lambda_b\vert  - \vert \lambda_a - \Lambda_b\vert - \vert \lambda_b - \Lambda_a\vert \\
& \geq |\omega \cdot k+\lambda_a+\lambda_b| - \vert \Omega'\cdot k-\omega\cdot k \vert - |\lambda_a - \Lambda_a| - |\lambda_b - \Lambda_b| \\
& \geq (  2 \nu^{1/2} - \nu^{1/2}- 2 \tilde{c} \nu )  (\langle a \rangle + \langle b \rangle) \geq \nu^{2/3} (\langle a \rangle+ \langle b \rangle).
\end{align*}
Let us now look at the last small divisor. Using \eqref{Lam}, we remark that, for $(a,b) \in \mathcal{L}^2$ with $ \vert a \vert \neq \vert b \vert$, we have:
\begin{align*}
\vert \Lambda_a - \Lambda_b - ( \lambda_a - \lambda_b) \vert & \leq \tilde{c}\nu \left|  \lambda_a^{-1} -  \lambda_b^{-1} \right| \\
& =\tilde{c}\nu \frac{\vert  \vert a \vert + \vert b \vert \vert }{(\lambda_a+\lambda_b)\lambda_a\lambda_b} \vert \vert a \vert - \vert b \vert \vert \\
& \leq \tilde{c}\nu (1+\vert \vert a \vert - \vert b \vert \vert).
\end{align*}
Which leads to
\begin{align*}
\vert \Omega'\cdot k+\Lambda_a-\Lambda_b \vert & \geq \vert \Omega' \cdot k+\lambda_a - \lambda_b\vert  - \vert \Lambda_a - \Lambda_b - ( \lambda_a - \lambda_b) \vert \\
& \geq |\omega \cdot k+\lambda_a-\lambda_b| - \vert \Omega'\cdot k-\omega\cdot k \vert - \vert \Lambda_a - \Lambda_b - ( \lambda_a - \lambda_b) \vert \\
& \geq (  2 \nu^{1/2} - \nu^{1/2}-  \tilde{c} \nu )  (1+\vert \vert a \vert - \vert b \vert \vert) \geq \nu^{2/3} (1+\vert \vert a \vert - \vert b \vert \vert).
\end{align*}
The proof is thus achieved.
\endproof
We have verified the non resonance hypotheses for $\vert k \vert \leq \nu^{-\gamma}$ , $\delta_0 = \frac{1}{2} C_{\mathcal{A}}^{-1} \nu$ and $\delta = \nu^{2/3}.$ For large $k$, i.e. $\vert k \vert > \nu^{-\gamma}$, we verify the separation conditions $A2$ on the derivatives in $\rho $ of the small divisors. More precisely, we have:
\begin{lemma} \label{Cond non resonance grand}
For $\gamma >0$ small enough, $k \in \mathbb{Z}^n$ with  $|k| > \nu^{-\gamma}$,  $(a,b) \in \mathcal{L}^2$, we consider $\delta_0 = \frac{1}{2} C_{\mathcal{A}}^{-1} \nu$. Then for all $\Omega' \in \mathcal{C}^1(\mathcal{D}, \mathbb{R}^n)$ that satisfies
\begin{equation} \nonumber
\vert \Omega - \Omega' \vert_{\mathcal{C}^1 \left( \mathcal{D} \right)  } < \delta_0,
\end{equation}
there exists a unit vector $z_k$, such that for all  $\rho \in \mathcal{D}$ we have
\begin{align*}
&(i)  \:\:\: \vert \langle \partial_\rho ( k\cdot\Omega'(\rho)),z_k \rangle \vert \geq C \nu^{1-\gamma} \geq \nu ,\\
&(ii) \:\:\vert \langle \partial_\rho ( k\cdot\Omega'(\rho)\pm \Lambda_{a}(\rho)),z_k \rangle \vert \geq \nu \quad \mbox{ for all } a \in \mathcal{L},\\
&(iii)\:\vert \langle \partial_\rho ( k\cdot\Omega'(\rho)+\Lambda_{a}(\rho)+\Lambda_{b}(\rho)),z_k \rangle \vert \geq \nu \quad \mbox{ for all } (a,b) \in \mathcal{L}^2,\\
&(iv) \:\: \vert \langle \partial_\rho ( k\cdot\Omega'(\rho)+\Lambda_{a}(\rho) - \Lambda_{b} ( \rho)),z_k \rangle \vert \geq \nu \quad \mbox{ for all } (a,b) \in \mathcal{L}^2.
\end{align*}
The constant $C$ depends on the admissible set $\mathcal{A}$.
\end{lemma}
\proof
Let us begin with $(i)$. First of all, we remark that:
\begin{equation} \nonumber
\langle \partial_\rho ( k.\Omega'(\rho)),z_k \rangle = \langle \partial_\rho ( k.\Omega(\rho)),z_k \rangle + \langle \partial_\rho ( k.\Omega'(\rho))-k.\Omega(\rho)),z_k \rangle.
\end{equation}
However
\begin{equation} \nonumber
\vert  \langle \partial_\rho ( k.\Omega'(\rho))-k.\Omega(\rho)),z_k \rangle \vert \leq  \frac{1}{2}  \vert k \vert C_{\mathcal{A}}^{-1} \nu.
\end{equation}
Recall that the matrix $M$ is symmetric and $\vert k \vert > \nu^{-\gamma}$. Assume that $z_k=\frac{Mk}{\left| Mk \right|}$. Then, we have:
\begin{align*}
\vert \langle \partial_\rho ( k\cdot\Omega'(\rho)),z_k \rangle \vert & \geq \nu \left| \langle Mk , \frac{Mk}{\left| Mk \right|}\rangle \right| - \left| \langle \partial_\rho ( k.\Omega'(\rho))-k.\Omega(\rho)),z_k \rangle \right|\\
&  = \nu  \left| Mk \right| - \frac{1}{2} \vert k \vert C_{\mathcal{A}}^{-1} \nu \\
& \geq \nu \left( C_{\mathcal{A}}^{-1} - \frac{1}{2} C_{\mathcal{A}}^{-1} \right) \left| k \right| \\
& \geq \frac{1}{2} C_{\mathcal{A}}^{-1} \nu^{1-\gamma} \geq \nu .
\end{align*}
Let us now consider the second estimation $(ii)$. For $a$, $l \in \mathcal{L}$, we have:
\begin{equation} \nonumber
\left| \partial_{\rho_l} \Lambda_a ( \rho) \right| = \left| \nu \frac{3}{\pi} \frac{1}{\lambda_a \lambda_l} \right| \leq c \nu.
\end{equation}
Using the same unit vector $z_k$, we get:
\begin{align*}
\left|  \langle \partial_\rho ( k\cdot\Omega'(\rho)\pm \Lambda_{a}(\rho)),z_k \rangle   \right| & \geq \left| \langle \partial_\rho ( k\cdot\Omega'(\rho)),z_k \rangle \right| - \left|  \partial_\rho ( \Lambda_a(\rho)) \right| \left| z_k \right| \\
& \geq \frac{1}{2} C_{\mathcal{A}}^{-1} \nu^{1-\gamma} - c' \nu \geq \nu.\end{align*}
Applying the same principle for $(iii)$ and $(iv)$, we obtain:
\begin{align*}
\left| \langle \partial_\rho ( k\cdot\Omega'(\rho)+\Lambda_{a}(\rho)\pm \Lambda_{b}(\rho)),z_k \rangle \right| & \geq   \left| \langle \partial_\rho ( k\cdot\Omega'(\rho)),z_k \rangle \right| - \left|  \partial_\rho ( \Lambda_a(\rho)) \right| \left| z_k \right| \\
 & -  \left|  \partial_\rho ( \Lambda_b(\rho)) \right| \left| z_k \right| \\
& \geq  \frac{1}{2} C_{\mathcal{A}}^{-1} \nu^{1-\gamma} - 2c' \nu \geq \nu,
\end{align*}
and the proof is thus concluded.
\endproof
To finish the verification of the transversality  condition, it remains to consider the cases where $(k,a)$ is $D_1$ resonant and $(k,a,b)$ is $D_2$ or $D_3$ resonant.
\begin{lemma}
Let $k \in \mathbb{Z}^n$ and  $(a,b) \in \mathcal{L}^2$. Consider $\delta_0 = \frac{1}{4} \nu \breve C_{\mathcal{A}}$, where $ \breve C_{\mathcal{A}}$ is a constant that depends on the admissible set $\mathcal{A}$. Then for all $\Omega' \in \mathcal{C}^1(\mathcal{D}, \mathbb{R}^n)$ that satisfies
\begin{equation} \nonumber
\vert \Omega - \Omega' \vert_{\mathcal{C}^1 \left( \mathcal{D} \right)  } < \delta_0,
\end{equation}
and all $\rho \in \mathcal{D}$, we have:
\begin{itemize}
\item[] \begin{equation} \nonumber
|\Omega' \cdot k+\Lambda_a|\geq  \breve C_{\mathcal{A}} \nu \langle a \rangle,
\end{equation}
if $(k,a)$ is $D_1$ resonant;
\item[] \begin{equation} \nonumber
|\Omega'\cdot k+\Lambda_a+\Lambda_b|\geq  \breve C_{\mathcal{A}} \nu (\langle a \rangle+\langle b \rangle),
\end{equation}
if $ (k,a,b)$ is $D_2$ resonant;
\item[] \begin{equation} \nonumber
|\Omega'\cdot k+\Lambda_a-\Lambda_b|\geq   \breve C_{\mathcal{A}} \nu (1+\vert \vert a \vert - \vert b \vert \vert),
\end{equation}
where $ \vert a \vert \neq \vert b \vert$  and $(k,a,b)$ is $ D_3$  resonant.
\end{itemize}
\end{lemma}
\proof
Consider $a \in \mathcal{A}$. From \eqref{Om} and \eqref{Lam}, we have
\begin{equation} \nonumber
\tilde{\omega}_a - \tilde{\lambda}_a = \frac{3}{2\pi} \frac{1}{\lambda_a} \sum_{l \in \mathcal{A}} \frac{2-3\delta_{l,a}}{\lambda_l}\rho_l.
\end{equation}
Assume that $(k,a)$ is $D_1$ resonant. Then we have
\begin{equation} \nonumber
\Omega' \cdot k+\Lambda_a = (\Omega_a' - \Omega_a) - \nu (\tilde{\omega}_a - \tilde{\lambda}_a).
\end{equation}
We remark that there exists a constant $C_1$, that depends on the admissible set $\mathcal{A}$, such that:
\begin{equation} \nonumber
\vert \tilde{\omega}_a - \tilde{\lambda}_a \vert \geq \frac{C_1}{\lambda_a}.
\end{equation}
Consider $ \delta_0 \leq \dfrac{1}{2} \nu \dfrac{C_1}{\max(\mathcal{A})^2+2}:=\dfrac{1}{2} \nu \tilde{C}_1$, then we have:
\begin{align*}
\vert \Omega' \cdot k+\Lambda_a \vert & \geq C_1 \nu \frac{1}{\lambda_a} - \frac{1}{2} \tilde C_1 \nu \\
& \geq \tilde{C}_1  \nu \langle a \rangle  - \frac{1}{2} \tilde C_1 \nu  \\
& \geq \frac{1}{2} \tilde{C}_1  \nu \langle a \rangle.
\end{align*}
Let us consider now the case where $(k,a,b)$ is $D_2$ resonant. There is a constant $C_2$, that depends on the admissible set  $ \mathcal{A}$, such that for all $(a,b) \in \mathcal{A}^2$,  we have:
\begin{align*}
\vert (\tilde{\omega}_a - \tilde{\lambda}_a ) + (\tilde{\omega}_b - \tilde{\lambda}_b ) \vert & \geq C_2 ( \frac{1}{\lambda_a} + \frac{1}{\lambda_b}) \\
& \geq \tilde{C}_2 (\langle a \rangle + \langle b \rangle),
\end{align*}
where $  \tilde{C}_2:= \dfrac{C_2}{\max(\mathcal{A})^2+2}$. So, if $(k,a,b)$ is $D_2$ resonant and $\delta_0 \leq \dfrac{1}{4} \nu \tilde{C}_2$, then we have:
\begin{align*}
|\Omega'\cdot k+\Lambda_a+\Lambda_b| & = \vert (\Omega' - \Omega) \cdot k - \nu ((\tilde{\omega}_a - \tilde{\lambda}_a ) + (\tilde{\omega}_b - \tilde{\lambda}_b )) \vert \\
& \geq  \nu \vert (\tilde{\omega}_a - \tilde{\lambda}_a ) + (\tilde{\omega}_b - \tilde{\lambda}_b ) \vert - 2 \vert \Omega' - \Omega \vert \\
& \geq \frac{1}{2} \tilde{C}_2 \nu (\langle a \rangle + \langle b \rangle).
\end{align*}
It remains to look at the last small divisor in the case where $(k,a,b)$ is $ D_3$ resonant. We note that there exists a constant $C_3$, which depends on the admissible set $ \mathcal{A}$, such that for all $(a,b) \in \mathcal{A}^2$, we have:
\begin{align*}
\vert (\tilde{\omega}_a - \tilde{\lambda}_a ) - (\tilde{\omega}_b - \tilde{\lambda}_b ) \vert & \geq C_3 \left| \frac{1}{\lambda_a} - \frac{1}{\lambda_b} \right| \\
& \geq \tilde{C}_3 (1+ \vert \vert a \vert - \vert b \vert \vert),
\end{align*}
where $ \tilde{C}_3:= \dfrac{C_3}{8(\max(\mathcal{A})^2+2)}$. So, if $(k,a,b)$ is $D_3$ resonant and $\delta_0 \leq \dfrac{1}{4} \nu \tilde{C}_3$, then we have:
\begin{align*}
|\Omega'\cdot k+\Lambda_a-\Lambda_b| & = \vert k \cdot (\Omega' - \Omega) - \nu ((\tilde{\omega}_a - \tilde{\lambda}_a ) - (\tilde{\omega}_b - \tilde{\lambda}_b )) \vert \\
& \geq  \nu \vert (\tilde{\omega}_a - \tilde{\lambda}_a ) - (\tilde{\omega}_b - \tilde{\lambda}_b ) \vert - 2 \vert \Omega' - \Omega \vert \\
& \geq \frac{1}{2} \tilde{C}_3 \nu (1+ \vert \vert a \vert - \vert b \vert \vert).
\end{align*}
We conclude the proof by choosing $\breve C_{\mathcal{A}} = \frac{1}{2} \min(\tilde{C}_1,\tilde{C}_2,\tilde{C}_3)$.
\endproof

The last hypothesis to verify to apply the KAM theorem is the second Melnikov condition.
Recall that $n=\operatorname{Card}( \mathcal{A})$.
\begin{lemma} \label{2end Melnikov for wave}
For $\delta = \nu$ , $\delta_0 \leq \delta$, $\tau=\frac{1}{3}$ and $\iota = n+\frac{3}{2}+\frac{2}{3\gamma}$, the second Melnikov condition is satisfied.
\end{lemma}
\proof
Consider $\gamma >0 $ small enough and $N \geq 0$. If $N \leq \nu^{-\gamma}$, then, using Lemma~\ref{Cond non resonance petit}, the second Melnikov condition is satisfied for all $\rho \in \mathcal{D}$. If $N > \nu^{-\gamma}$, then from Lemma~\ref{Cond non resonance petit} and for $|k| \leq \nu^{-\gamma}$ the second Melnikov condition is satisfied for all $\rho \in \mathcal{D}$.
Assume now that $\vert k \vert > \nu^{-\gamma}$, then from Lemma~\ref{Cond non resonance grand} there is a unit vector $z_k$ such that:
\begin{equation} \nonumber
\vert \langle \partial_\rho ( k\cdot\Omega'(\rho)+\Lambda_{a}(\rho) - \Lambda_{b} ( \rho)),z_k \rangle \vert \geq \nu.
\end{equation}
For $0< \kappa < \nu$, consider the set
\begin{equation} \nonumber
J(k,a,b) = \lbrace \rho  \in\mathcal{D} \mid \left| \Omega'(\rho)\cdot k+\Lambda_a(\rho)-\Lambda_b(\rho)  \right| < \kappa \rbrace.
\end{equation}
Then we have
\begin{equation} \nonumber
\operatorname{mes} J(k,a,b) \leq C \kappa \nu^{-1},
\end{equation}
where $C$ is a constant that depends on $\mathcal{D}$.
For $p \in \mathbb{Z}$ and $k \in \mathcal{Z}^n$, we consider the following set
\begin{equation} \nonumber
W(k,p) = \lbrace \rho \in \mathcal{D} \mid \left| \Omega' \cdot k + p \right| < 5 \kappa^{1/3} \rbrace.
\end{equation}
By the first estimate $(i)$ from Lemma~\ref{Cond non resonance grand}, we have
\begin{equation} \nonumber
\operatorname{mes} W(k,p) \leq C \kappa^{1/3} \nu^{-1},
\end{equation}
where $C$ is a constant that depends on $\mathcal{D}$. Consider
\begin{equation} \nonumber
W = \lbrace \rho \in \mathcal{D} \mid \left| \Omega' \cdot k + p \right| < 5 \kappa^{1/3} \rbrace.
\end{equation}
We remark that $\left| \Omega' \cdot k + p \right| < 5 \kappa^{1/3}$ for $\vert k \vert < N$, so we have $|p|\leq C|k|< CN$. This leads to
\begin{equation} \nonumber
\operatorname{mes} \left(  W \right) \leq \sum \limits_{\underset{\vert k \vert \leq  N}{k \in \mathbb{Z}^n}} \sum \limits_{\underset{\vert p \vert < CN}{p \in \mathbb{Z}}}  W(k,p) \leq C N^{n+1} \kappa^{1/3} \nu^{-1}.
\end{equation}
For $a \in  \mathcal{L}$, we recall that $\Lambda_a(\rho)=\lambda_a+ \nu C(\rho)\lambda_a^{-1}$. Consider $\nu$ small enough, then for $\rho \in \mathcal{D}$, we have:
\begin{equation} \nonumber
\left|  \Lambda_a( \rho) - \vert a \vert  \right| = \left|  \lambda_a( \rho) - \vert a \vert + \nu \frac{C(\rho)}{\lambda_a} \right| = \left|  \frac{m}{\lambda_a + \vert a \vert} + \nu \frac{C(\rho)}{\lambda_a} \right| \leq  \frac{2}{\vert a \vert}.
\end{equation}
If $ \vert a \vert > \vert b \vert > k^{-1/3}$, we obtain
\begin{equation} \nonumber
\left| \Lambda_a(\rho) - \Lambda_b(\rho) - \left( \vert a  \vert - \vert b \vert  \right)    \right| \leq \frac{4}{\vert b \vert} \leq 4 \kappa^{1/3}.
\end{equation}
So, if $\rho \in \mathcal{D} \setminus W$ and $\vert a  \vert > \vert b \vert > \kappa^{-1/3}$, we obtain
\begin{align*}
\left| k\cdot\Omega'(\rho)+\Lambda_{a}(\rho) - \Lambda_{b} ( \rho) \right| & \geq \left| k\cdot\Omega'(\rho)+ \left(  \vert a \vert - \vert b \vert  \right)  \right|  - \left| \Lambda_a(\rho) - \Lambda_b(\rho) - \left( \vert a  \vert - \vert b \vert  \right) \right| \\
 & \geq 5 \kappa^{1/3} - 4 \kappa^{1/3} = \kappa^{1/3}.
\end{align*}
It remains to look at the cases where $\vert a \vert \leq \kappa^{-1/3}$ or $\vert b \vert \leq \kappa^{-1/3}$.
Then, there is $ k \in \mathbb{Z}^n$ such that
\begin{equation} \nonumber
\left| k\cdot\Omega'(\rho)+\Lambda_{a}(\rho) - \Lambda_{b} ( \rho) \right| < 1,
\end{equation}
for $\nu^{-\gamma} < \vert k \vert < N$. We remark that, in those cases, $\left| \vert a \vert - \vert b \vert \right| \leq CN$.  Consider the set:
\begin{equation} \nonumber
\mathcal{Q} = \left\lbrace (a,b) \in \mathbb{Z}^2 \mid \min( \vert a \vert, \vert b \vert)  \leq \kappa^{-1/3} \:\text{et} \: \left| \vert a \vert - \vert b \vert \right| \leq CN   \right\rbrace.
\end{equation}
We have \begin{equation} \nonumber
\operatorname{Card} \left(  \mathcal{Q} \right)  \leq C N \kappa^{-2/3}.
\end{equation}
Let
\begin{equation} \nonumber
\mathcal{D}' = \mathcal{D} \setminus (  W \bigcup_{\stackrel{\vert k \vert \leq  N}{(a,b) \in \mathcal{Q}}} \left( J (k,a,b)\right)) .
\end{equation}
Then for all $\rho \in \mathcal{D}'$ we have:
\begin{equation} \nonumber
\left| k\cdot\Omega'(\rho)+\Lambda_{a}(\rho) - \Lambda_{b} ( \rho) \right| \geq \kappa.
\end{equation}
Moreover,
\begin{align*}
\operatorname{mes} \left(  \mathcal{D} \setminus \mathcal{D}' \right) & \leq \operatorname{mes}  \left(  W \right)  +  \sum \limits_{\underset{\vert k \vert \leq  N}{k \in \mathbb{Z}^n}} \sum_{(a,b) \in \mathcal{Q}} \operatorname{mes} J(k,a,b) \\
& \leq C N^{n+1} \kappa^{1/3} \nu^{-1} + CN^{n} N \kappa^{-2 /3} \kappa \nu^{-1} \\
& \leq CN^{n+1} (\kappa \nu^{-1})^{1/3}.
\end{align*}
Recall that $N> \nu^{-\gamma}$. This leads to
\begin{equation} \nonumber
\operatorname{mes} \left(  \mathcal{D} \setminus \mathcal{D}' \right) \leq CN^{n+1+2/3\gamma} \nu^{-2/3} (\kappa \nu^{-1})^{1/3}.
\end{equation}
Now it remains to show that for all $ \rho \in \mathcal{D}'$ we have:
\begin{equation} \nonumber
\vert \Omega'( \rho) \cdot k+\Lambda_{a}(\rho)-\Lambda_{b}(\rho) \vert \geq \kappa(1+ \vert a \vert - \vert b \vert \vert).
\end{equation}
We will prove this estimation in two steps. Assume at first that $\vert  \vert a \vert - \vert b \vert \vert \geq 16 \vert \Omega' \cdot k \vert $. Then, using the second separation condition from Lemma~\ref{separation freq}, we have:
\begin{align*}
\vert \Omega'( \rho) \cdot k+\Lambda_{a}(\rho)-\Lambda_{b}(\rho) \vert & \geq \vert \Lambda_{a}  - \Lambda_{b}  \vert -  \vert \Omega' \cdot k \vert\\
& \geq \frac{1}{16} \vert  \vert a \vert - \vert b \vert \vert \geq \frac{1}{32}  ( 1 + \vert  \vert a \vert - \vert b \vert \vert) \\
& \geq \kappa ( 1 + \vert  \vert a \vert - \vert b \vert \vert),
\end{align*}
for all  $ 0 \leq \kappa \leq \frac{1}{32}$ and all $\rho \in \mathcal{D}$.

Assume now that $\vert  \vert a \vert - \vert b \vert \vert < 16 \vert \Omega' \cdot k \vert < 16 C_{\mathcal{A}}N$. Then for all $\rho \in \mathcal{D}'$, $ 0< \vert k \vert < N$ and $a$,$b\in \mathcal{L}$ with $|a| \neq |b|$, we have:
\begin{align*}
\vert \Omega'\cdot k +\Lambda_a-\Lambda_b  \vert &  \geq \frac{\kappa}{1+16C_\mathcal{A}N} (1+ \vert a \vert - \vert b \vert ) \\
& \geq \tilde{\kappa} (1+ \vert a \vert - \vert b \vert ),
\end{align*}
where
\begin{equation} \nonumber
\operatorname{mes} \left(  \mathcal{D} \setminus \mathcal{D}' \right)  \leq  C N^{n+3/2+2/3\gamma}\left(\tilde\kappa \nu^{-7/5}   \right)^{1/2},
\end{equation}
which ends the proof and concludes the verification of the assumptions made on the frequencies.
\endproof

\subsubsection{Application of KAM theorem~\ref{theoreme kam}}
Using lemmas \ref{separation freq}-\ref{2end Melnikov for wave}, the separation conditions, the Transversality condition and the second Melnikov condition are satisfied for $\delta = \nu$ and $\delta_0 \leq \delta$.
From Theorem~\ref{theo-F.N}, for $0< \mu \leq 1$, we have:
\begin{equation} \label{taille-perturbation-onde}
\lc f^T \rc^{\alpha,1/2,\kappa} _{\frac\sigma2,\frac\mu2,\mathcal{D}}  \leq 2 C  \nu^{3/2}.
\end{equation}
Moreover we have:
\begin{equation} \nonumber
\lc f \rc^{\alpha,\gamma,\kappa} _{\sigma,\mu,\mathcal{D}} = O \left(  \left( \lc f^T \rc^{\alpha,1/2,\kappa} _{\frac\sigma2,\frac\mu2,\mathcal{D}}  \right)^{3/2} \right)  .
\end{equation}
To use Theorem~\ref{theoreme kam}, we need that:
\begin{equation} \nonumber
\lc f \rc^{\alpha,\gamma,\kappa} _{\sigma,\mu,\mathcal{D}} \ll  \delta= \nu^{7/5}
\end{equation}
and
\begin{equation} \nonumber
\lc f \rc^{\alpha,\gamma,\kappa} _{\sigma,\mu,\mathcal{D}} = O \left(  \left( \lc f^T \rc^{\alpha,1/2,\kappa} _{\frac\sigma2,\frac\mu2,\mathcal{D}}  \right)^{\tau} \right),  \mbox{ for } \tau \in \left] 0,1\right[ .
\end{equation}
These last two conditions are indeed verified. We can therefore apply the Theorem~\ref{theoreme kam}.
\section{Proof of Theorem~\ref{théorème onde}}
Now we have all the tools to prove Theorem~\ref{théorème onde}. For $m \in \mathcal{U}$ and $\rho \in \left( \mathcal{D} \setminus \mathcal{D}' \right) $,  consider:
\begin{equation} \nonumber
\Pi_\rho=\Psi_\rho \circ \Phi_\rho = \tau \circ \chi_\rho \circ \Upsilon \circ \Phi_\rho,
\end{equation}
where
\begin{itemize}
\item $\tau$ is the Birkhoff change of variable constructed in the Proposition~\ref{Bir.Nor.For}.
\item $\chi_\rho$ is the rescaling defined in \ref{res-calling};
\item $\Upsilon$ is the transition to real variables;
\item $\Phi_\rho$ is the KAM change of variable from Theorem~\ref{theoreme kam}.
\end{itemize}
So $\Pi_\rho$ is a real holomorphic symplectic  change variable
\begin{equation} \nonumber
\Pi_\rho: \mathcal{O}^\alpha \left(   \frac{\sigma}{4} , \frac{\mu}{4} \right)  \rightarrow \mathfrak{T}_\rho(\nu,\sigma,\mu,\alpha),
\end{equation}
that transforms the perturbed Hamiltonian \eqref{hamiltonien} into
\begin{equation} \nonumber
H \circ \Pi_\rho=\omega'(\rho)\cdot r + \frac{1}{2}  \langle \zeta,A(\rho)\zeta \rangle + \tilde{f} (\theta,r,\zeta; \rho).
\end{equation}
where $ A:\mathcal{L}\times\mathcal{L} \to \mathcal{M}_{2\times 2}(\mathbb{R})$ is a block diagonal symmetric infinite matrix in $\mathcal{M}_\beta$ (ie $A_{[a]}^{[b]}=0$ if $[a] \neq [b]$). Moreover $\partial_r \tilde{f}=\partial_\zeta \tilde{f}=\partial_{\zeta\zeta}^2 \tilde{f}=0$ for $r=\zeta=0$. From \eqref{matrice frequence} and \eqref{freq_kam}, the internal frequencies are given by
\begin{equation} \nonumber
\omega' = \omega + \nu M \rho  + O(\nu^{3/2}).
\end{equation}
For the following, let $I=\nu \rho$ and $\mathcal{D} = [ \nu, 2 \nu ]^n$. Then, for $m \in \mathcal{U}$, there is a Borel set $\mathcal{D}' \subset [ \nu, 2 \nu ]^n$ such that:
\begin{equation} \nonumber
\operatorname{mes}\left( \left[ \nu , 2 \nu \right] ^n \setminus  \mathcal{D}'  \right) \leq \nu^{\gamma+n},
\end{equation}
for $\gamma>0$ and depends on $n$. For $X=(\theta, r , \zeta)$, we denotes $[X]_\alpha = \left(  \underset{a \in \mathcal{A}}{ \sum} \vert r_a e^{2i\theta_a} \vert \right) ^{1/2} + \Vert \zeta \Vert_\alpha$. Let $X \in \mathbb{T}^n \times \lbrace I \rbrace \times \lbrace 0 \rbrace$. So $\chi_\rho \circ \Upsilon(X) \in \mathfrak{T}_\rho(\nu,[X]_\alpha,\alpha)$, and we have:
\begin{equation} \nonumber
\operatorname{dist}_\alpha ( \chi_\rho \circ \Upsilon(X),X) \leq 2 \nu^{1/2} [X]_\alpha \leq 4 \nu^{3/2}.
\end{equation}
Using \eqref{chang.birk.id}, we obtain that
\begin{equation}  \nonumber
\operatorname{dist}_\alpha ( \tau \circ \chi_\rho \circ \Upsilon(X), X) \leq  \nu^{3/2} + \operatorname{dist}_\alpha ( \chi_\rho \circ \Upsilon(X),X) \leq 5 \nu^{3/2}.
\end{equation}
Then, thanks to \eqref{kam_id} and \eqref{taille-perturbation-onde}, we have
\begin{equation} \nonumber
\operatorname{dist}_\alpha ( \Pi_\rho(X), X) \leq C \nu^{4/5},
\end{equation}
where $C$ is an absolute constant. For $m \in \left( \left[ 1,2 \right] \setminus \mathcal{U} \right)$ and $I \in \mathcal{D}'$, consider $(\tilde{\theta},\tilde{I},\tilde{\zeta})=\Pi^{-1}(X)$;
 and let
\begin{equation} \nonumber
u\left( \tilde \theta,\tilde I,x\right) =\sum_{a \in \mathcal{A}}\sqrt{\tilde I_a} \frac{e^{-i\tilde \theta_a}\varphi_a(x)+e^{i\theta_a}\varphi_{-a}(x)}{\sqrt{2}\omega'^{1/4}}.
\end{equation}
The function
$
t\mapsto u(\tilde \theta+t\omega', \tilde I,x)
$
is a quasi-periodic solution of the wave equation \eqref{eq onde}. Let $\zeta_{I,\theta}=(\xi_{I,\theta},\eta_{I,\theta})$ where
\begin{align*} \nonumber
& \left( \xi_{I,\theta} \right)_a =  \sqrt{I_a} e^{i\theta_a}, \quad \left( \eta_{I,\theta} \right)_a =  \sqrt{I_a} e^{-i\theta_a}, \quad \mbox{if } a \in \mathcal{A}, \\
& \left( \xi_{I,\theta} \right)_s = \left( \eta_{I,\theta} \right)_s= 0 , \quad \mbox{if } s \in \mathcal{L}.
\end{align*}
Then we have:
\begin{equation} \nonumber
 \sup_{\theta_0 \in \mathbb{T}^n} \parallel u(\theta_0,I,.)-u_{I,m}(\theta_0,.) \parallel_{H^\alpha} \leq  \Vert \Pi_\rho(X)-\zeta_{I,\theta}\Vert_\alpha \leq C \nu^{4/5}.
\end{equation}
\bibliographystyle{plain}
\nocite{*}
\bibliography{biblio1}
\end{document}